\documentclass[final,3p,11pt]{my-article} 

\usepackage{graphicx}
\usepackage{subfigure}
\usepackage[colorlinks=true,citecolor=blue]{hyperref}
\usepackage{mathtools}

\usepackage{amssymb}

\usepackage{amsmath}

\usepackage{lineno}
\usepackage{marginnote}
\usepackage{color,soul}
\usepackage{setspace}
\usepackage{multirow}
\usepackage{url}
\usepackage{colortbl} 
\usepackage{bm}	
\usepackage{enumerate}
\usepackage{isomath} 

\newcommand{\eref}[1]{(\ref{#1})}
\newcommand{\fref}[1]{Fig.~\ref{#1}}

\newcommand{\vm}[1]{\mathbf{#1}}
\newcommand{\bsym}[1]{\bm{#1}}

\renewcommand{\Re}{{\rm{I\!R}}}
\newcommand{\diffx}{d\vm{x}}
\newcommand{\diffs}{ds}
\newcommand{\transpose}{\mathsf{T}}

\newcommand{\cons}{\mathrm{c}}
\newcommand{\stab}{\mathrm{s}}
\newcommand{\trace}[1]{\mathrm{trace}{\,\left(#1\right)}}
\newcommand{\smat}[2][ccccccccccccccccccccccccccccccccccccccccccccccccccc]{\left
[\arraycolsep=4.0pt\def\arraystretch{0.9}\begin{array}{#1}#2 \\ \end{array} \right]}

\newcommand{\alejandro}[1]{{\color{black}{#1}}}
\newcommand{\alejandrog}[1]{{\color{black}{#1}}}

\biboptions{sort&compress}

\usepackage{siunitx}
\usepackage{textcomp}

\usepackage{xcolor}

\begin{document}

\begin{frontmatter}


\title{
A node-based uniform strain virtual element method for compressible and nearly incompressible elasticity
}

\author[adr1]{A. Ortiz-Bernardin\corref{cor1}}
\ead{aortizb@uchile.cl}
\author[adr1]{R. Silva-Valenzuela}
\author[adr2]{S. Salinas-Fern\'andez}
\author[adr2]{N.~Hitschfeld-Kahler}
\author[adr1]{S. Luza}
\author[adr1]{B. Rebolledo}

\cortext[cor1]{Corresponding author. Tel: +56 2 2978 4664, Fax: +56 2 2689 6057,}
\address[adr1]{Computational and Applied Mechanics Laboratory, Department of Mechanical Engineering, Universidad de Chile, Av. Beauchef 851, Santiago 8370456, Chile.}
\address[adr2]{Department of Computer Science, Universidad de Chile, Av. Beauchef 851, Santiago 8370456, Chile.}

\doublespacing
\begin{abstract} 

We propose a combined nodal integration and virtual element method for compressible
and nearly incompressible elasticity, wherein the strain is averaged at the nodes from 
the strain of surrounding virtual elements. For the strain averaging procedure, 
a nodal averaging operator is constructed using a generalization to virtual elements 
of the node-based uniform strain approach for finite elements. We refer to the proposed
technique as the node-based uniform strain virtual element method (NVEM). No additional degrees of
freedom are introduced in this approach, thus resulting in a displacement-based formulation.
A salient feature of the NVEM is that the stresses and strains become nodal variables 
just like displacements, which can be exploited in nonlinear simulations. Through several 
benchmark problems in compressible and nearly incompressible elasticity
\alejandro{as well as in elastodynamics}, we demonstrate 
that the NVEM is accurate, optimally convergent and devoid of volumetric locking.

\end{abstract}


\end{frontmatter}

\doublespacing

\section{Introduction}
\label{sec:intro}
\alejandro{The virtual element method~\cite{BeiraoDaVeiga-Brezzi-Cangiani-Manzini-Marini-Russo:2013} 
(VEM) is a generalization of the finite element method (FEM) for approximation of field variables 
on very general meshes formed by elements with arbitrary number of edges/faces 
(convex or nonconvex polytopes) known as virtual elements.
In brief, the method consists in the construction of an algebraic
(exact) representation of the stiffness matrix without computation of basis
functions (basis functions are \textit{virtual}). In this process, a decomposition
of the stiffness matrix into a consistency part and a stability part that
ensures convergence of the method~\cite{Cangiani-Manzini-Russo-Sukumar:2015}
is realized. Although the VEM is relatively modern technology, 
nowadays its development includes applications such as 
elastic and inelastic solids~\cite{BeiraodaVeiga-Brezzi-Marini:2013,
BeiraodaVeiga-Lovadina-Mora:2015,Artioli-BeiraoDaVeiga-Lovadina-Sacco:2017,
Artioli-BeiraoDaVeiga-Lovadina-Sacco:2017b,Park-Chi-Paulino:2021,
Gain-Talischi-Paulino:2014,Daltri-deMiranda-Patruno-Sacco:2021,
Tang-Liu-Zhang-Feng:2020,Cihan-Hudobivnik-Aldakheel-Wriggers:2021}, 
elastodynamics~\cite{Park-Chi-Paulino:2020,Park-Chi-Paulino:2021b}, 
finite deformations~\cite{Wriggers-Rust:2019,DeBellis-Wriggers-Hudobivnik:2019,
Aldakheel-Hudobivnik-Wriggers:2019,Wriggers-Hudobivnik:2017,Wriggers-Reddy-Rust-Hudobivnik:2017,
Hudobivnik-Aldakheel-Wriggers:2019,Zhang-Chi-Paulino:2020,Chi-BeiraoDaVeiga-Paulino:2017,
vanHuyssteen-Reddy:2020}, contact mechanics~\cite{Wriggers-Rust:2019,Wriggers-Rust-Reddy:2016,
Aldakheel-Hudobivnik-Artioli-BeiraoDaVeiga-Wriggers:2020}, 
fracture mechanics~\cite{Aldakheel-Hudobivnik-Wriggers:2019b,
Hussein-Aldakheel-Hudobivnik-Wriggers-Guidault-Allix:2019,
NguyenThanh-Zhuang-NguyenXuan-Rabczuk-Wriggers:2018,
Benedetto-Caggiano-Etse:2018,Artioli-Marfia-Sacco:2020}, 
fluid mechanics~\cite{BeiraoDaVeiga-Pichler-Vacca:2021,Chen-Wang:2019,
Gatica-Munar-Sequeira:2018,BeiraoDaVeiga-Lovadina-Vacca:2018,
Chernov-Marcati-Mascotto:2021}, geomechanics~\cite{Andersen-Nilsen-Raynaud:2017,
Lin-Zheng-Jiang-Li-Sun:2020} and 
topology optimization~\cite{Zhang-Chi-Paulino:2020,Paulino-Gain:2014,
Chi-Pereira-Menezes-Paulino:2020}. 

For linear elasticity, which is the focus in this paper, 
with the displacement field as the only variable (displacement-based approach), 
the simplest virtual element provides a linear approximation of the displacement field. 
As it occurs in traditional linear displacement-based finite elements 
(for instance, three-node triangles and four-node quadrilaterals), severe stiffening 
is exhibited by these virtual elements when the material is close to being incompressible, 
which occurs when the Poisson's ratio approaches $1/2$. This is a numerical
artifact known as volumetric locking. Currently, there are various
approaches to deal with this issue in the VEM literature (details
are provided below). Most of them are generalizations of well-established
FEM approaches to VEM. In this paper, we propose an 
approach for linear displacement-based virtual elements 
that is usable for compressible and nearly incompressible elasticity.
In the proposed approach, the strain is averaged at the nodes from
the strain of surrounding virtual elements, a technique that is 
very popular in meshfree and finite element nodal integration 
(nodal integration is discussed below) and results in locking-free formulations.} 
The proposed method differs from existing virtual element methods 
for compressible and nearly incompressible 
elasticity in the following aspects:
\begin{itemize} 
\item In the proposed approach, the volumetric locking is alleviated 
by the strain averaging at the nodes and the stabilization matrix chosen. 
No additional degrees of freedom are introduced in this approach, 
thus resulting in a displacement-based formulation. In the existing 
VEM approaches, the volumetric locking is alleviated by the B-bar formulation~\cite{Park-Chi-Paulino:2020,Park-Chi-Paulino:2021},
mixed formulation~\cite{Artioli-deMiranda-Lovadina-Patruno:2017}, enhanced strain
formulation~\cite{Daltri-deMiranda-Patruno-Sacco:2021}, hybrid formulation~\cite{Dassi-Lovadina-Visinoni:2021},
nonconforming formulations~\cite{Kwak-Park:2022,Zhang-Zhao-Yang-Chen:2019,Yu:2021}, or addition of 
degrees of freedom related to the normal components of the displacement field on the 
element's edges to satisfy the inf-sup condition~\cite{Tang-Liu-Zhang-Feng:2020}.
\item The strains and stresses in the proposed approach are projected to the nodes, thus
becoming nodal quantities just like displacements, whereas in the existing approaches 
these are projected inside the elements, thus becoming element quantities.
\end{itemize}

Nodal integration can be traced back to the work of Beissel and Belytschko~\cite{Beissel-Belytschko:1996} 
for meshfree Galerkin methods, where nodal integration emerged as an alternative
to Gauss integration to accelerate the numerical integration of the weak form integrals.
The main difference between these two types of numerical integration is that
in nodal integration the weak form integrals are sampled at the nodes, whereas in
Gauss integration the integrals are numerically integrated in the interior of the elements
using cubature. Thus, in nodal integration, the state and history-dependent variables in the
weak form integrals naturally become nodal quantities, whereas in Gauss integration
these variables become element quantities. The obvious advantages of nodal integration
over Gauss integration can be summarized as follows. There are much less nodes than
Gauss points in the mesh, which results in huge computational savings when nodal integration
is employed. In nonlinear computations with nodally integrated weak form integrals, 
the state and history-dependent variables are tracked only at the nodes. 
This feature can be exploited to avoid mesh remapping of these variables in 
Lagrangian large deformation simulations with remeshing 
(see, for instance, \cite{Puso-Chen-Zywicz-Elmer:2008}), which is not possible
when Gauss integration is used.

After the work of Beissel and Belytschko~\cite{Beissel-Belytschko:1996} nodal integration 
attracted great interest from the meshfree community (e.g., see Refs.~\cite{Belytschko-Xiao:2002,Fries-Belytschko:2008,Belytschko-Guo-Liu-Xiao:2000,Hillman-Chen:2016,Chen-Hillman-Ruter:2013,Puso-Solberg:2006,Chen-Wu-Yoon-You:2001,Chen-Yoon-Wu:2002,Elmer-Chen-Puso-Taciroglu:2012,Silva-Ortiz-Sukumar-Artioli-Hitschfeld:2020,Duan-Wang-Gao-Li:2014}) 
due to the improved robustness of this type of integration
over cell-based Gauss integration. Even though nodal integration greatly improves the computational
cost of meshfree Galerkin simulations, still there exists the computational burden of 
meshfree basis functions (e. g., moving least-squares~\cite{Lancaster-Salkauskas:1981,Belytschko-Lu-Gu:1994}, 
reproducing kernel~\cite{Liu-Jun-Zhang:1995} and 
maximum-entropy~\cite{Sukumar:2004,Sukumar-Wright:2007,Arroyo-Ortiz:2006} approximants), 
which includes the solution of an optimization problem and a neighbor search at every 
integration node. Besides, meshfree basis functions generally do not vanish on the boundary,
which requires additional procedures to impose Dirichlet boundary conditions.
Regarding finite elements, there are various nodal integration approaches
that use some form of pressure or strain averaging at the nodes; for instance, the average 
nodal pressure tetrahedral element~\cite{Bonet-Burton:1998}, node-based uniform 
strain triangular and tetrahedral elements~\cite{Dohrmann-Heinstein-Jung-Key-Witkowski:2000}, 
the averaged nodal deformation gradient linear tetrahedral element~\cite{Bonet-Marriott-Hassan:2001}, 
and the family of nodally integrated continuum elements (NICE) 
and its derivative approaches~\cite{Krysl-Zhu:2008,Castellazzi-Krysl:2009,Krysl-Kagey:2012,Castellazzi-Krysl:2012,Castellazzi-Krysl:2012,Broccardo-Micheloni-Krysl:2009,Castellazzi-Krysl-Bartoli:2013,Artioli-Castellazzi-Krysl:2014}. 
Recently, the node-based uniform 
strain approach~\cite{Dohrmann-Heinstein-Jung-Key-Witkowski:2000} was adopted in the nodal particle 
finite element method (N-PFEM)~\cite{Franci-Cremonesi-Perego-Onate:2020,Meng-Zhang-Utili-Onate:2021}.

Nodal integration is prone to instabilities, and thus it requires stabilization. 
Stabilization is also an important aspect within the theoretical framework of the 
VEM~\footnote{Although stabilization-free virtual elements 
were recently proposed~\cite{BerroneBorioMarcon:2022,ChenSukumar:2023}.}, 
and thus this represents the main connection between 
nodal integration and the virtual element method. A summary of the stability issue in 
nodal integration is provided in Ref.~\cite{Puso-Chen-Zywicz-Elmer:2008}, 
where also a gradient-based stabilization is proposed for small and Lagrangian 
large deformation simulations.

The VEM allows direct imposition of Dirichlet boundary conditions
and its computational cost is comparable to that of the finite element method~\cite{Ortiz-Alvarez-Hitschfeld-Russo-Silva-Olate:2019}. 
In addition, it is equipped with a stabilization procedure like in nodal integration. 
Regarding its sensitivity to mesh distortions, several authors have reported very robust 
solutions under mesh distortions~\cite{Artioli-BeiraoDaVeiga-Lovadina-Sacco:2017,Artioli-BeiraoDaVeiga-Lovadina-Sacco:2017b,Park-Chi-Paulino:2020,Wriggers-Reddy-Rust-Hudobivnik:2017,BeiraoDaVeiga-Dassi-Russo:2017,BeiraoDaVeiga-Brezzi-Dassi-Marini-Russo:2017} 
and highly nonconvex elements~\cite{Park-Chi-Paulino:2021,Park-Chi-Paulino:2021b,Chi-BeiraoDaVeiga-Paulino:2017,Paulino-Gain:2014}.
\alejandro{Since the VEM provides approximations 
that are interpolatory, stiffness matrix stabilization, 
and robustness under mesh distortion without the computational 
burden of meshfree methods}, it is our opinion that the VEM will play 
a key role in the development of robust nodal integration techniques, especially in
extreme large deformations. 

\alejandrog{In Ref.~\cite{Silva-Ortiz-Sukumar-Artioli-Hitschfeld:2020},
where a nodal integration for meshfree Galerkin methods is proposed, 
the meshfree stiffness matrix was nodally integrated 
at the Voronoi sites of a background Voronoi tessellation.
To render the nodal integration consistent and stable,
the nodally integrated meshfree stiffness matrix was constructed
using a meshfree version of the consistency and stability 
virtual element stiffness matrix decomposition.
The nodal integration at Voronoi sites, however, is not possible
in the standard VEM} \alejandro{because of the following two issues:
(a) it requires the evaluation of basis functions
from surrounding Voronoi sites at the edges of the Voronoi cell,
which cannot be performed in the VEM,
and (b)} \alejandrog{the Voronoi sites contain the degrees of freedom} 
\alejandro{while the VEM requires the degrees of freedom at the vertices of the
Voronoi cell. Therefore, to enable nodal integration in the} \alejandrog{standard} 
\alejandro{VEM a different approach must be devised.}

In this exploratory paper, a combined nodal integration and virtual element method
is proposed using ideas taken from the node-based uniform strain approach of 
Dohrmann et al.~\cite{Dohrmann-Heinstein-Jung-Key-Witkowski:2000} for finite elements,
wherein the strain is averaged at the nodes from the strain of surrounding elements. 
We refer to the proposed approach as the node-based uniform strain virtual element method (NVEM).
The governing equations for linear elastostatics
are provided in Section~\ref{sec:goveqns}. Section~\ref{sec:vem} summarizes the 
theoretical framework of the VEM for linear elasticity. The NVEM is developed 
in Section~\ref{sec:nbvem}. Stabilizations for the 
NVEM are discussed in Section~\ref{sec:stabilization}. In Section~\ref{sec:numexp}, 
the numerical behavior of the NVEM is examined through several benchmark problems 
in compressible and nearly incompressible elasticity,
\alejandro{and is also assessed in linear elastodynamics.}
Some concluding remarks along with ongoing and future research directions are discussed 
in Section~\ref{sec:conclusions}. To ease the reading of the equations, a list
of main symbols used throughout the paper is provided in~\ref{appA}.

\section{Governing equations}
\label{sec:goveqns}

Consider an elastic body that occupies the
open domain $\Omega \subset \Re^2$ and is
bounded by the one-dimensional surface $\Gamma$ whose
unit outward normal is $\vm{n}_{\Gamma}$. The boundary is assumed
to admit decompositions $\Gamma=\Gamma_D\cup\Gamma_N$ and
$\emptyset=\Gamma_D\cap\Gamma_N$, where $\Gamma_D$ is the Dirichlet
boundary and $\Gamma_N$ is the Neumann boundary. The closure of
the domain is $\overline{\Omega}=\Omega\cup\Gamma$. Let
$\vm{u}(\vm{x}) : \overline{\Omega} \rightarrow \Re^2$ be
the displacement field at a point of the elastic body 
with position vector $\vm{x}$ when the body is subjected to external tractions
$\vm{t}_N(\vm{x}):\Gamma_N\rightarrow \Re^2$ and body forces $\vm{b}(\vm{x}):\Omega\rightarrow\Re^2$.
The imposed Dirichlet (essential) boundary conditions are
$\vm{u}_D(\vm{x}):\Gamma_D\rightarrow \Re^2$. The boundary-value
problem for linear elastostatics reads: find
$\vm{u}(\vm{x}): \overline{\Omega} \rightarrow \Re^2$
such that
\begin{subequations}\label{eq:strongform}
\begin{align}
\bsym{\nabla} \cdot \bsym{\sigma} + \vm{b} &= 0 \quad \mathrm{in}\,\,\Omega, \\ 
\vm{u} &= \vm{u}_D \quad \mathrm{on}\,\,\Gamma_D, \\ 
\bsym{\sigma} \cdot \vm{n}_{\Gamma} &= \vm{t}_N \quad \mathrm{on}\,\,\Gamma_N,
\end{align}
\end{subequations}
where $\bsym{\sigma}$ is the Cauchy stress
tensor.

The Galerkin weak formulation of the above problem is:
find $\vm{u}(\vm{x})\in \mathcal{V}$ such that
\begin{equation}\label{eq:weakform}
\begin{split}
a(\vm{u},\vm{v}) & = \ell(\vm{v}) \quad \forall \vm{v}(\vm{x})\in \mathcal{W}, \\
a(\vm{u},\vm{v})=\int_{\Omega}\bsym{\sigma}(\vm{u}):\bsym{\varepsilon}(\vm{v})\,\diffx,
& \quad
\ell(\vm{v}) = \int_{\Omega}\vm{b}\cdot\vm{v}\,\diffx + \int_{\Gamma_N}\vm{t}_N\cdot\vm{v}\,\diffs,
\end{split}
\end{equation}
where $\mathcal{V}$ denotes the space of admissible displacements 
and $\mathcal{W}$ the space of its variations, and $\bsym{\varepsilon}$ is the 
small strain tensor that is given by
\begin{equation}\label{eq:symstrain}
\bsym{\varepsilon}(\vm{u})= 
\frac{1}{2}\left(\vm{u}\otimes\bsym{\nabla}+\bsym{\nabla}\otimes\vm{u}\right).
\end{equation}

\section{Virtual element method for linear elasticity}
\label{sec:vem}

The virtual element method (VEM)~\cite{BeiraoDaVeiga-Brezzi-Cangiani-Manzini-Marini-Russo:2013} 
is a generalization of the finite element method (FEM) that can
deal with very general meshes formed by elements with arbitrary number of edges
(convex or non-convex polygons); the mesh can even include elements with coplanar
edges and collapsing nodes, and non-matching elements. In the VEM, the displacement 
field is decomposed into a polynomial part and other terms using a projection operator. 
Altogether, this procedure results in \alejandro{an algebraic construction of the 
element stiffness matrix without the explicit construction of basis functions} 
(basis functions are \textit{virtual}). The decomposition of the displacement 
field leads to the following decomposition of the virtual element stiffness matrix:
\begin{equation}\label{eq:matrixdecomp}
\vm{K}_E = \vm{K}_E^\cons + \vm{K}_E^\stab,
\end{equation}
where $\vm{K}_E^\cons$ is the consistency part of the element stiffness matrix, which provides consistency
to the method (i.e., ensures patch test satisfaction), and $\vm{K}_E^\stab$ is the stability
part of the element stiffness matrix, which provides stability.

\begin{figure}[!bth]
\centering
\includegraphics[width=0.4\linewidth]{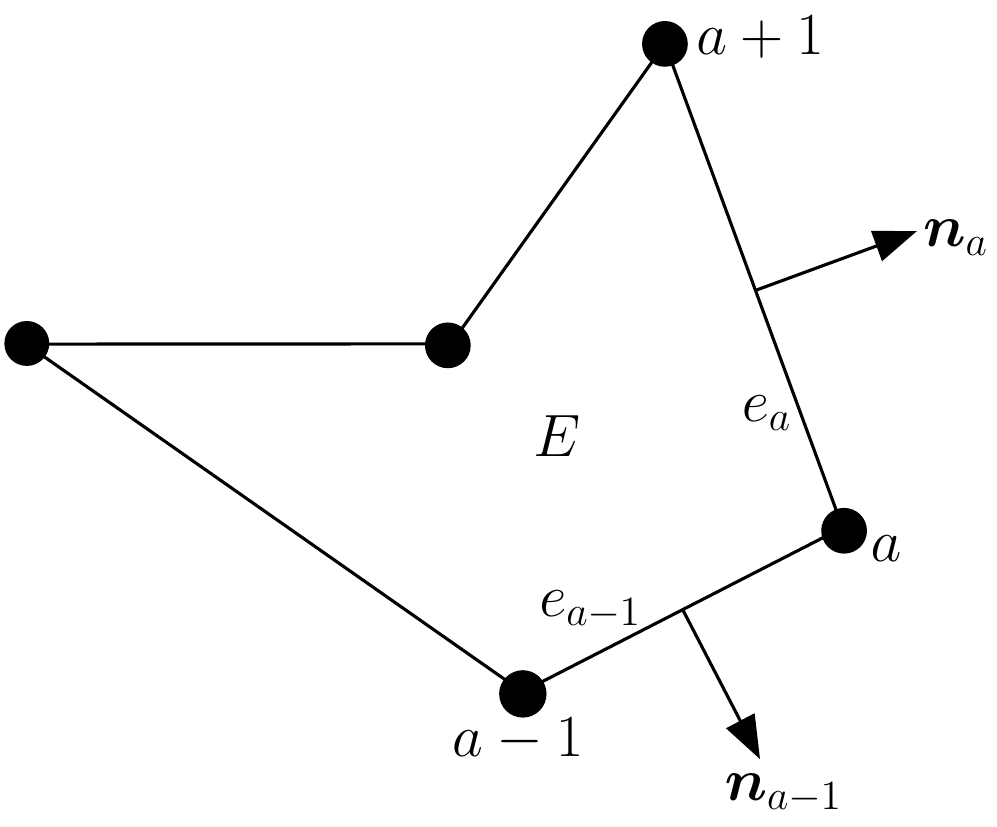}
\caption{Schematic representation of a polygonal element of $N_E^V=5$ edges.}
\label{fig:polyelement}
\end{figure}

A brief summary of the VEM for linear elasticity is given next. All the derivations are
presented at the element level. Concerning the element, the following notation is used:
the element is denoted by $E$ and its boundary by $\partial E$. The area of the element 
is denoted by $|E|$. The number of edges/nodes of an element are denoted by $N_E^V$. 
The unit outward normal to the element boundary in the Cartesian 
coordinate system is denoted by $\vm{n}=[n_1 \quad n_2]^\transpose$. 
\fref{fig:polyelement} depicts an element with five edges ($N_E^V=5$), where 
the edge $e_a$ of length $|e_a|$ and the edge $e_{a-1}$ of length $|e_{a-1}|$ are the 
element edges incident to node $a$, and $\vm{n}_a$ and $\vm{n}_{a-1}$ are the unit 
outward normals to these edges, respectively.

\subsection{Projection operator}
\label{sec:vem:projoperator}

The fundamental procedure in the VEM is the decomposition of the displacement field into its 
polynomial part and some additional terms that are generally nonpolynomials. 
To this end, a projection operator $\Pi$ is defined such that
\begin{equation}\label{eq:vem_decomp}
\vm{u}_h=\Pi\vm{u}_h+(\vm{u}_h-\Pi\vm{u}_h),
\end{equation}
where $\vm{u}_h$ is the displacement approximation on the element (trial function),
$\Pi\vm{u}_h$ is the polynomial part of $\vm{u}_h$ and $\vm{u}_h-\Pi\vm{u}_h$ contains the 
nonpolynomial terms of $\vm{u}_h$.

In this paper, we restrict ourselves to low-order approximations; that is, $\Pi$ is meant
to project onto the space of linear polynomials. Let $\vm{x}_a=\left[x_{1a} \quad x_{2a}\right]^\transpose$
and $\vm{u}_a=\left[u_{1a} \quad u_{2a}\right]^\transpose$ be the coordinates and displacements of 
node $a$, respectively. At the element level, the projection operator onto the linear displacements can 
be defined as~\cite{Gain-Talischi-Paulino:2014,Ortiz-Russo-Sukumar:2017,Ortiz-Alvarez-Hitschfeld-Russo-Silva-Olate:2019,Silva-Ortiz-Sukumar-Artioli-Hitschfeld:2020}
\begin{equation}\label{eq:operator}
\Pi\vm{u}_h = 
\smat{(x_1-\bar{x}_1) & 0 & \frac{1}{2}(x_2-\bar{x}_2) & 1 & 0 & \frac{1}{2}(x_2-\bar{x}_2)\\
       0 & (x_2-\bar{x}_2) & \frac{1}{2}(x_1-\bar{x}_1) & 0 & 1 & -\frac{1}{2}(x_1-\bar{x}_1)}
\smat{\widehat{\varepsilon}_{11} \\ \widehat{\varepsilon}_{22} 
       \\ 2\,\widehat{\varepsilon}_{12} \\ \bar{u}_{1} \\ \bar{u}_{2} 
       \\ 2\,\widehat{\omega}_{12} },
\end{equation}
where $\bar{x}_1$ and $\bar{x}_2$ are the components of the mean of the values that
the position vector $\vm{x}$ takes over the vertices of the element; i.e.,
\begin{equation}\label{eq:bar_x}
\bar{\vm{x}}=\smat{\bar{x}_1 \\ \bar{x}_2}=\frac{1}{N_E^V}\sum_{a=1}^{N_E^V}\vm{x(\vm{x}_a)},
\end{equation}
and $\bar{u}_{1}$ and $\bar{u}_{2}$ are the components of the mean of the values
that the displacement $\vm{u}$ takes over the vertices of the element; i.e.,
\begin{equation}\label{eq:bar_disp}
\bar{\vm{u}}=\smat{\bar{u}_{1} \\ \bar{u}_{2}}=\frac{1}{N_E^V}\sum_{a=1}^{N_E^V}\vm{u}(\vm{x}_a).
\end{equation}
In other words, $\bar{\vm{x}}$ and $\bar{\vm{u}}$ represents the geometric center of
the element and its associated displacement vector, respectively;
the terms $\widehat{\varepsilon}_{ij}$ are components of the element average
$\widehat{\bsym{\varepsilon}}=\frac{1}{|E|}\int_E \bsym{\varepsilon}\,\diffx$, 
and $\widehat{\omega}_{12}$ is the component of the element average
$\widehat{\bsym{\omega}}=\frac{1}{|E|}\int_E \bsym{\omega}\,\diffx$,
where $\bsym{\omega}$ is the skew-symmetric tensor that represents rotations. 
These element averages are evaluated on the boundary of $E$ by invoking the divergence
theorem, which gives the following identities:

\begin{equation}\label{eq:avg_strain}
\widehat{\bsym{\varepsilon}}(\vm{u}_h)=\frac{1}{|E|}\int_E \bsym{\varepsilon}(\vm{u}_h)\,\diffx = \frac{1}{2|E|}\int_{\partial
E}\left(\vm{u}_h\otimes\vm{n}+\vm{n}\otimes\vm{u}_h\right)\,\diffs,
\end{equation}

\begin{equation}\label{eq:avg_skewstrain}
\widehat{\bsym{\omega}}(\vm{u}_h)=\frac{1}{|E|}\int_E \bsym{\omega}(\vm{u}_h)\,\diffx=\frac{1}{2|E|}\int_{\partial
E}\left(\vm{u}_h\otimes\vm{n}-\vm{n}\otimes\vm{u}_h\right)\,\diffs.
\end{equation}

\subsection{Discretization of field variables}
\label{sec:vem:trialtestfunc}

Following a Galerkin approach, the element trial and test functions, $\vm{u}_h$ and $\vm{v}_h$,
respectively, are assumed to be discretized as
\begin{equation}\label{eq:discrete_field}
\vm{u}_h=\smat{u_{1h} \\ u_{2h}}=\sum_{a=1}^{N_E^V}\phi_a(\vm{x})\vm{u}_a, \quad
\vm{v}_h=\smat{v_{1h} \\ v_{2h}}=\sum_{a=1}^{N_E^V}\phi_a(\vm{x})\vm{v}_a,
\end{equation}
where $\{\phi_a(\vm{x})\}_{a=1}^{N_E^V}$ are basis functions that form a partition of unity.
Eqns. \eref{eq:bar_disp}--\eref{eq:avg_skewstrain} reveal that all
what is needed regarding the basis functions is the knowledge of their behavior on the element
boundary. Thus, in the linear VEM the following behavior for the basis functions is assumed
on the element boundary: 
\begin{itemize}
	\item basis functions are piecewise linear (edge by edge);
	\item basis functions are continuous on the element edges.
\end{itemize}
This means that the basis functions possess the Kronecker-delta property on the element edges,
and hence they behave like the one-dimensional hat function. 

When establishing the element matrices, we will see that these matrices are written 
in terms of the basis functions on the element boundary through~\eref{eq:avg_strain} 
and the discretization of the displacements using~\eref{eq:discrete_field}.
The element matrices will be then computed algebraically using the above assumptions 
in such a way that the basis functions are never computed. Thus, it is stressed 
that in the VEM the basis functions are not (and do not need to be) known explicitly.

\subsection{VEM bilinear form}
\label{sec:vem:bilinearform}
The projection operator \eref{eq:operator} is derived\footnote{The derivation is given in Ref.~\cite{Silva-Ortiz-Sukumar-Artioli-Hitschfeld:2020}} 
from the following orthogonality condition at the element level:
\begin{equation}\label{eq:orthogonal}
a_E(\vm{u}_h-\Pi\vm{u}_h,\vm{p})=a_E(\vm{p},\vm{v}_h-\Pi\vm{v}_h)=0 \quad \forall\vm{p} \in [\mathcal{P}(E)]^2,
\end{equation}
where $[\mathcal{P}(E)]^2$ is the space of linear displacements over the element. 
This orthogonality condition states that the nonpolynomial terms $\vm{u}_h-\Pi\vm{u}_h$ 
in $E$, measured in the energy norm, are orthogonal to $[\mathcal{P}(E)]^2$.
Using the decomposition \eref{eq:vem_decomp} along with condition \eref{eq:orthogonal},
and noting that $\Pi\vm{u}_h$ and $\Pi\vm{v}_h$ $\in [\mathcal{P}(E)]^2$,
the following VEM representation of the bilinear form is obtained:
\begin{equation}\label{eq:VEM_bilinearform1}
a_E(\vm{u}_h,\vm{v}_h)=a_E(\Pi\vm{u}_h,\Pi\vm{v}_h)+a_E(\vm{u}_h-\Pi\vm{u}_h,\vm{v}_h-\Pi\vm{v}_h).
\end{equation}

The VEM bilinear form can be further elaborated. For convenience, the projection operator
\eref{eq:operator} is rewritten as
\begin{equation}\label{eq:projmat1}
\Pi\vm{u}_h = \vm{h}(\vm{x})\widehat{\bsym{\varepsilon}}(\vm{u}_h)+\vm{g}(\vm{x})\vm{r}(\vm{u}_h),
\end{equation}
where
\begin{equation}\label{eq:mtrxhg}
\vm{h}(\vm{x}) = \smat{(x_1-\bar{x}_1) & 0 & \frac{1}{2}(x_2-\bar{x}_2)\\ 0 & (x_2-\bar{x}_2) & \frac{1}{2}(x_1-\bar{x}_1)},\quad
\vm{g}(\vm{x}) = \smat{1 & 0 & \frac{1}{2}(x_2-\bar{x}_2)\\ 0 & 1 & -\frac{1}{2}(x_1-\bar{x}_1)},
\end{equation}
\begin{equation}\label{eq:mtrxr}
\vm{r}(\vm{u}_h) = \smat{\bar{u}_{1} & \bar{u}_{2} & 2\,\widehat{\omega}_{12}}^\transpose,
\end{equation}
and
\begin{equation}\label{eq:vectoravgstrain}
\widehat{\bsym{\varepsilon}}(\vm{u}_h) = 
\smat{\widehat{\varepsilon}_{11} 
      & \widehat{\varepsilon}_{22} 
      & 2\,\widehat{\varepsilon}_{12}}^\transpose
\end{equation}
is the vector form of the element average strain. From now on, we always use
this vector form when referring to the element average strain $\widehat{\bsym{\varepsilon}}$.

Using the bilinear form written in terms of the vector form of the strain --- that is,
$a(\vm{u}_h,\vm{v}_h)=\int_\Omega \bsym{\varepsilon}^\transpose(\vm{v}_h)\,\vm{D}\,\bsym{\varepsilon}(\vm{u}_h)\,\diffx$,
where $\vm{D}$ is the constitutive matrix --- and
noting that $\vm{r}(\vm{u}_h)$ represents rigid body translations and rotation, and that these
do not produce strain energy, the first term on the right-hand side of \eref{eq:VEM_bilinearform1}
can be written, as follows:
\begin{align}\label{eq:VEM_consistencybilinearform}
a_E(\Pi\vm{u}_h,\Pi\vm{v}_h)&=
\int_E \smat{\frac{\partial\Pi v_{1h}}{\partial x_1} 
             & \frac{\partial\Pi v_{2h}}{\partial x_2} 
             & \frac{\partial\Pi v_{1h}}{\partial x_2}+\frac{\partial\Pi v_{2h}}{\partial x_1}}
             \,\vm{D}\,
       \smat{\frac{\partial\Pi u_{1h}}{\partial x_1} \\
             \frac{\partial\Pi u_{2h}}{\partial x_2} \\
             \frac{\partial\Pi u_{1h}}{\partial x_2}+\frac{\partial\Pi u_{2h}}{\partial x_1}}\,\diffx
\nonumber\\
&=\int_E \widehat{\bsym{\varepsilon}}^\transpose(\vm{v}_h)\,\vm{D}\,\widehat{\bsym{\varepsilon}}(\vm{u}_h)\,\diffx.
\end{align}
For an isotropic linear elastic material, the constitutive matrix is given by
\begin{equation}\label{eq:pstrainmat}
\vm{D} =
\frac{E_\mathrm{Y}}{(1+\nu)(1-2\nu)}\smat{1-\nu & \nu & 0\\ \nu & 1-\nu & 0\\ 0 & 0 & \frac{1-2\nu}{2}}
\end{equation}
for plane strain condition, and
\begin{equation}\label{eq:pstressmat}
\vm{D} =
\frac{E_\mathrm{Y}}{1-\nu^2}\smat{1 & \nu & 0 \\ \nu & 1 & 0 \\ 0 & 0 & \frac{1-\nu}{2}}
\end{equation}
for plane stress condition, where $E_\mathrm{Y}$ is the Young's modulus and $\nu$ is 
the Poisson's ratio.

The second term on the right-hand side of \eref{eq:VEM_bilinearform1}
can be written as
\begin{equation}\label{eq:VEM_stabilitybilinearform}
a_E(\vm{u}_h-\Pi\vm{u}_h,\vm{v}_h-\Pi\vm{v}_h) = (\vm{1}-\Pi)^\transpose a_E(\vm{v}_h,\vm{u}_h)(\vm{1}-\Pi),
\end{equation}
which reveals that $a_E(\vm{u}_h-\Pi\vm{u}_h,\vm{v}_h-\Pi\vm{v}_h)$ is not computable
since $a_E(\vm{v}_h,\vm{u}_h)$ would need to be evaluated inside the element using cubature, 
but the VEM basis functions are not known inside the element. 
However, since $a_E(\vm{u}_h-\Pi\vm{u}_h,\vm{v}_h-\Pi\vm{v}_h)$ is meant to only provide
stability, $a_E(\vm{v}_h,\vm{u}_h)$ can be conveniently approximated by a bilinear 
form $s_E(\cdot,\cdot)$ that is computable while preserving the coercivity
of the system. In other words, $s_E(\vm{v}_h,\vm{u}_h)$ 
must be positive definite and scale like $a_E(\cdot,\cdot)$ in 
the kernel of $\Pi$~\cite{BeiraoDaVeiga-Brezzi-Cangiani-Manzini-Marini-Russo:2013,Gain-Talischi-Paulino:2014}.

Using~\eref{eq:VEM_consistencybilinearform} and \eref{eq:VEM_stabilitybilinearform} with
the approximation $a_E(\vm{v}_h,\vm{u}_h)=s_E(\vm{v}_h,\vm{u}_h)$, the VEM bilinear
form for linear elasticity is written as
\begin{equation}\label{eq:VEM_bilinearform_planeelast}
a_{h,E}(\vm{u}_h,\vm{v}_h)=
\int_E\widehat{\bsym{\varepsilon}}^\transpose(\vm{v}_h)\,\vm{D}\,\widehat{\bsym{\varepsilon}}(\vm{u}_h)\,\diffx
\,+\,(\vm{1}-\Pi)^\transpose s_E(\vm{v}_h,\vm{u}_h)(\vm{1}-\Pi)
\end{equation}

\subsection{VEM element stiffness matrix}
\label{sec:vem:stiffness}

The discrete element average strain is computed using~\eref{eq:discrete_field}, as follows:
\begin{equation}\label{eq:disc_avgstrain}
\widehat{\bsym{\varepsilon}}_h(\vm{u}_h)=\widehat{\bsym{\varepsilon}}\left(\sum_{a=1}^{N_E^V}\phi_a(\vm{x})\vm{u}_a\right)=\vm{B}\vm{d},
\end{equation}
where
\begin{equation}\label{eq:vecnodaldisp}
\vm{d} = \smat{u_{11} & u_{21} & \cdots & u_{1a} & u_{2a} & \cdots & u_{1N_E^V} & u_{2N_E^V}}^\transpose
\end{equation}
is the column vector of element nodal displacements and 

\begin{equation}\label{eq:matrix_B}
\vm{B}=\smat{
\vm{B}_1 & \cdots \vm{B}_a \cdots & \vm{B}_{N_E^V}}, \quad
\vm{B}_a = \smat{q_{1a} & 0 \\ 0 & q_{2a} \\ q_{2a} & q_{1a}},
\end{equation}
where $q_{ia}=\frac{1}{|E|}\int_{\partial E}\phi_a(\vm{x})n_i\,\diffs$.
Since the VEM basis functions are assumed to be piecewise linear (edge by edge)
and continuous on the element boundary, $q_{ia}$ can be exactly computed using a trapezoidal
rule, which gives the following algebraic expression:
\begin{equation}\label{eq:qia}
q_{ia}=\frac{1}{2|E|}\left(|e_{a-1}| n_{i(a-1)}+|e_a| n_{ia}\right),\quad i=1,2,
\end{equation}
where $n_{ia}$ is the $i$-th component of $\vm{n}_a$ and $|e_a|$ is the length of
the edge incident to node $a$ as defined in \fref{fig:polyelement}. Using~\eref{eq:disc_avgstrain}
and noting that $\widehat{\bsym{\varepsilon}}$ and $\vm{D}$ are constant over the element,
the discrete version of the first term on the right-hand side of 
\eref{eq:VEM_bilinearform_planeelast} is written as
\begin{equation}\label{eq:discrete_consistency}
\int_E\widehat{\bsym{\varepsilon}}_{h}^\transpose(\vm{v}_h)\,\vm{D}\,\widehat{\bsym{\varepsilon}}_{h}(\vm{u}_h)\,\diffx
=|E|\,\widehat{\bsym{\varepsilon}}_{h}^\transpose(\vm{v}_h)\,\vm{D}\,\widehat{\bsym{\varepsilon}}_{h}(\vm{u}_h)=
\vm{q}^\transpose\left(|E|\vm{B}^\transpose \vm{D}\vm{B}\right)\vm{d},
\end{equation}
where $\vm{q}$ is a column vector similar to~\eref{eq:vecnodaldisp} 
that contains the element nodal values associated with $\vm{v}_h$.

To obtain the discrete version of the second term on the right-hand side of 
\eref{eq:VEM_bilinearform_planeelast}, the discretization
of $\Pi$ is required. This demands the discretization of $h(\vm{x})$, $\widehat{\bsym{\varepsilon}}(\vm{u}_h)$,
$g(\vm{x})$ and $r(\vm{u}_h)$. The discretization of $\widehat{\bsym{\varepsilon}}(\vm{u}_h)$ is given in \eref{eq:disc_avgstrain}.
The remainder terms are discretized, as follows:
\begin{equation}\label{eq:disc_hx}
\vm{h}_h(\vm{x}) = \sum_{a=1}^{N_E^V}\phi_a(\vm{x})\vm{h}(\vm{x}_a)=\vm{N}\vm{H},
\end{equation}
\begin{equation}\label{eq:disc_gx}
\vm{g}_h(\vm{x}) = \sum_{a=1}^{N_E^V}\phi_a(\vm{x})\vm{g}(\vm{x}_a)=\vm{N}\vm{G},
\end{equation}
\begin{equation}\label{eq:disc_ruh}
\vm{r}_h=\vm{r}\left(\sum_{a=1}^{N_E^V}\phi_a(\vm{x})\vm{u}_a\right)=\vm{R}\vm{d},
\end{equation}
where
\begin{equation}\label{eq:matrix_N}
\vm{N} =
\smat{\vm{N}_1 & \cdots & \vm{N}_a & \cdots & \vm{N}_{N_E^V}}, \quad \vm{N}_a=\smat{\phi_a(\vm{x}) & 0 \\ 0 &
\phi_a(\vm{x})},
\end{equation}

\begin{equation}\label{eq:matrix_H}
\vm{H}=\smat{
\vm{H}_1\\ \vdots\\ \vm{H}_a\\ \vdots\\ \vm{H}_{N_E^V}}, \quad
\vm{H}_a = \smat{(x_{1a}-\bar{x}_1) & 0 & \frac{1}{2}(x_{2a}-\bar{x}_2) \\ 0 & (x_{2a}-\bar{x}_2) & \frac{1}{2}(x_{1a}-\bar{x}_1)},
\end{equation}

\begin{equation}\label{eq:matrix_G}
\vm{G}=\smat{
\vm{G}_1\\ \vdots\\ \vm{G}_a\\ \vdots\\ \vm{G}_{N_E^V}}, \quad
\vm{G}_a = \smat{1 & 0 & \frac{1}{2}(x_{2a}-\bar{x}_2) \\ 0 & 1 & -\frac{1}{2}(x_{1a}-\bar{x}_1)},
\end{equation}
and

\begin{equation}\label{eq:matrix_R}
\vm{R}=\smat{
\vm{R}_1 & \cdots & \vm{R}_a & \cdots & \vm{R}_{N_E^V}}, \quad
\vm{R}_a = \smat{\frac{1}{N_E^V} & 0 \\ 0 & \frac{1}{N_E^V} \\ q_{2a} & -q_{1a}}.
\end{equation}
Therefore, the discrete $\Pi$ can be written as
\begin{equation}\label{eq:discrete_pi}
\Pi_h\vm{u}_h=\vm{N}\vm{H}\vm{B}\vm{d}+\vm{N}\vm{G}\vm{R}\vm{d}=\vm{N}(\vm{H}\vm{B}+\vm{G}\vm{R})\vm{d},
\end{equation}
which defines the projection matrix as
\begin{equation}\label{eq:pimatrix}
\vm{P}=\vm{H}\vm{B}+\vm{G}\vm{R}.
\end{equation}

Using the projection matrix and the discrete field $\vm{u}_h=\sum_{a=1}^{N_E^V}\phi_a(\vm{x})\vm{u}_a=\vm{N}\vm{d}$,
the discrete version of the second term on the right-hand side of \eref{eq:VEM_bilinearform_planeelast} is
\begin{equation}\label{eq:discrete_stability}
\vm{q}^\transpose(\vm{I}-\vm{P})^\transpose s_E(\vm{N}^\transpose,\vm{N})(\vm{I}-\vm{P})\vm{d}
=\vm{q}^\transpose(\vm{I}-\vm{P})^\transpose\,\vm{S}\,(\vm{I}-\vm{P})\vm{d}.
\end{equation}

Finally, using \eref{eq:discrete_consistency} and \eref{eq:discrete_stability}, the 
discrete VEM bilinear form is written as
\begin{equation}\label{eq:VEM_discrete_bilinearform}
a_{h,E}(\vm{u}_h,\vm{v}_h) = 
\vm{q}^\transpose\left\{|E|\,\vm{B}^\transpose\,\vm{D}\,\vm{B}+(\vm{I}-\vm{P})^\transpose\,\vm{S}\,(\vm{I}-\vm{P})\right\}\vm{d},
\end{equation}
which defines the element stiffness matrix as the sum of the element consistency stiffness matrix, $\vm{K}_E^\cons$, and
the element stability stiffness matrix, $\vm{K}_E^\stab$, as follows:
\begin{equation}\label{eq:VEM_stiffness_matrix}
\vm{K}_E = \vm{K}_E^\cons + \vm{K}_E^\stab,\quad
\vm{K}_E^\cons=|E|\,\vm{B}^\transpose\,\vm{D}\,\vm{B},\quad
\vm{K}_E^\stab = (\vm{I}-\vm{P})^\transpose\,\vm{S}\,(\vm{I}-\vm{P}).
\end{equation}

\subsection{VEM element force vector}

For linear fields, the VEM version of the loading term in the weak form \eref{eq:weakform} is computed
at the element level, as follows~\cite{BeiraoDaVeiga-Brezzi-Cangiani-Manzini-Marini-Russo:2013,BeiraodaVeiga-Brezzi-Marini:2013,Artioli-BeiraoDaVeiga-Lovadina-Sacco:2017}: 
\begin{equation}\label{eq:VEM_loading_b}
  \ell_{E}^b(\vm{v}_h) = \int_E\widehat{\vm{b}}\cdot\bar{\vm{v}}_h\,\diffx = |E|\,\widehat{\vm{b}}\cdot\bar{\vm{v}}_h,\quad \widehat{\vm{b}}=\frac{1}{|E|}\int_E\vm{b}\,\diffx, \quad \bar{\vm{v}}_h=\frac{1}{N_E^V}\sum_{a=1}^{N_E^V}\vm{v}_h(\vm{x}_a)
\end{equation}
for the loading term associated with the body forces, and similarly 
\begin{equation}\label{eq:VEM_loading_t}
  \ell_{e}^t(\vm{v}_h) = |e|\,\widehat{\vm{t}}_{N}\cdot\bar{\vm{v}}_h,\quad \widehat{\vm{t}}_{N}=\frac{1}{|e|}\int_e\vm{t}_N\,\diffs,
  \quad \bar{\vm{v}}_h=\frac{1}{N_e^V}\sum_{a=1}^{N_e^V}\vm{v}_h(\vm{x}_a)
\end{equation}
for the loading term associated with the tractions applied on an edge located on the Neumann boundary, 
where $e$ and $|e|$ are as defined in \fref{fig:polyelement} and $N_e^V$ is the 
number of nodes of the edge. 
The discrete versions of these loading terms are
$\ell_{h,E}^b(\vm{v}_h)=\vm{q}^\transpose\vm{f}_E^b$ and $\ell_{h,e}^t(\vm{v}_h)=\vm{q}^\transpose\vm{f}_e^t$,
where $\vm{f}_E^b$ and $\vm{f}_e^t$ are the element force vectors given, respectively, by
\begin{equation}\label{eq:VEM_forcevector_b}
 \vm{f}_E^b = |E|\,\bar{\vm{N}}^\transpose \widehat{\vm{b}},
\end{equation}
where
\begin{equation}\label{eq:bar_N}
  \bar{\vm{N}} = \smat{\bar{\vm{N}}_1 & \cdots & \bar{\vm{N}}_a & \cdots & \bar{\vm{N}}_{N_E^V}},
  \quad \bar{\vm{N}}_a = \smat{\frac{1}{N_E^V} & 0 \\ 0 & \frac{1}{N_E^V}},
\end{equation}
and
\begin{equation}\label{eq:VEM_forcevector_t}
 \vm{f}_e^t = |e|\,\bar{\vm{N}}_\Gamma^\transpose\, \widehat{\vm{t}}_{N},
\end{equation}
where
\begin{equation}\label{eq:bar_N_gamma}
  \bar{\vm{N}}_\Gamma = \smat{\frac{1}{N_e^V} & 0 & \frac{1}{N_e^V} & 0 \\ 0 & \frac{1}{N_e^V} & 0 & \frac{1}{N_e^V}} =
  \smat{\frac{1}{2} & 0 & \frac{1}{2} & 0 \\ 0 & \frac{1}{2} & 0 & \frac{1}{2}}
\end{equation}
since the number of nodes of an edge is $N_e^V=2$.

\section{Node-based uniform strain virtual element method (NVEM)}
\label{sec:nbvem}

The virtual element method for linear elasticity described in Section~\ref{sec:vem} is
prone to volumetric locking in the limit $\nu \to 1/2$. 
\alejandro{Using the virtual element mesh and considering a typical
nodal vertex $I$,} we propose to use a nodal 
averaging operator $\pi_I$ in~\eref{eq:VEM_bilinearform1} that is designed to preclude
volumetric locking without introducing additional degrees of freedom. 
This nodal averaging operator leads to a nodal version of the VEM bilinear form,
as follows:
\begin{equation}\label{eq:VEM_nodal_bilinearform}
a_{h,I}(\vm{u}_h,\vm{v}_h)=a_I(\pi_I[\Pi\vm{u}_h],\pi_I[\Pi\vm{v}_h])
+s_I(\pi_I[\vm{u}_h-\Pi\vm{u}_h],\pi_I[\vm{v}_h-\Pi\vm{v}_h]),
\end{equation}
where according to the VEM theory, 
$s_I(\pi_I[\vm{u}_h-\Pi\vm{u}_h],\pi_I[\vm{v}_h-\Pi\vm{v}_h])$ 
is a computable approximation to $a_I(\pi_I[\vm{u}_h-\Pi\vm{u}_h],\pi_I[\vm{v}_h-\Pi\vm{v}_h])$
and is meant to provide stability. \alejandro{The notations $a_I$ and $s_I$ are
introduced as the nodal counterparts of the element bilinear form $a_E$ 
and its stabilization term $s_E$, respectively}.

\subsection{Construction of the nodal averaging operator}

\begin{figure}[!bth]
\centering
\includegraphics[width=0.45\linewidth]{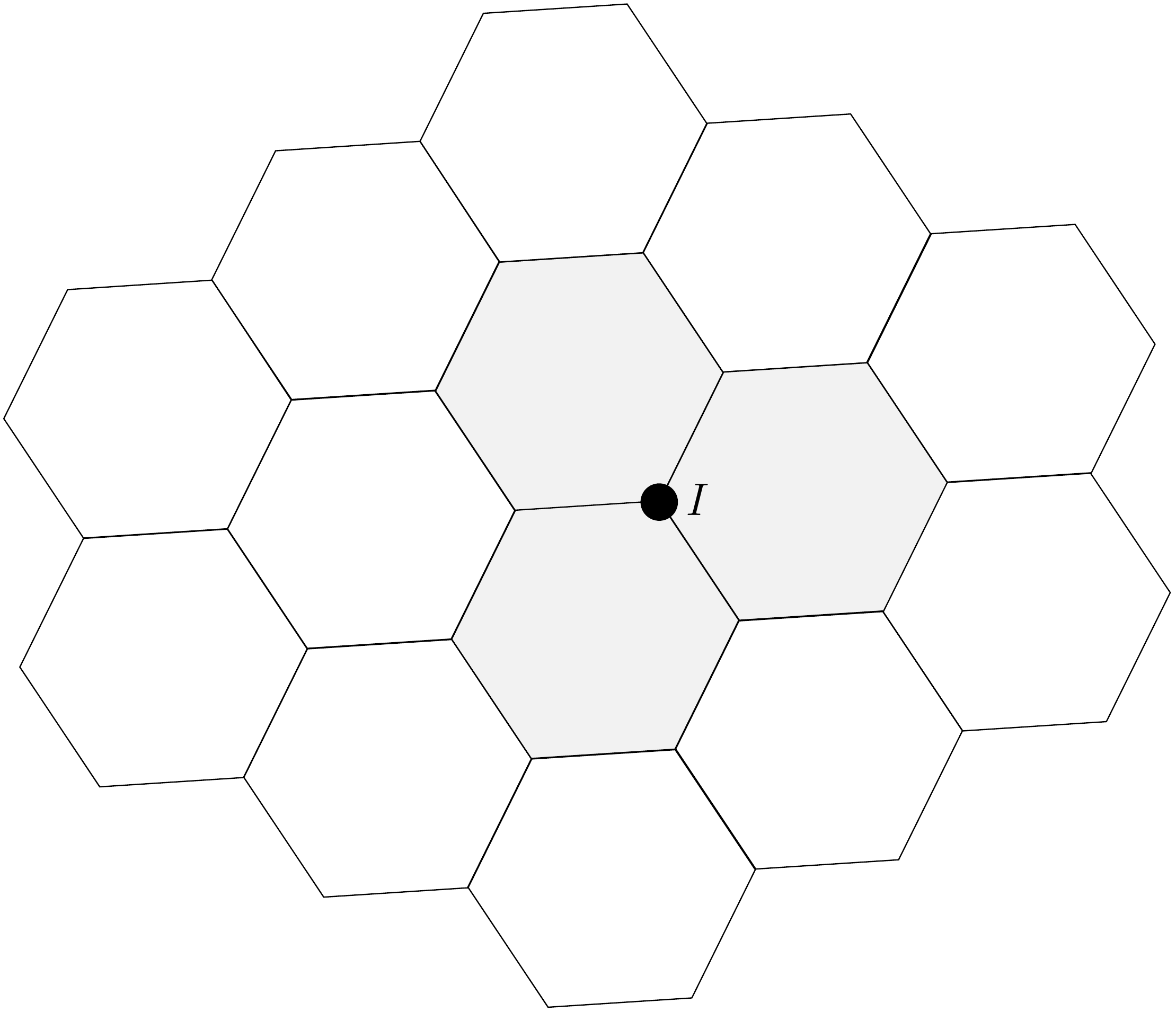}
\caption{Nodal patch $\mathcal{T}_I$ (shaded elements) formed by the virtual elements that are connected to node $I$.}
\label{fig:nodalpatch}
\end{figure}

Each node of the mesh is associated with their own patch of virtual elements.
The patch for node $I$ is denoted by $\mathcal{T}_I$ and is defined as the set of 
virtual elements connected to node $I$ (see~\fref{fig:nodalpatch}). Each node of a virtual element $E$ 
in the patch is assigned the area $\frac{1}{N_E^V}|E|$; that is,
the area of an element is uniformly distributed among its nodes. The representative 
area of node $I$ is denoted by $|I|$ and is computed by addition of all the areas 
that are assigned to node $I$ from the elements in $\mathcal{T}_I$; that is,
\begin{equation}\label{eq:nodal_area}
|I|=\sum_{E\in\mathcal{T}_I}\frac{1}{N_E^V}|E|.
\end{equation}
Similarly, each node of a virtual element $E$ is uniformly assigned the 
strain $\frac{1}{N_E^V}\widehat{\bsym{\varepsilon}}(\vm{u}_h)$. On considering each 
strain assigned to node $I$ from the elements in $\mathcal{T}_I$, we define 
the \textit{node-based uniform strain}, as follows:
\begin{equation}\label{eq:nodal_strain}
\widehat{\bsym{\varepsilon}}_I(\vm{u}_h)=
\frac{1}{|I|}\sum_{E\in\mathcal{T}_I}|E|\frac{1}{N_E^V}\widehat{\bsym{\varepsilon}}(\vm{u}_h).
\end{equation}
Since $\widehat{\bsym{\varepsilon}}(\vm{u}_h)$ is by definition given at the element level, then
from Eq.~\eref{eq:nodal_strain} the following nodal averaging operator is proposed:
\begin{equation}\label{eq:nodal_avg_operator}
\pi_I[\,\cdot\,]=\frac{1}{|I|}\sum_{E\in\mathcal{T}_I}|E|\frac{1}{N_E^V}[\,\cdot\,]_{{}_E},
\end{equation}
where $[\,\cdot\,]_{{}_E}$ denotes evaluation over the element $E$. 
Using this nodal averaging operator, the nodal representation of a function $F$
is obtained as
\begin{equation}\label{eq:nodal_avg_term}
F_I=\pi_I[\,F\,]=\frac{1}{|I|}\sum_{E\in\mathcal{T}_I}|E|\frac{1}{N_E^V}[\,F\,]_{{}_E}.
\end{equation}

The idea of having a nodal strain (like the one in~\eref{eq:nodal_strain}) 
that is representative of the nodal patch $\mathcal{T}_I$ has long been used in meshfree 
and finite element nodal integration techniques~\cite{Bonet-Burton:1998,Dohrmann-Heinstein-Jung-Key-Witkowski:2000,
Bonet-Marriott-Hassan:2001,Andrade-Souza-Cuesta:2004,Puso-Solberg:2006,Puso-Chen-Zywicz-Elmer:2008,Krysl-Zhu:2008,
Broccardo-Micheloni-Krysl:2009,Castellazzi-Krysl:2009,Castellazzi-Krysl:2012}. In particular, 
the node-based uniform strain approach was developed for three-node triangular and 
four-node tetrahedral elements in Ref.~\cite{Dohrmann-Heinstein-Jung-Key-Witkowski:2000}, 
and proved to perform well for nearly incompressible solids. However, that approach presented 
stability issues in some particular instances and required stabilization~\cite{Puso-Chen-Zywicz-Elmer:2008}. 
Within the VEM framework, stabilization is obtained by construction of the method, and hence 
in combining the node-based uniform strain approach with the VEM no stability issues will
arise (\alejandro{we demonstrate this through numerical experiments in Section~\ref{sec:numexp}}).

\subsection{NVEM bilinear form}
Using the nodal averaging operator~\eref{eq:nodal_avg_operator}, the nodal representation 
of the VEM bilinear form for linear elasticity (Eq.~\eref{eq:VEM_bilinearform_planeelast}) 
is written as
\begin{equation}\label{eq:VEM_nodal_VEM_bilinearform}
a_{h,I}(\vm{u}_h,\vm{v}_h)=|I|\,\widehat{\bsym{\varepsilon}}_I^\transpose(\vm{v}_h)\,\vm{D}\,\widehat{\bsym{\varepsilon}}_I(\vm{u}_h)
+(\vm{1}-\Pi)_I^\transpose s_I(\vm{v}_h,\vm{u}_h) (\vm{1}-\Pi)_I.
\end{equation}

Let $N$ be the total number of nodes in the mesh. The global bilinear form is 
obtained by summing through the $N$ nodes, as follows:
\begin{equation}\label{eq:VEM_global_nodal_VEM_bilinearform}
a_h(\vm{u}_h,\vm{v}_h)=\sum_{I=1}^N\Bigg[|I|\,\widehat{\bsym{\varepsilon}}_I^\transpose(\vm{v}_h)\,\vm{D}\,\widehat{\bsym{\varepsilon}}_I(\vm{u}_h)
+(\vm{1}-\Pi)_I^\transpose s_I(\vm{v}_h,\vm{u}_h) (\vm{1}-\Pi)_I\Bigg].
\end{equation}

\subsection{NVEM nodal stiffness matrix}
The nodal stiffness matrix that arises from~\eref{eq:VEM_nodal_VEM_bilinearform}
is the nodal version of~\eref{eq:VEM_stiffness_matrix}. This gives
\begin{equation}\label{eq:nodal_VEM_stiffness_matrix}
\vm{K}_I = \vm{K}_I^\cons + \vm{K}_I^\stab,\quad
\vm{K}_I^\cons=|I|\,\vm{B}_I^\transpose\,\vm{D}\,\vm{B}_I,\quad
\vm{K}_I^\stab = (\vm{I}-\vm{P})_I^\transpose\,\vm{S}_I\,(\vm{I}-\vm{P})_I,
\end{equation}
where $\vm{B}_I=\pi_I[\,\vm{B}\,]$, $(\vm{I}-\vm{P})_I=\pi_I[\,\vm{I}-\vm{P}\,]$, and
\alejandro{the stability matrix} $\vm{S}_I$ is defined in Section~\ref{sec:stabilization}.

\subsection{NVEM nodal force vector}

The nodal force vector associated with the body forces
is the nodal representation of~\eref{eq:VEM_forcevector_b} and is given by
\begin{equation}\label{eq:nodal_VEM_forcevector_b}
 \vm{f}_I^b = |I|\,\bar{\vm{N}}_I^\transpose\widehat{\vm{b}}_I,
\end{equation}
where $\bar{\vm{N}}_I=\pi_I[\,\bar{\vm{N}}\,]$ and $\widehat{\vm{b}}_I=\pi_I[\,\widehat{\vm{b}}\,]$.
Similarly, the nodal force vector associated with the tractions is given
by the nodal representation of~\eref{eq:VEM_forcevector_t}, as follows:
\begin{equation}\label{eq:nodal_VEM_forcevector_t}
 \vm{f}_I^t = |I|\,\bar{\vm{N}}_{\Gamma,I}^\transpose\widehat{\vm{t}}_{N,I},
\end{equation}
where the nodal components are now computed with respect to the one-dimensional domain
on the Neumann boundary; that is, $|I|=\sum_{e\in\mathcal{T}_I}\frac{1}{2}|e|$, where $e$ is
an element's edge located on the Neumann boundary and $|e|$ its length;
$\mathcal{T}_I$ now represents the set of edges connected to node $I$ on the Neumann boundary.
Proceeding as mentioned, the remainder nodal matrices are
\begin{equation}\label{eq:nodal_N_gamma}
\bar{\vm{N}}_{\Gamma,I}=\frac{1}{|I|}\sum_{e\in\mathcal{T}_I}|e|\,\frac{1}{2}[\,\bar{\vm{N}}_\Gamma\,]_{{}_e},
\end{equation}
and
\begin{equation}\label{eq:nodal_traction}
\widehat{\vm{t}}_{N,I}=\frac{1}{|I|}\sum_{e\in\mathcal{T}_I}|e|\,\frac{1}{2}[\,\widehat{\vm{t}}_N\,]_{{}_e}.
\end{equation}

\section{Stabilization for the node-based uniform strain virtual element method}
\label{sec:stabilization}

Within the VEM framework, stabilization is one of the key ingredients to guarantee 
convergence. However, in nodal integration, stabilization can make the formulation 
somewhat stiff in the nearly incompressible limit~\cite{Puso-Solberg:2006}. To mitigate
this, a modified constitutive matrix $\widetilde{\vm{D}}$ can be used in lieu of the standard
$\vm{D}$ when computing the stability matrix. We opt for $\widetilde{\vm{D}}$
given as follows~\cite{Puso-Chen-Zywicz-Elmer:2008}:

\begin{equation}\label{eq:Dtilde}
\widetilde{\vm{D}}=\vm{D}\big(\widetilde{E},\widetilde{\nu}\big),
\end{equation}	
where 
\begin{equation}\label{eq:modyoungnu}
\widetilde{E}=\frac{\widetilde{\nu}\left(3\widetilde{\lambda}+2\widetilde{\mu}\right)}{\widetilde{\lambda}+\widetilde{\mu}},
\quad
\widetilde{\nu}=\frac{\widetilde{\lambda}}{2\left(\widetilde{\lambda}+\widetilde{\mu}\right)};		
\end{equation}
$\widetilde{\lambda}$ and $\widetilde{\mu}$ are, respectively, the modified Lam\'e's first and second parameters, which
are calculated as follows: 
\begin{equation}\label{eq:modlame}
\widetilde{\mu}\coloneqq\mu, \quad \widetilde{\lambda}\coloneqq\min{\left(\lambda,25\widetilde{\mu}\right)},
\end{equation}	
where	$\lambda$ and $\mu$ are the Lam\'e parameters of the problem to solve. Similarly,
another possibility to mitigate the stiff response is offered by the use of a modified constitutive
matrix $\vm{D}_{\mu}$ given by
\begin{equation}\label{eq:Dmu1}
\vm{D}_{\mu} =
\smat{2\mu & 0 & 0 \\ 0 & 2\mu & 0 \\ 0 & 0 & \mu},
\end{equation}
which omits the pressure parameter $\lambda(\varepsilon_{11}+\varepsilon_{22})$ 
that is responsible of the possible stiff behavior. Alternatively,
a simpler definition of $\vm{D}_{\mu}$ that produces practically identical results 
is obtained using the constitutive matrix that is related to the shear deformations; i.e.,
\begin{equation}\label{eq:Dmu2}
\vm{D}_{\mu} =
\smat{0 & 0 & 0 \\ 0 & 0 & 0 \\ 0 & 0 & \mu}.
\end{equation}

We point out that the foregoing stabilization procedure with $\vm{D}_{\mu}$ as given 
in \eref{eq:Dmu2} is similar to one of the stabilization approaches in the B-bar 
VEM~\cite{Park-Chi-Paulino:2021} (see Appendix A therein).

Using the modified constitutive matrices, we define the stability matrix as the nodal 
version of the \textit{D-recipe} stabilization technique~\cite{BeiraoDaVeiga-Dassi-Russo:2017,Mascotto:2018}. 
This gives two diagonal stability matrices whose entries are given by
\begin{equation}\label{eq:stabDtilde}
\left(\vm{S}_I\right)_{i,i}=
\max{\left(1,\big(|I|\,\vm{B}_I^\transpose\,\widetilde{\vm{D}}\,\vm{B}_I\big)_{i,i}\right)}
\end{equation}
and
\begin{equation}\label{eq:stabDmu}
\left(\vm{S}_I\right)_{i,i}=
\max{\left(1,\big(|I|\,\vm{B}_I^\transpose\,\vm{D}_{\mu}\,\vm{B}_I\big)_{i,i}\right)}.
\end{equation}
We refer to \eref{eq:stabDtilde} and \eref{eq:stabDmu}, respectively, 
as $\widetilde{\vm{D}}$ and $\vm{D}_{\mu}$ stabilizations.

\section{Numerical experiments}
\label{sec:numexp}

In this section, numerical experiments are conducted to demonstrate the accuracy,
convergence and stability of the NVEM for compressible and nearly incompressible 
\alejandro{linear elastostatics. 
In addition, the performance of the NVEM in linear elastodynamics is also demonstrated.
For the elastostatic numerical experiments, the numerical results are compared 
with the B-bar VEM~\cite{Park-Chi-Paulino:2021} and the standard linear VEM.
For the elastodynamic numerical experiments, the numerical results are compared
with the finite element method (FEM) using eight-node quadrilateral elements (Q8).

To assess the convergence of the numerical solution of the NVEM approach,
we use the nodal values of the exact displacement, strain and 
pressure fields, that is, $\vm{u}(\vm{x}_I)$, $\bsym{\varepsilon}(\vm{x}_I)$ and $p(\vm{x}_I)$,
respectively, and define the following nodal errors at node $I$:
\begin{equation*}
\vm{e}_u=\vm{u}(\vm{x}_I)-\vm{u}_I, \quad \vm{e}_{\varepsilon}=\bsym{\varepsilon}(\vm{x}_I)-\widehat{\bsym{\varepsilon}}_I
\quad e_p=p(\vm{x}_I)-p_I,
\end{equation*}
where $\vm{u}_I$, $\widehat{\bsym{\varepsilon}}_I$ and $p_I=-\lambda\,\trace{\widehat{\bsym{\varepsilon}}_I}$ are the 
nodal values of the displacement, strain and pressure fields approximations, respectively. 
Using the nodal errors, the following
norms are evaluated through nodal integration in the NVEM approach:
\begin{equation*}\label{eq:nodal_norms_L2u}
\frac{\|\vm{u}(\vm{x})-\vm{u}_h(\vm{x})\|_{L^2(\Omega)}}{\|\vm{u}(\vm{x})\|_{L^2(\Omega)}}
=\sqrt{\frac{\sum_{I=1}^N \vm{e}_u^\transpose \vm{e}_u\,|I|}{\sum_{I=1}^N \vm{u}^\transpose(\vm{x}_I)\,\vm{u}(\vm{x}_I) \,|I|}},
\end{equation*}
\begin{equation*}\label{eq:nodal_norms_H1}
\frac{\|\vm{u}(\vm{x})-\vm{u}_h(\vm{x})\|_{H^1(\Omega)}}{\|\vm{u}(\vm{x})\|_{H^1(\Omega)}}
=\sqrt{\frac{\sum_{I=1}^N \vm{e}_{\varepsilon}^\transpose \,\vm{D}\,\vm{e}_{\varepsilon}\,|I|}{\sum_{I=1}^N  \bsym{\varepsilon}^\transpose(\vm{x}_I) \,\vm{D}\, \bsym{\varepsilon}(\vm{x}_I) \,|I|}},
\end{equation*}
\begin{equation*}\label{eq:nodal_norms_L2p}
\frac{\|p(\vm{x})-p_h(\vm{x})\|_{L^2(\Omega)}}{\|p(\vm{x})\|_{L^2(\Omega)}}
=\sqrt{\frac{\sum_{I=1}^N e_p^\transpose e_p\,|I|}{\sum_{I=1}^N p^\transpose(\vm{x}_I)\,p(\vm{x}_I) \,|I|}}.
\end{equation*}
For the standard VEM and the B-Bar VEM approaches, the above norms are evaluated using
cubature inside the element through partition of the element into triangles
and the errors defined at the element level as $\vm{e}_u=\vm{u}(\vm{x})-\Pi\vm{u}_h(\vm{x})$,
$\vm{e}_{\varepsilon}=\bsym{\varepsilon}(\vm{x})-\widehat{\bsym{\varepsilon}}_h(\vm{u}_h)$ and
$e_p=p(\vm{x})-p_h(\vm{x})$, where $p_h(\vm{x})=-\lambda\,\trace{\widehat{\bsym{\varepsilon}}_h(\vm{u}_h)}$.
}

\subsection{Patch test and numerical stability}
\begin{figure}[!bth]
        \centering
        \mbox{
        \subfigure[] {\label{fig:patchtestsmesh_a}\includegraphics[width=0.45\linewidth]{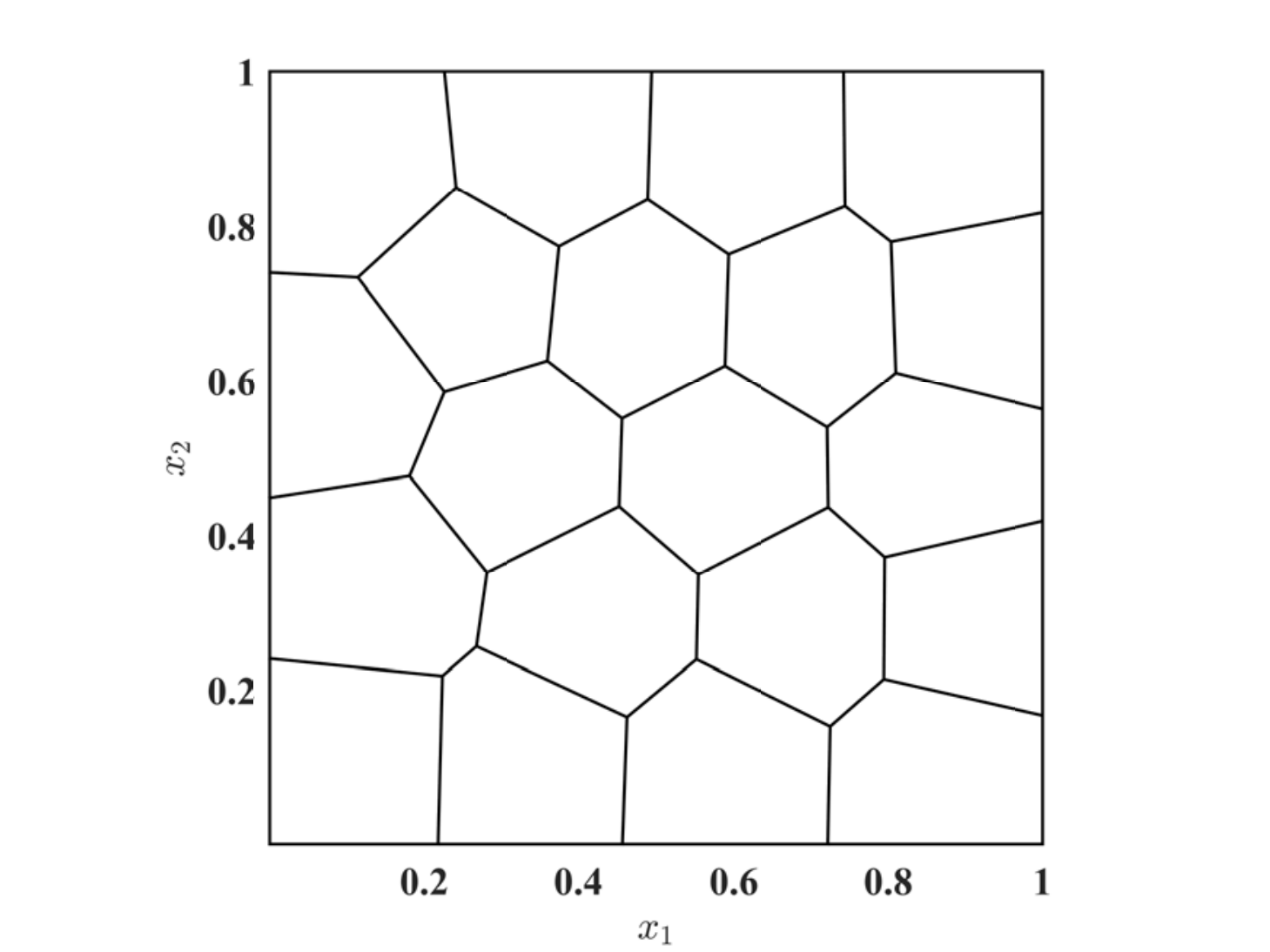}}
        \subfigure[] {\label{fig:patchtestsmesh_b}\includegraphics[width=0.45\linewidth]{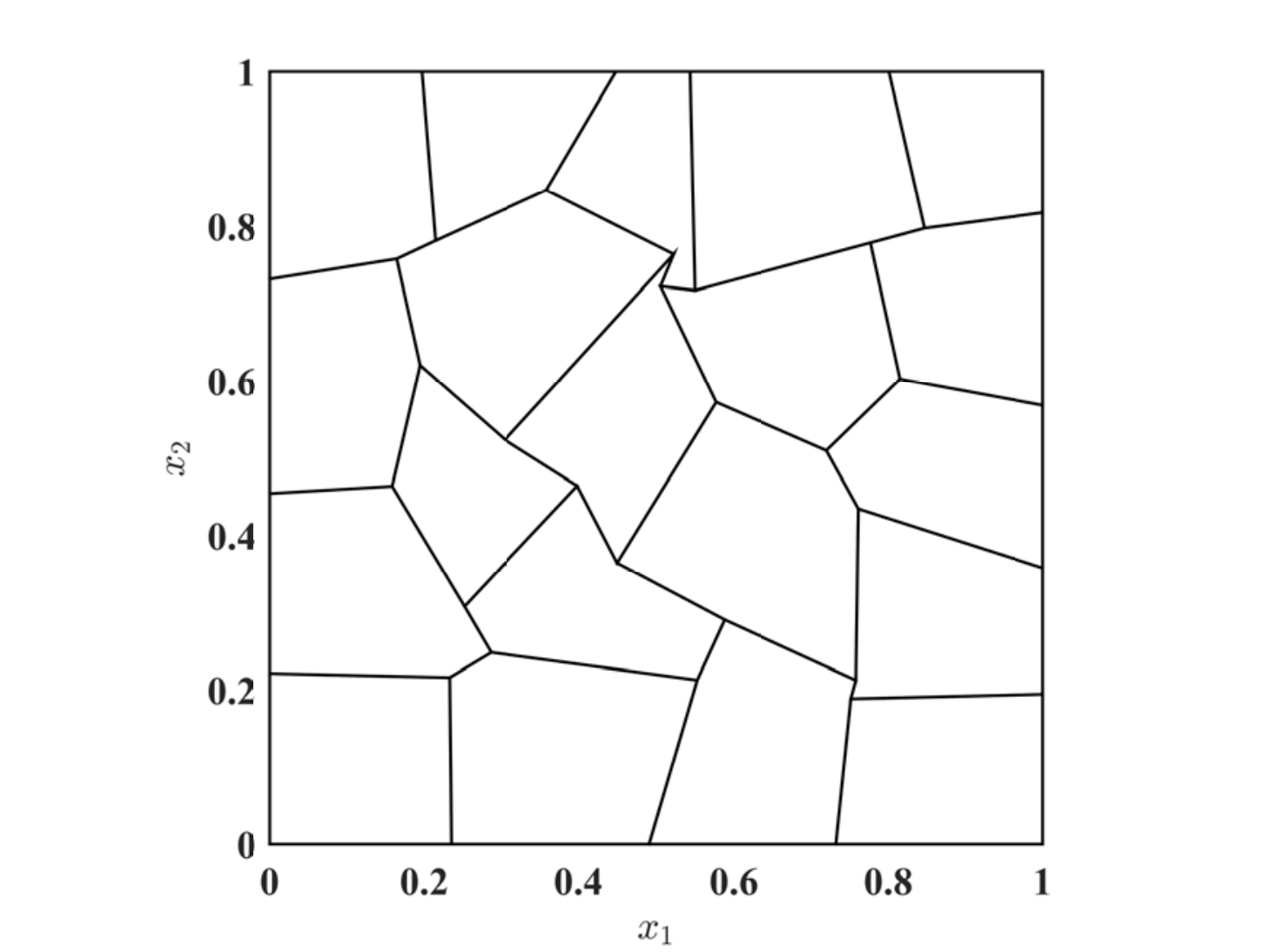}}
        }     
        \caption{Meshes used for the displacement patch tests and numerical stability. 
                 (a) Regular mesh and (b) distorted mesh.}
        \label{fig:patchtestsmesh}
\end{figure}

\alejandrog{With dimensions in inches,} \alejandro{we investigate the accuracy 
and numerical stability of the proposed NVEM approach by solving displacement 
patch tests and eigenvalue analyses} \alejandrog{on the meshes shown in \fref{fig:patchtestsmesh}
for the square domain $\Omega=(0,1)^2$.}
\alejandro{In both problems $\vm{b}=\vm{0}$ \alejandrog{lbf/$\mathrm{in}^3$} and plane strain condition is assumed with the 
following material parameters: $E=1\times 10^7$ psi and $\nu = 0.3$. For the displacement 
patch test, the Dirichlet boundary condition
$\vm{u}_D=\left\{x_1 \quad x_1+x_2\right\}^\transpose$ in} \alejandrog{inches} \alejandro{is prescribed along the 
entire boundary of the domain. Numerical results for the relative error 
in the $L^2$ norm and the $H^1$ seminorm are presented in
Tables~\ref{table:L2norm_patchtest_elast} and \ref{table:H1norm_patchtest_elast}, respectively,
where it is confirmed that the patch test is met to machine precision. 

The eigenvalue analyses delivers three zero eigenvalues, 
which corresponds to the three zero-energy rigid
body modes. The three mode shapes that follow the three rigid body modes are
depicted in \fref{fig:numstability} for the distorted mesh. These mode shapes are smooth, which is an
indication of the numerical stability of the method.
}

\begin{table}[!tbhp]
\setlength{\arrayrulewidth}{.15em}
\newcolumntype{a}{>{\columncolor{white}}c}
\begin{center}
\addtolength{\tabcolsep}{-1pt}
\caption{Relative error in the $L^{2}$ norm for the displacement patch tests.}
\vspace*{0pt}
\renewcommand{\arraystretch}{1.5}
\begin{tabular}{aaa}
\hline
\rowcolor{white}
Method & Regular & Distorted \\
\hline
\rowcolor{white}
NVEM, $\widetilde{\vm{D}}$ & $7.8 \times 10^{-16}$  & $9.7 \times 10^{-16}$  \\
\rowcolor{white}
NVEM, $\vm{D}_\mu$ & $4.8 \times 10^{-16}$  & $3.9 \times 10^{-16}$  \\
\hline
\end{tabular}
\label{table:L2norm_patchtest_elast}
\end{center}
\end{table}

\begin{table}[!tbhp]
\setlength{\arrayrulewidth}{.15em}
\newcolumntype{a}{>{\columncolor{white}}c}
\begin{center}
\addtolength{\tabcolsep}{-1pt}
\caption{Relative error in the $H^{1}$ seminorm for the displacement patch tests.}
\vspace*{0pt}
\renewcommand{\arraystretch}{1.5}
\begin{tabular}{aaa}
\hline
\rowcolor{white}
Method & Regular & Distorted \\
\hline
\rowcolor{white}
NVEM, $\widetilde{\vm{D}}$ & $5.9 \times 10^{-16}$  & $8.2 \times 10^{-16}$  \\
\rowcolor{white}
NVEM, $\vm{D}_\mu$ & $8.0 \times 10^{-16}$  & $6.6 \times 10^{-16}$  \\
\hline
\end{tabular}
\label{table:H1norm_patchtest_elast}
\end{center}
\end{table}

\begin{figure}[!bth]
        \centering
        \mbox{
        \subfigure[] {\label{fig:numstability_a}\includegraphics[width=0.3\linewidth]{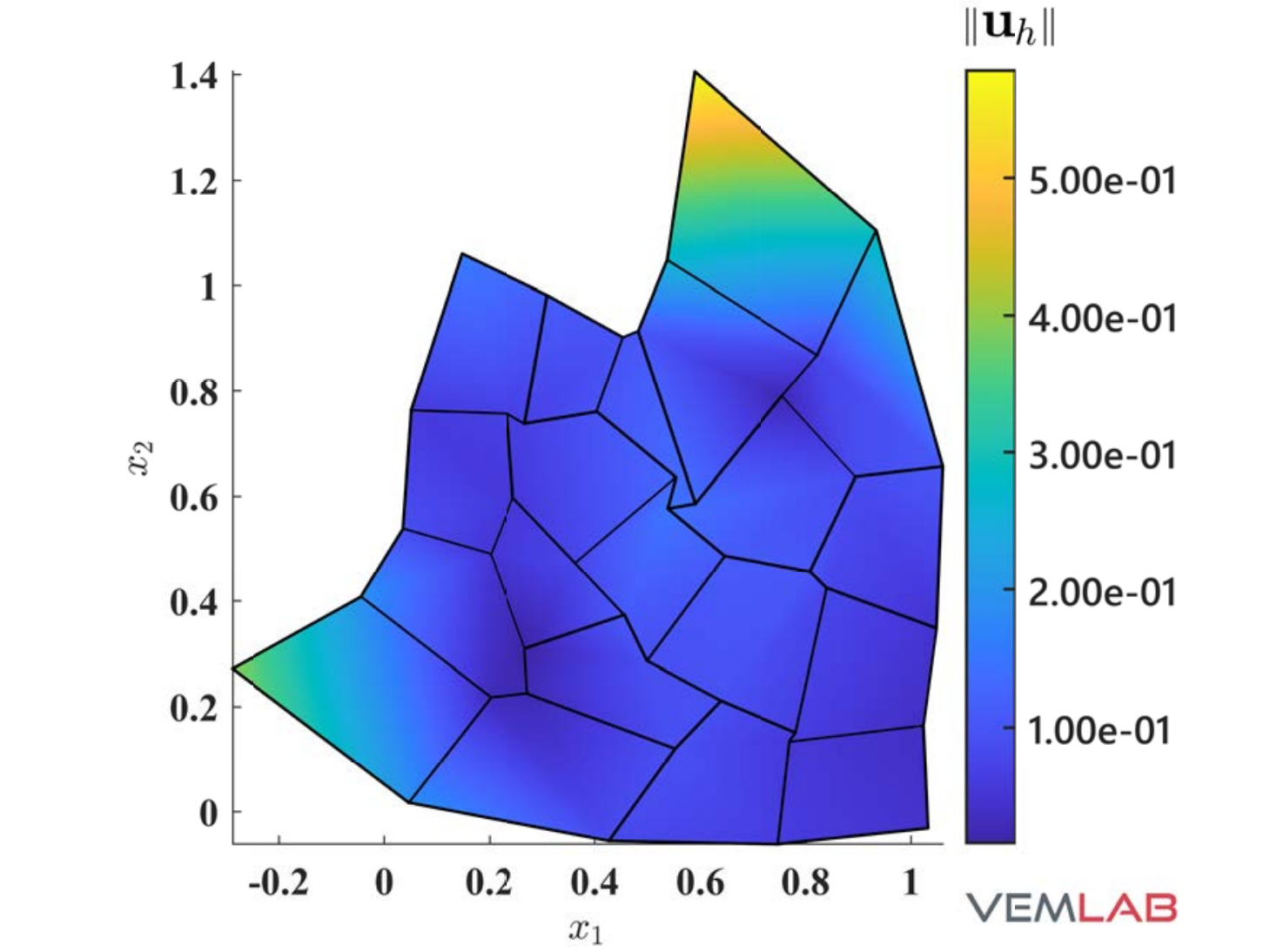}}
        \subfigure[] {\label{fig:numstability_b}\includegraphics[width=0.3\linewidth]{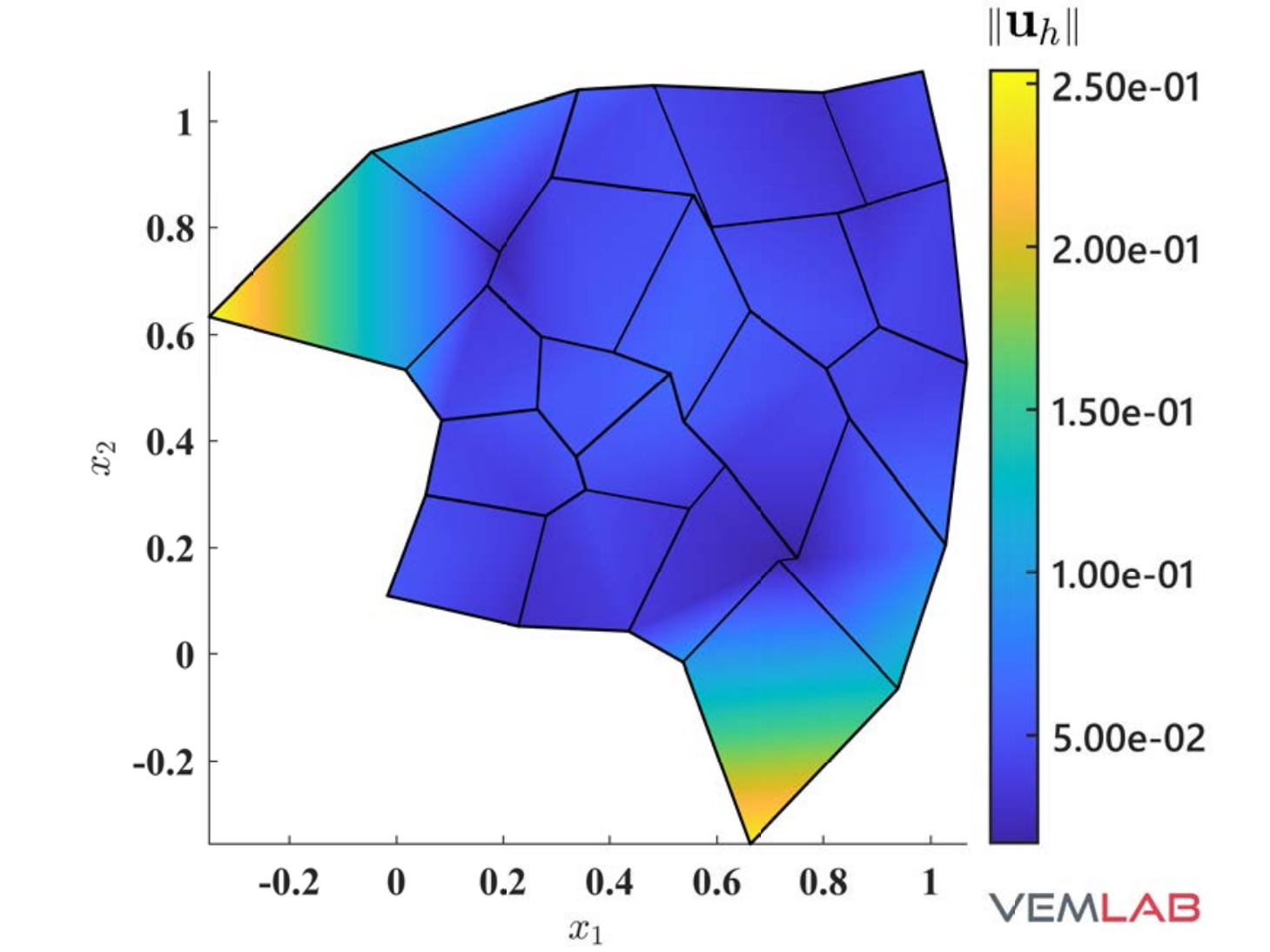}}
        \subfigure[] {\label{fig:numstability_c}\includegraphics[width=0.3\linewidth]{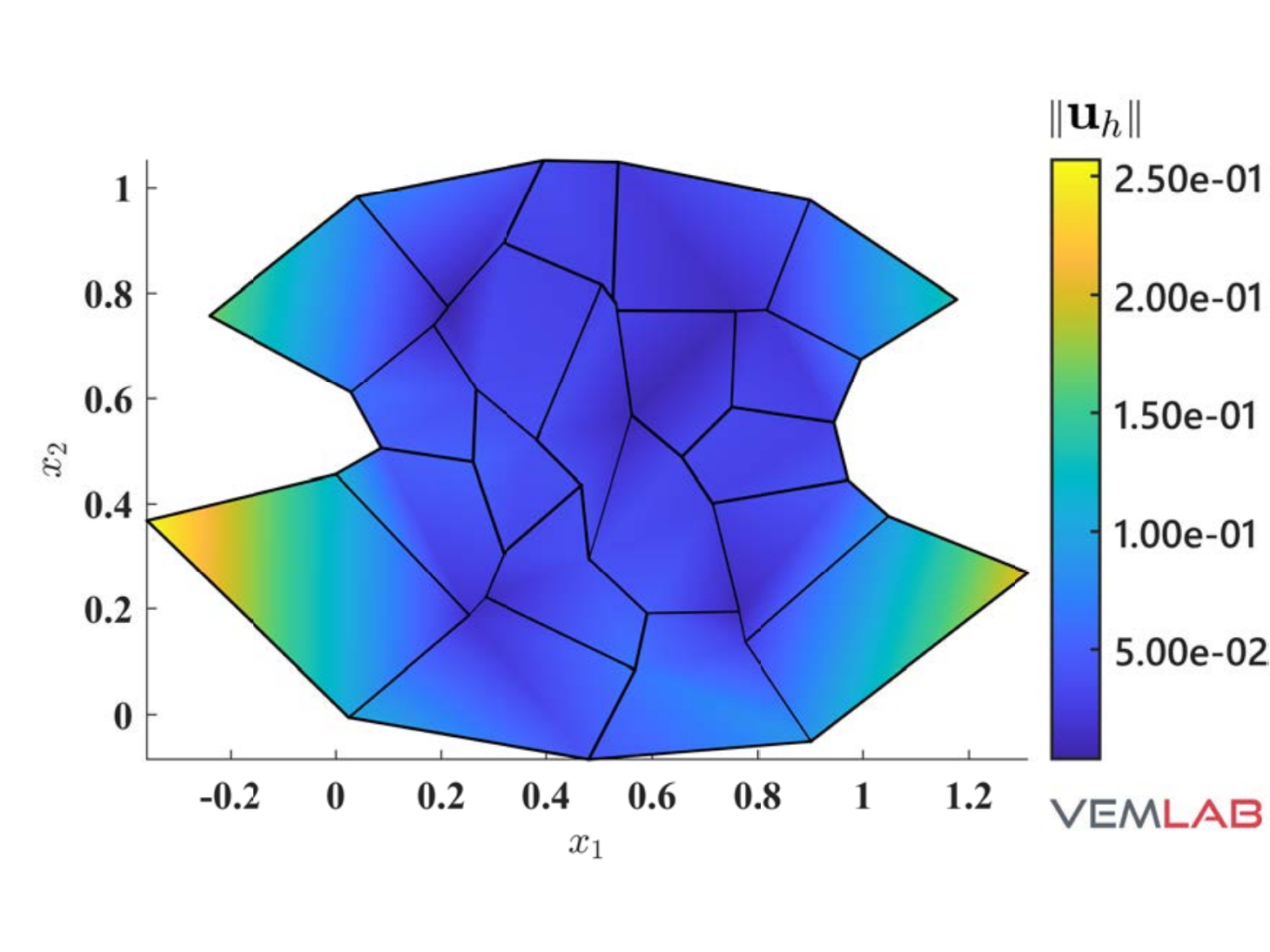}}
        } 
        \mbox{
        \subfigure[] {\label{fig:numstability_d}\includegraphics[width=0.3\linewidth]{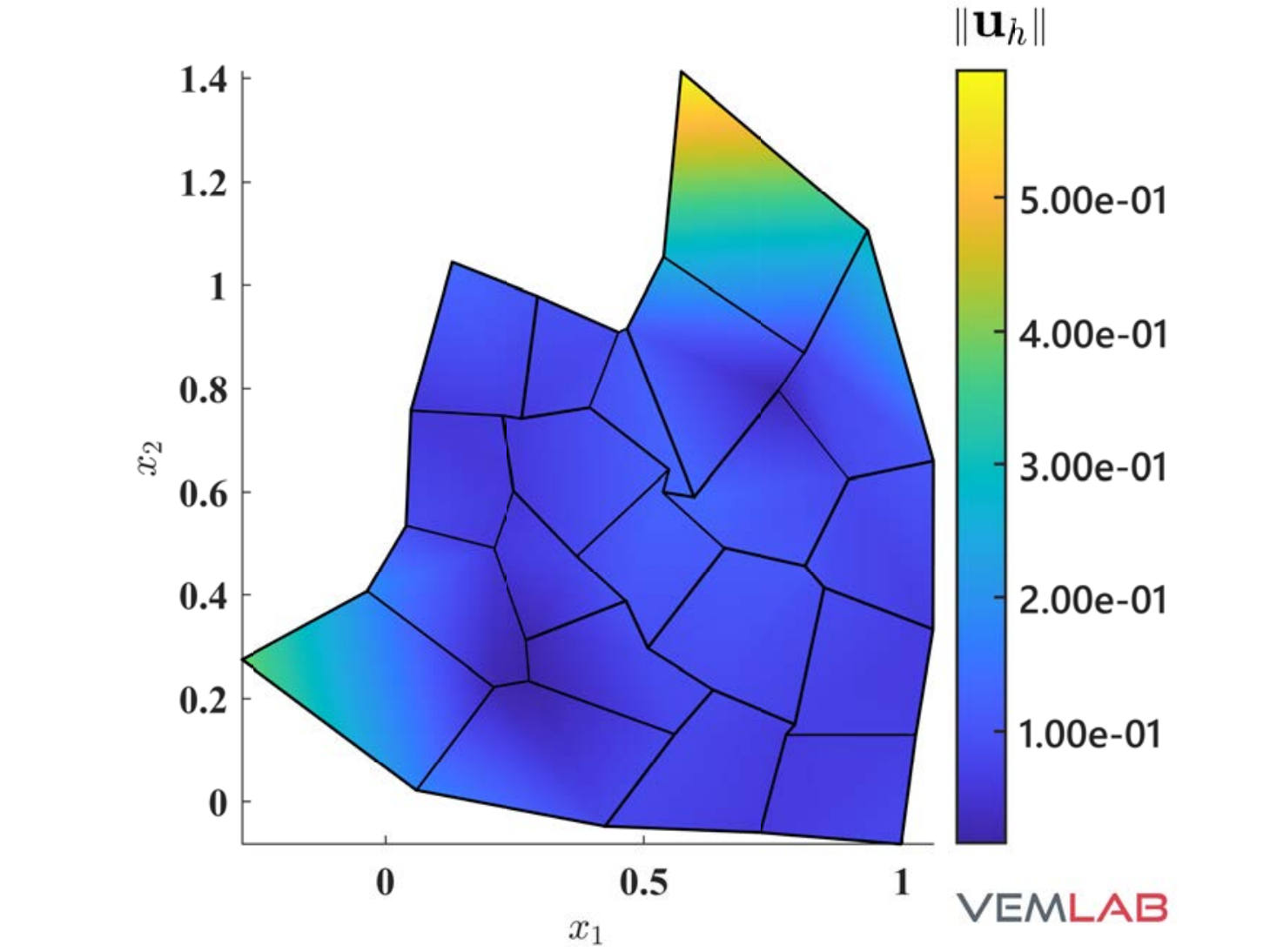}}
        \subfigure[] {\label{fig:numstability_e}\includegraphics[width=0.3\linewidth]{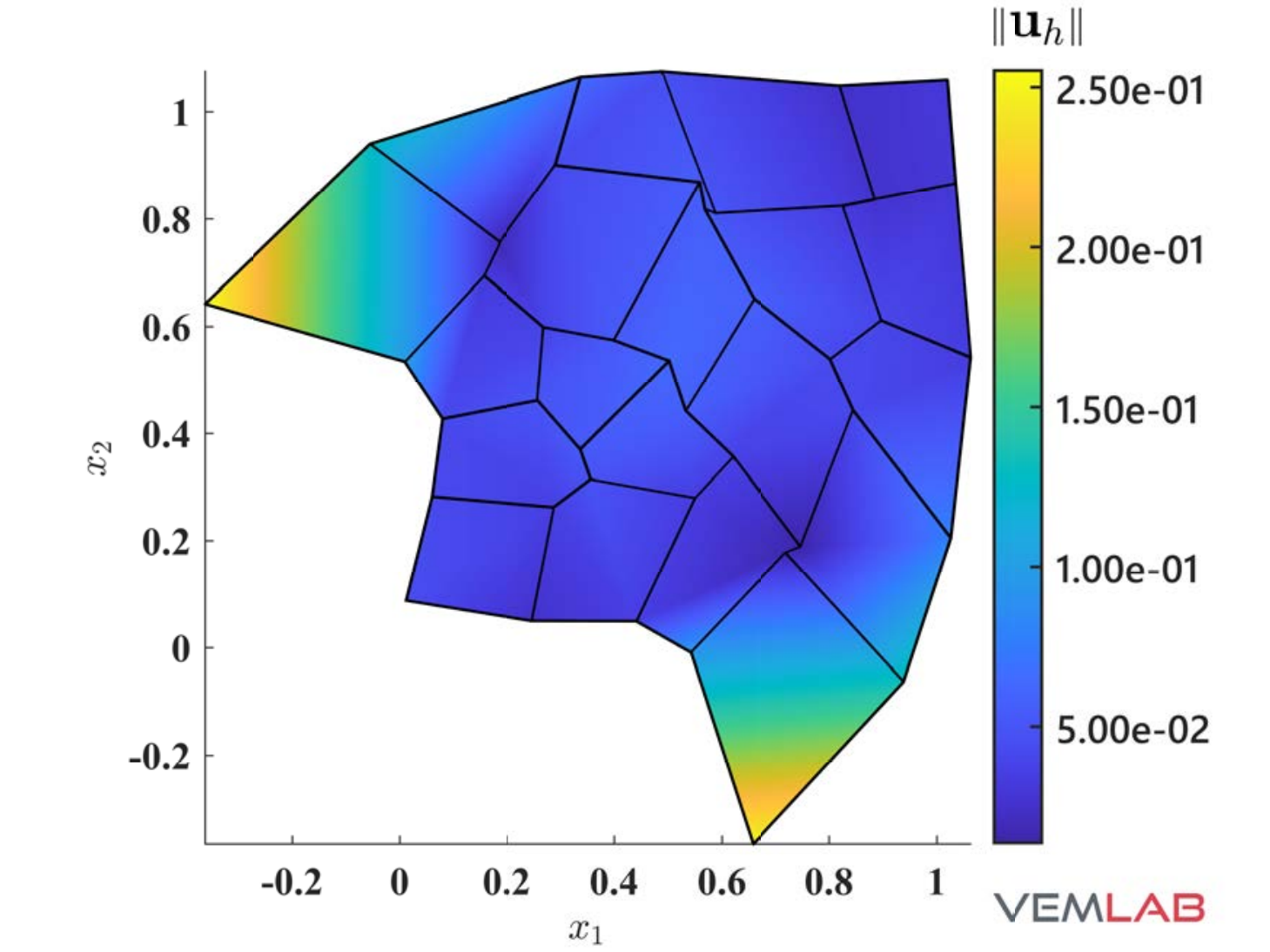}}
        \subfigure[] {\label{fig:numstability_f}\includegraphics[width=0.3\linewidth]{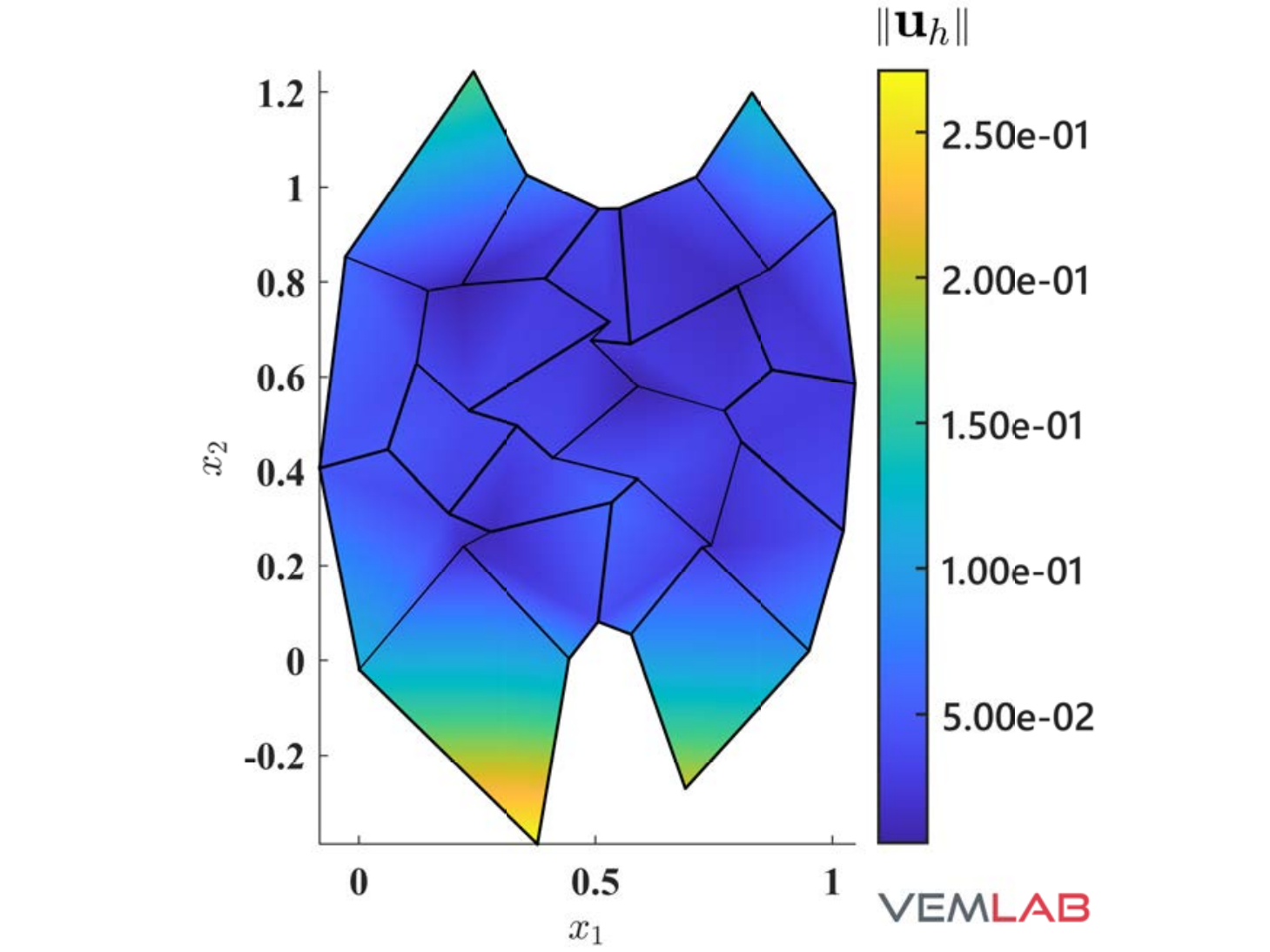}}
        }             
        \caption{Eigenvalue analyses on the distorted mesh. A pictorial of the three mode shapes that follow the three rigid body modes. 
                 (a)--(c) NVEM with $\widetilde{\vm{D}}$ stabilization and (d)--(f) NVEM with $\vm{D}_\mu$ stabilization.}
        \label{fig:numstability}
\end{figure}

\subsection{Colliding flow}
\label{sec:colliding}
\begin{figure}[!bth]
	\centering
	\includegraphics[width=0.6\linewidth]{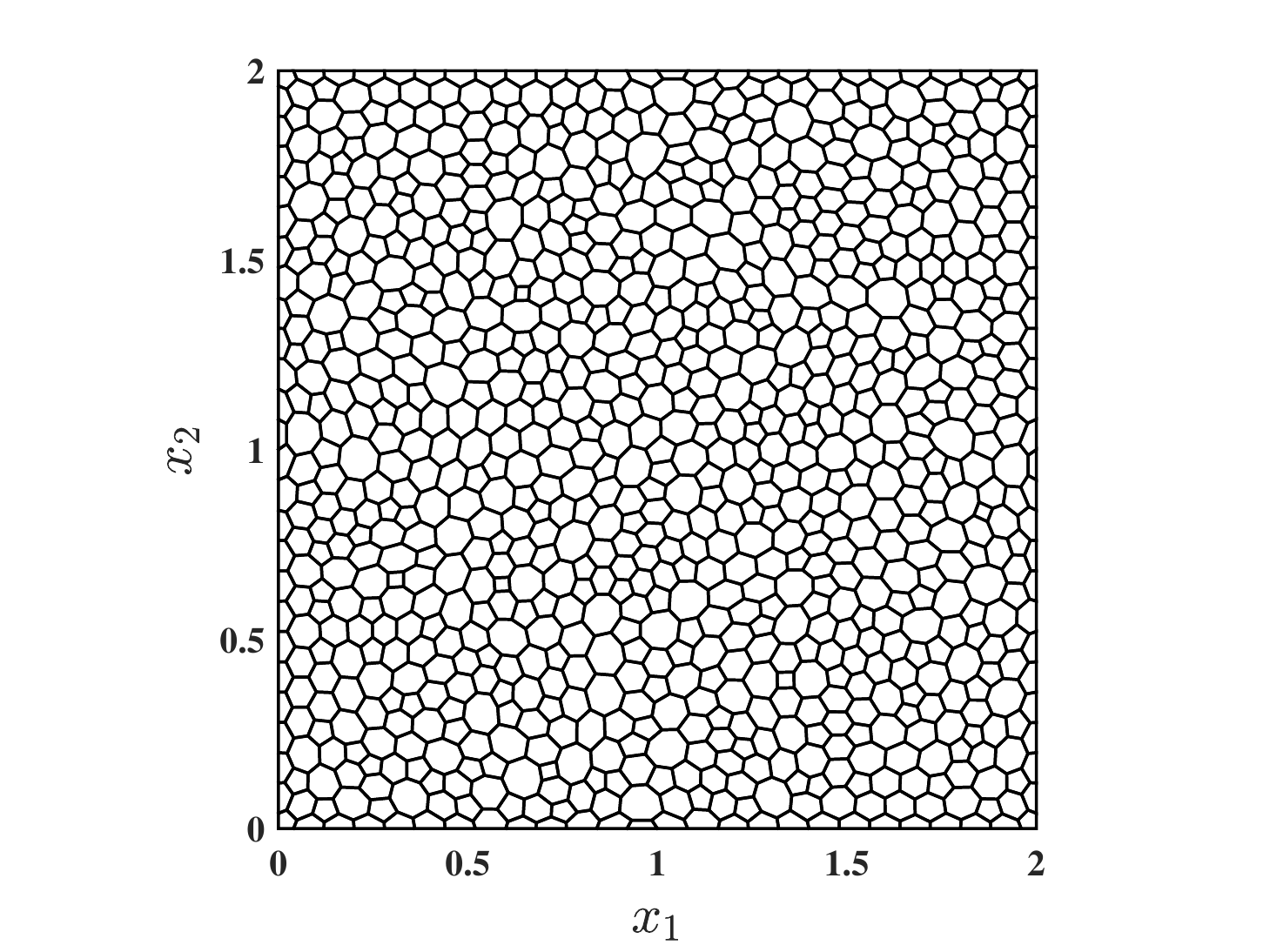}	
	\caption{Sample mesh for the colliding flow problem.}
	\label{fig:collidingmesh}
\end{figure}

We consider a simple model of colliding flow, which is a well-known 
standard benchmark problem (for example, see
Refs.~\cite{Elman-Silvester-Wathen:2006}) \alejandro{that is described by
the Stokes equations. Since the Stokes flow
equations coincide with the equations that govern the static incompressible
elasticity, the colliding flow} can be solved using
the linear elastostatic model with the following constitutive matrix:
\begin{equation*}
\vm{D} =
\smat{\lambda+2\mu & \lambda & 0 \\ \lambda & \lambda+2\mu & 0 \\ 0 & 0 & \mu},
\end{equation*}
where the Lam\'e parameters are set to $\lambda=5\times 10^7$ \alejandro{psi}
and $\mu=1$ \alejandro{psi}, which
results in a Poisson's ratio $\nu=0.49999999$. \alejandro{For dimensions in inches}, 
the problem is defined on the square domain $\Omega=(0,2)^2$ and \alejandro{Dirichlet
boundary conditions are imposed along the entire boundary of the domain using}
the following analytical solution~\cite{Elman-Silvester-Wathen:2006}:
\begin{subequations}\label{colliding_soln}
\begin{align}
 u_1 &= 20(x_1-1)(x_2-1)^3 \label{colliding_soln1}\\
 u_2 &= 5(x_1-1)^4-5(x_2-1)^4 \label{colliding_soln2} \\
 p &= 60(x_1-1)^2(x_2-1)-20(x_2-1)^3 \label{colliding_soln3},
\end{align}
\end{subequations}
\alejandro{where $u_1$ and $u_2$ are the components of the velocity field
$\vm{u}$} \alejandrog{in inches,} \alejandro{and $p$ is the pressure field} \alejandrog{in psi}. 
The stabilizations presented in 
Section~\ref{sec:stabilization} are considered to assess the convergence of the method
upon mesh refinement. A sample polygonal mesh used in the convergence study 
is shown in~\fref{fig:collidingmesh}. 

\begin{figure}[!bth]
	\centering
	\mbox{
	\subfigure[] {\label{fig:collidingsolutions_a}\includegraphics[width=0.5\linewidth]{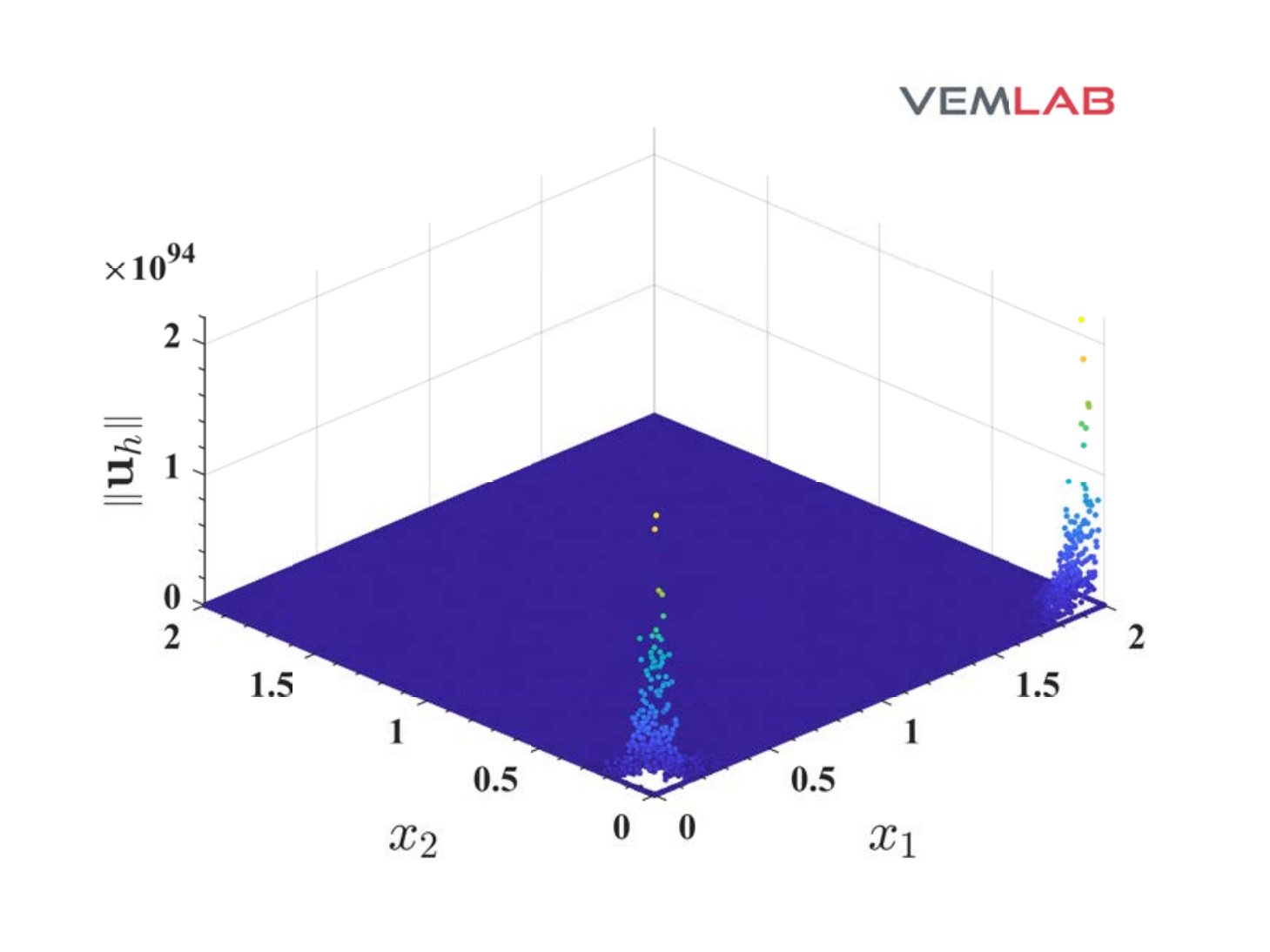}}
	\subfigure[] {\label{fig:collidingsolutions_b}\includegraphics[width=0.5\linewidth]{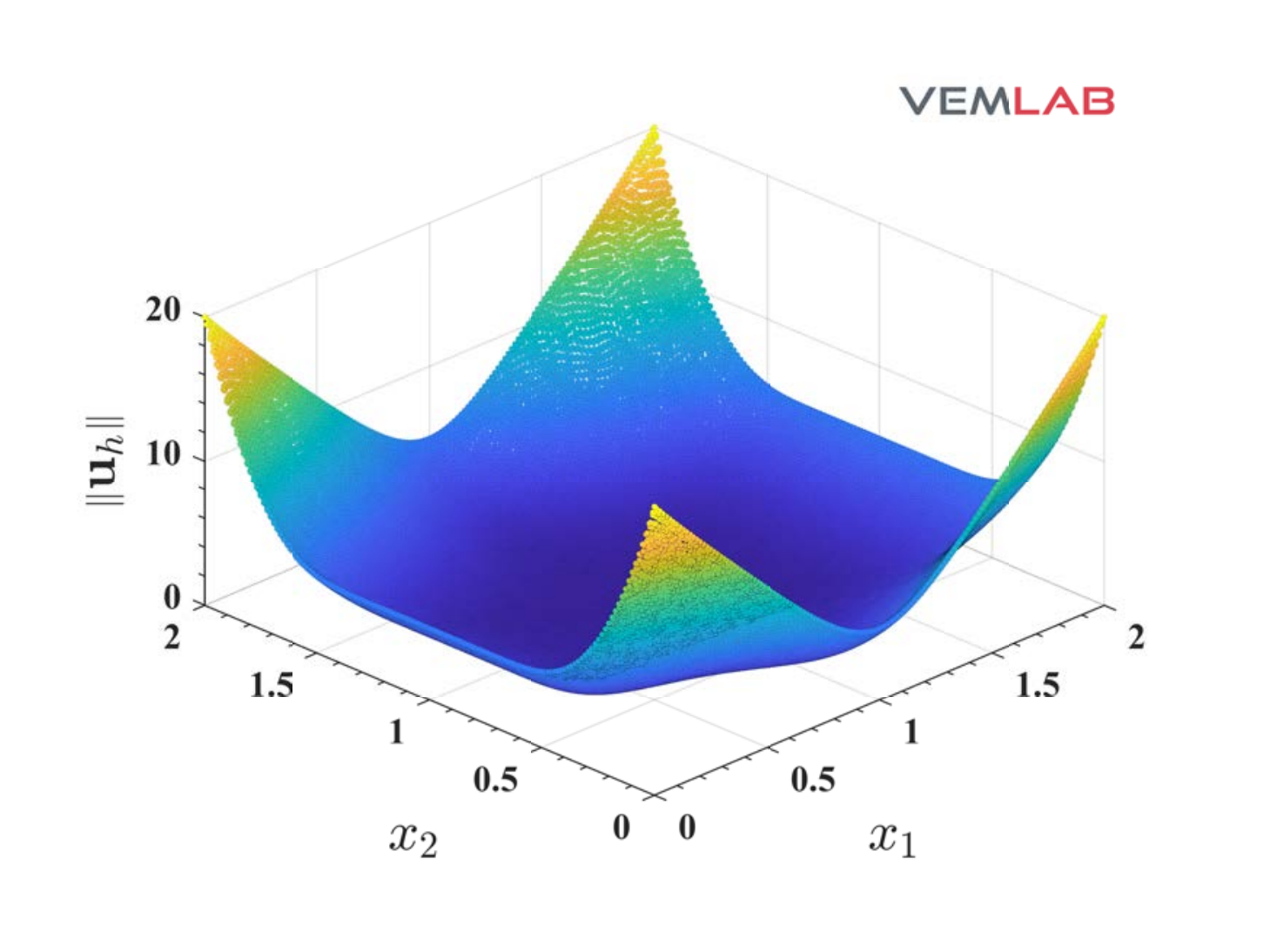}}
	}
	\mbox{
	\subfigure[] {\label{fig:collidingsolutions_c}\includegraphics[width=0.5\linewidth]{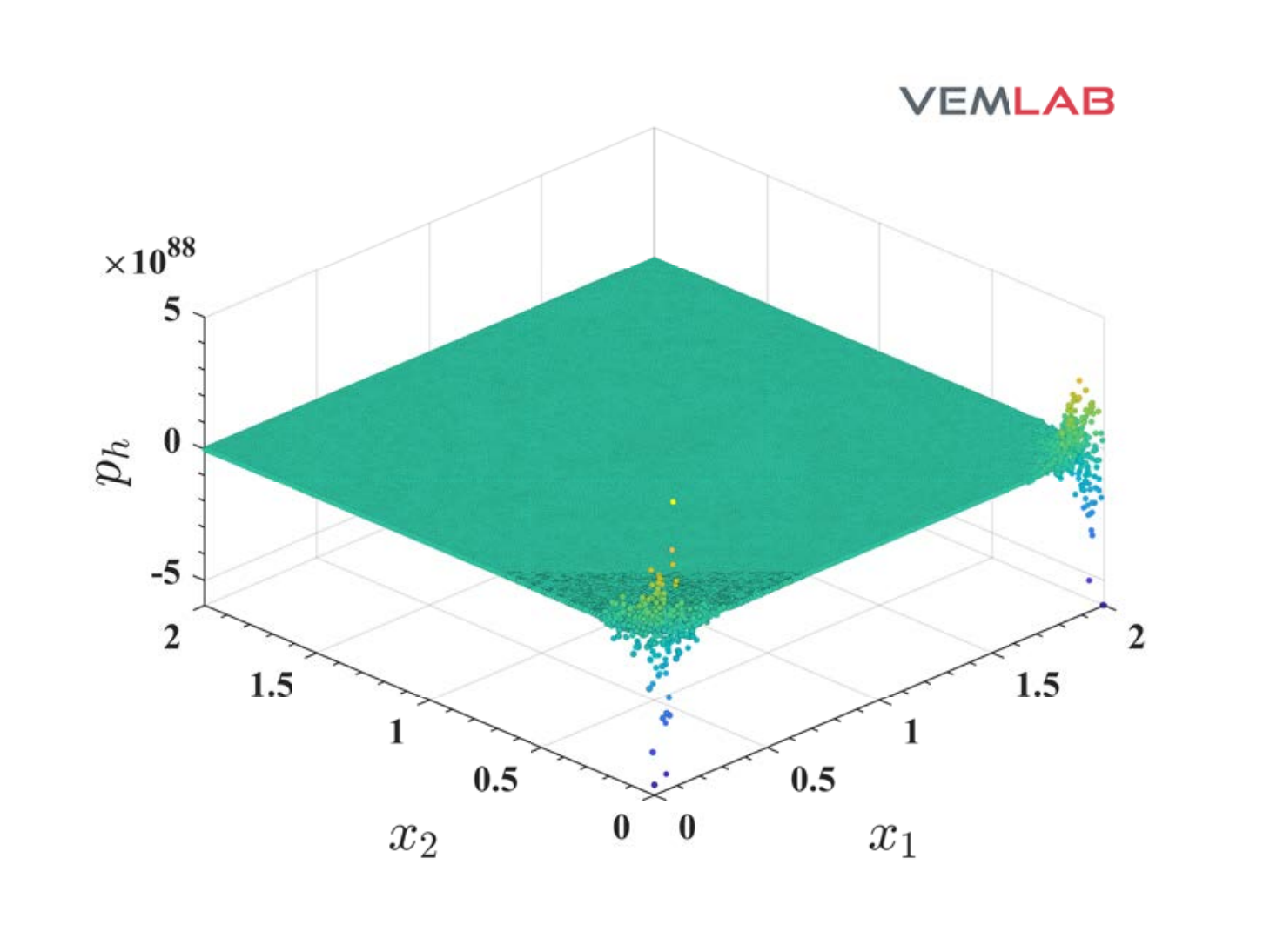}}
	\subfigure[] {\label{fig:collidingsolutions_d}\includegraphics[width=0.5\linewidth]{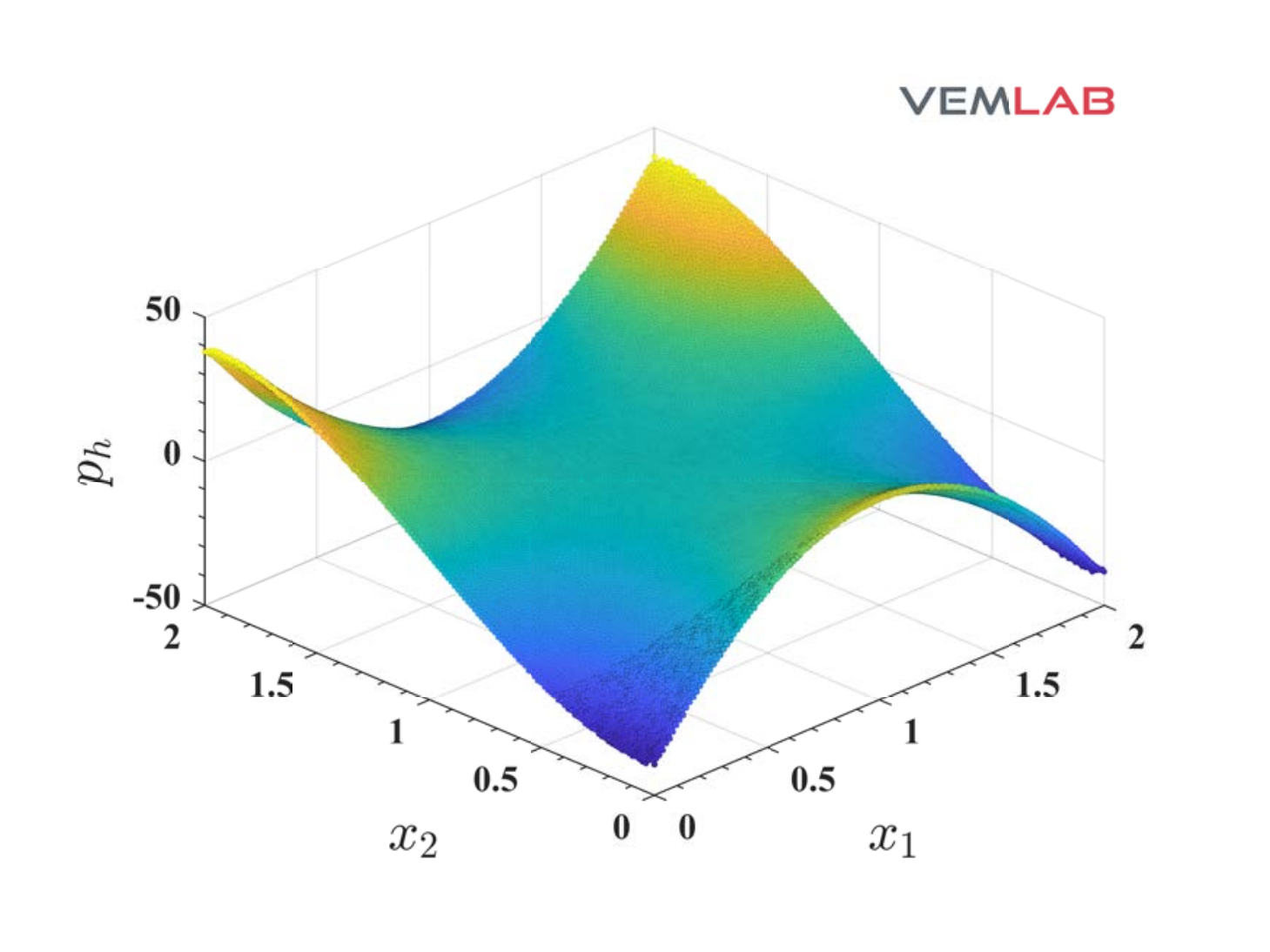}}	
	}	
	\caption{Comparison between unstabilized and stabilized solutions in the NVEM for 
	the colliding flow problem. (a) Unstabilized $\|\vm{u}_h\|$ \alejandro{(in)},  (b) stabilized $\|\vm{u}_h\|$ \alejandro{(in)}, 
	(c) unstabilized $\|p_h\|$ \alejandro{(psi)} and (d) stabilized $\|p_h\|$ \alejandro{(psi)}.}
	\label{fig:collidingsolutions}
\end{figure}

We start by showing the need for stability
in the NVEM approach. \fref{fig:collidingsolutions} presents a comparison between
unstabilized and stabilized solutions for the nodal fields $\|\vm{u}_h\|$
and $p_h$ --- scatter plots are used for the NVEM as in this approach 
the field variables are known at the nodes. The unstabilized solutions exhibit 
marked spurious oscillations at the corners $(0,0)$ and $(2,0)$ (Figs.~\ref{fig:collidingsolutions_a} 
and~\ref{fig:collidingsolutions_c}), whereas the stabilized solutions become smooth
on the whole domain (Figs.~\ref{fig:collidingsolutions_b} and~\ref{fig:collidingsolutions_d}).

\begin{figure}[!bth]
	\centering
	\mbox{
	\subfigure[] {\label{fig:collidingrates_a}\includegraphics[width=0.5\linewidth]{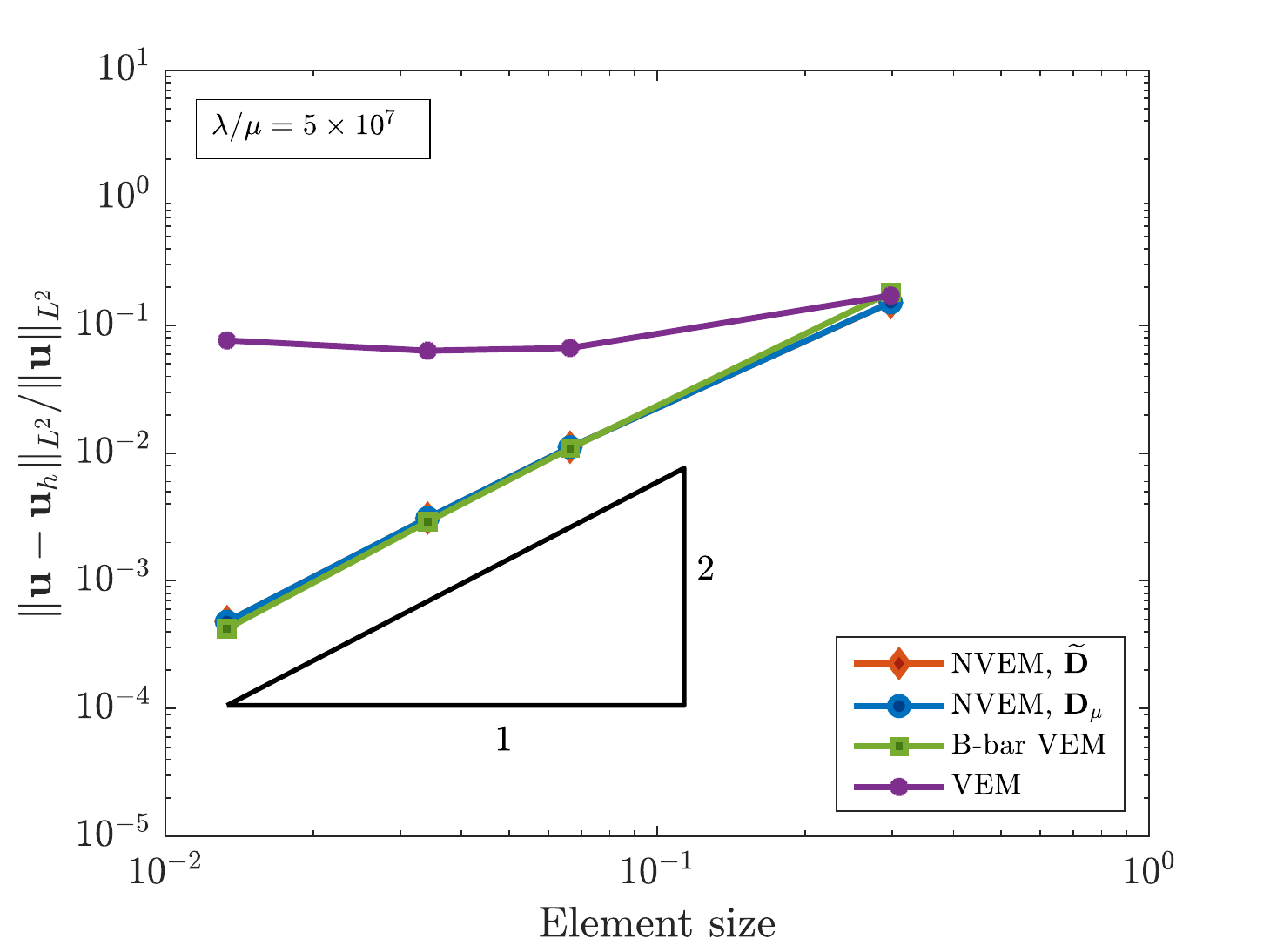}}
	\subfigure[] {\label{fig:collidingrates_b}\includegraphics[width=0.5\linewidth]{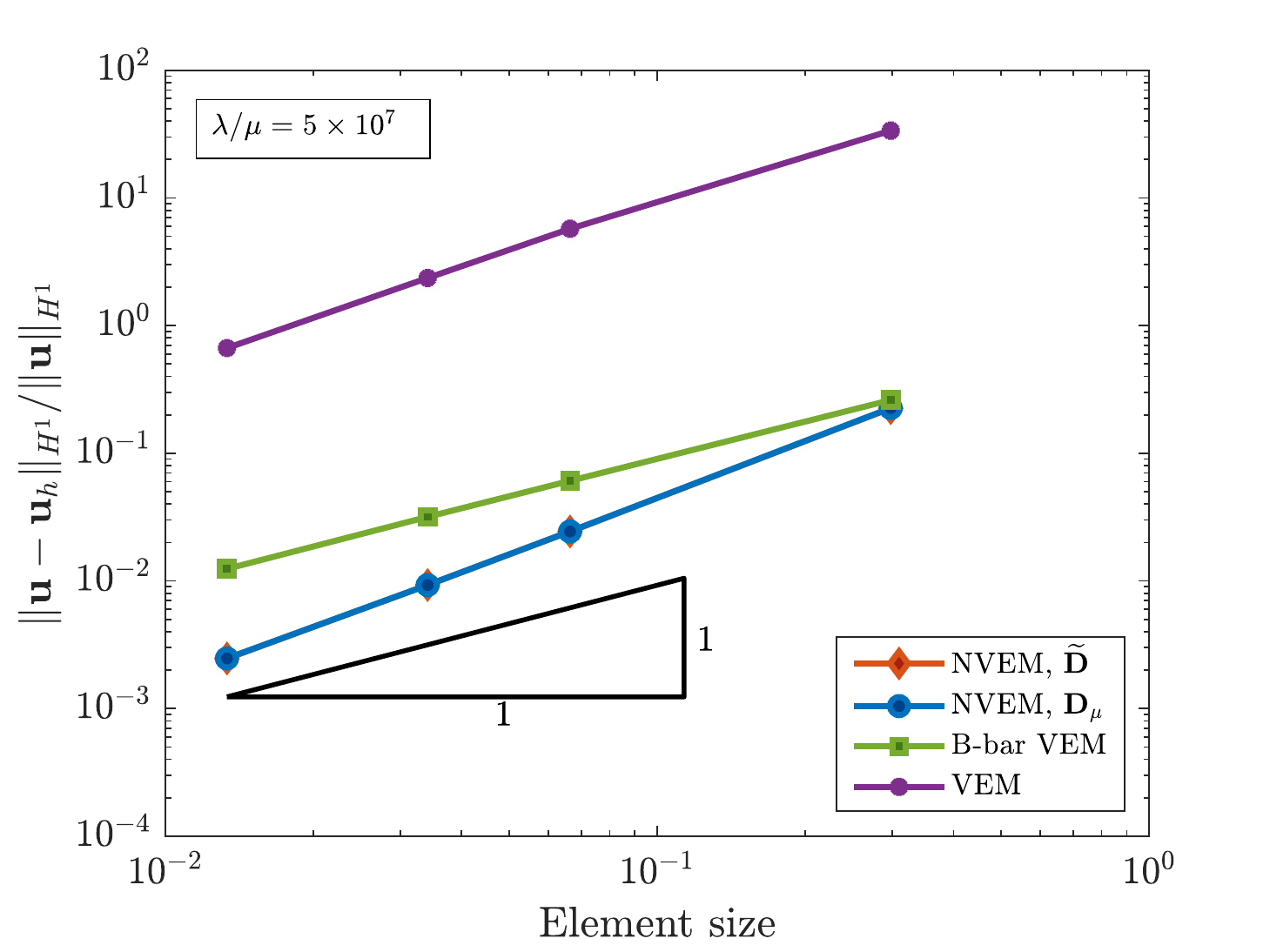}}
	}
	\subfigure[] {\label{fig:collidingrates_c}\includegraphics[width=0.5\linewidth]{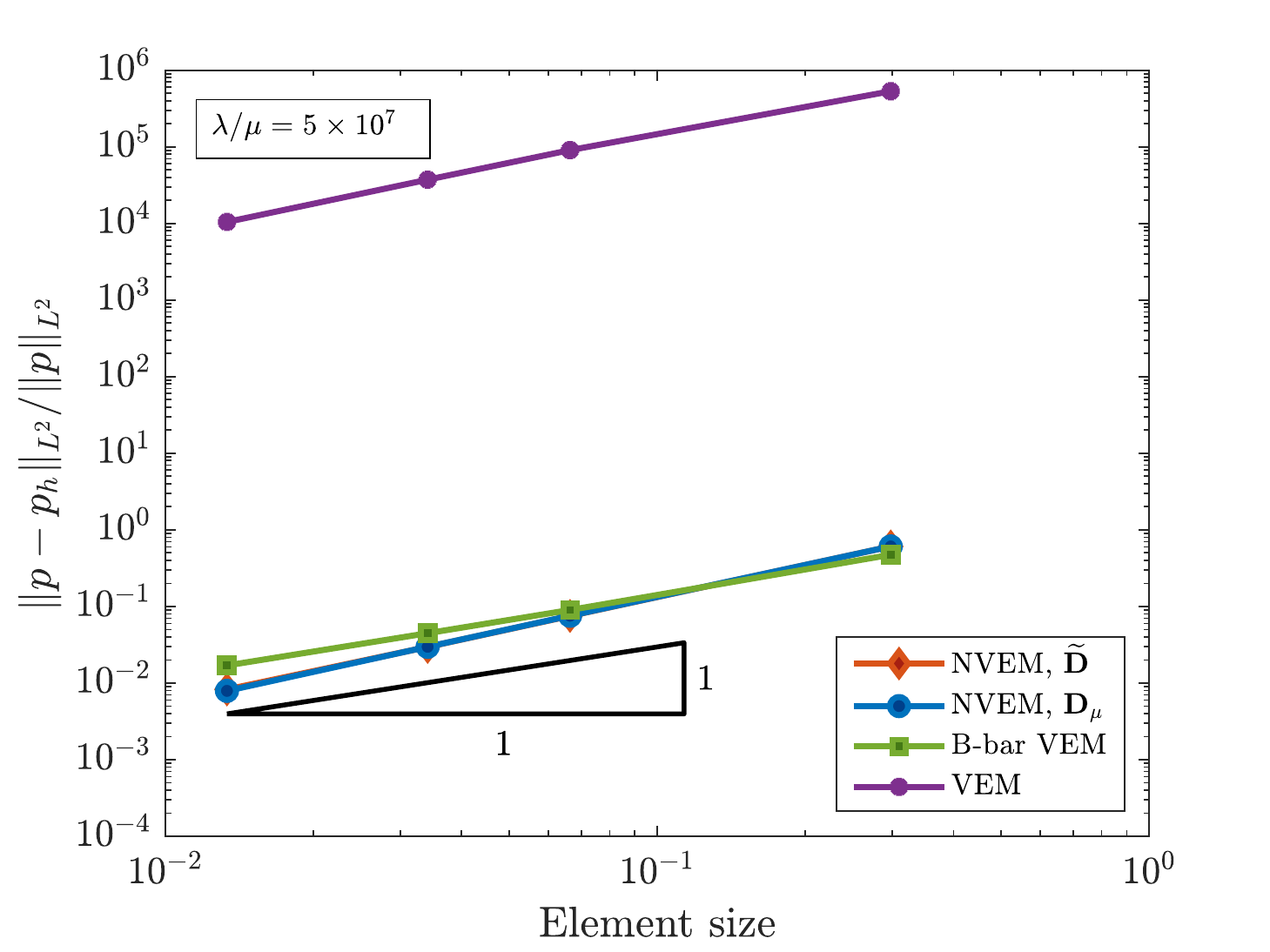}}	
	\caption{Colliding flow problem. Convergence rates in the (a) $L^2$ norm of the 
	displacement error, (b) $H^1$ seminorm of the displacement error and 
	(c) $L^2$ norm of the pressure error for the VEM, B-bar VEM and NVEM.}
	\label{fig:collidingrates}
\end{figure}

The convergence and accuracy of the NVEM with $\widetilde{\vm{D}}$
and $\vm{D}_\mu$ stabilizations are demonstrated in~\fref{fig:collidingrates},
where the $L^2$ norm and the $H^1$ seminorm of the displacement error, and
the $L^2$ norm of the pressure error indicate accurate solutions with 
optimal convergence rates of 2, 1 and 1, respectively. 
In comparison with the B-bar VEM approach, which is also accurate and optimally
convergent, the NVEM approach is slightly more accurate
in the $H^1$ seminorm of the displacement error and the $L^2$ norm of the 
pressure error. Finally, as expected, in these convergence plots the standard 
VEM solution behaves very inaccurate due to volumetric locking.

\subsection{Cantilever beam subjected to a parabolic end load}
\label{sec:beam}

\begin{figure}[!bth]
	\centering
	\mbox{
	\subfigure[] {\label{fig:beamproblem_a}\includegraphics[width=0.5\linewidth]{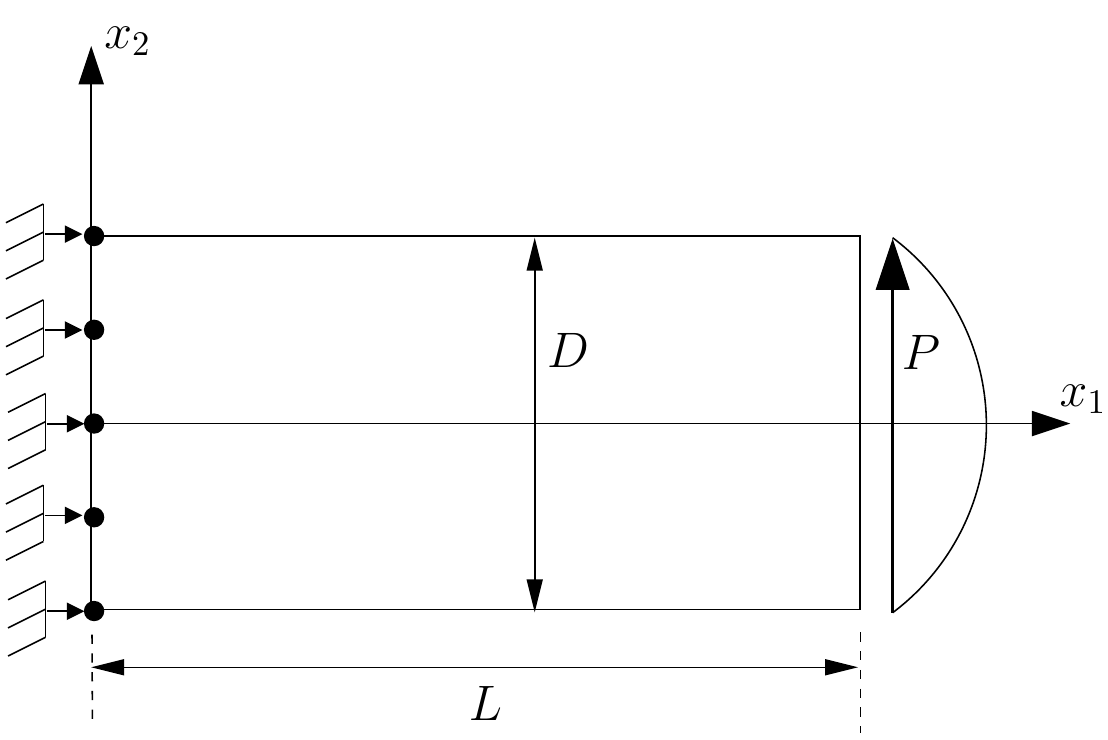}}
	\subfigure[] {\label{fig:beamproblem_b}\includegraphics[width=0.42\linewidth]{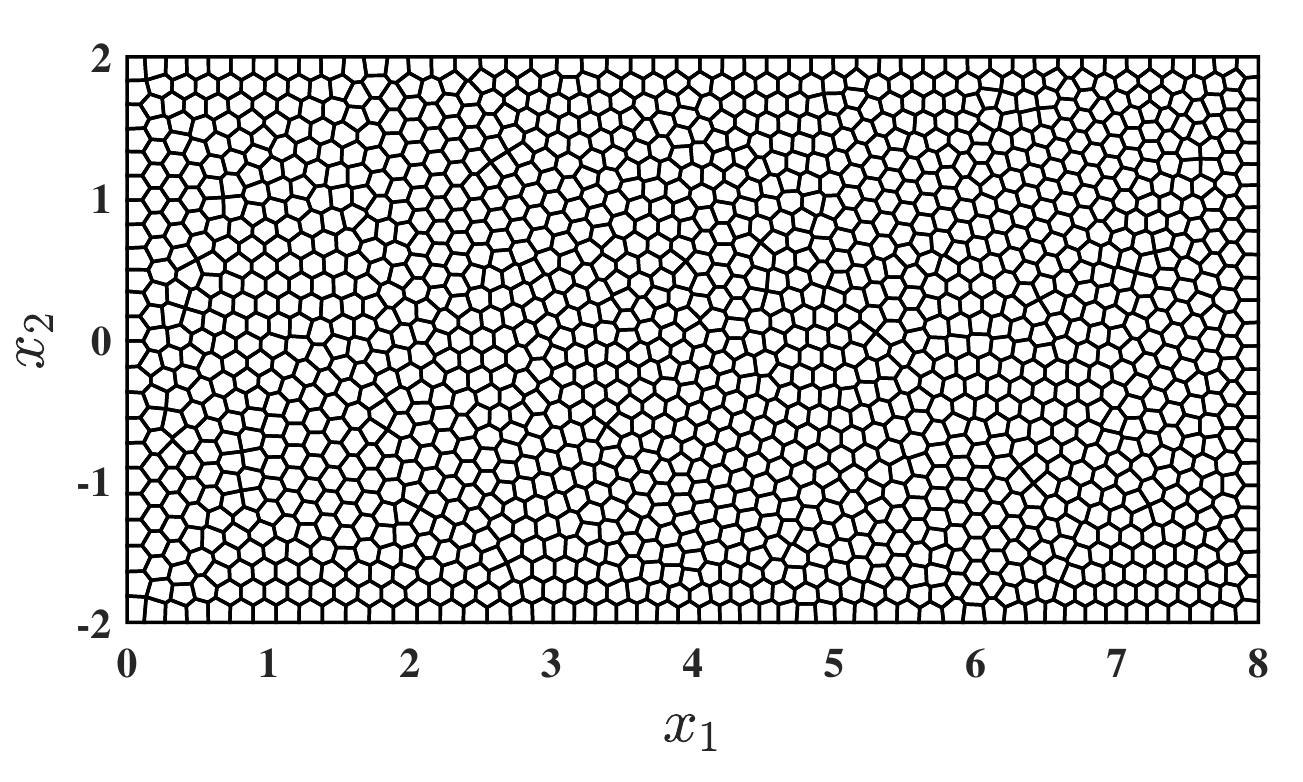}}
	}
	\caption{Cantilever beam (a) geometry and boundary conditions, and (b) a sample mesh.}
	\label{fig:beamproblem}
\end{figure}

\alejandro{With dimensions in inches,} a cantilever beam of unit thickness 
whose domain of analysis is $\Omega=(0,8)\times (-2,2)$
is subjected to a parabolic end load $P$ that is applied on the edge $x=8$. 
The Dirichlet boundary conditions are applied on the edge $x=0$ according to the following exact 
solution given by Timoshenko and Goodier~\cite{Timoshenko-Goodier:1970} for plane strain state:
\begin{subequations}\label{eq:beam_exact_sol}
\begin{align*}
u_1 &=  -\frac{Px_2}{6\overline{E}_Y I}\left((6L-3x_1)x_1 + (2+\overline{\nu})x_2^2 - \frac{3D^2}{2}(1+\overline{\nu})\right),\\
u_2 &= \frac{P}{6\overline{E}_Y I}\left(3\overline{\nu}x_2^{2}(L-x_1)+(3L-x_1)x_1^{2}\right)
\end{align*}
\end{subequations}
with $\overline{E}_Y=E_Y/\left(1-\nu^{2}\right)$ and $\overline{\nu}=\nu/\left(1-\nu\right)$, 
where $E_Y$ and $\nu$ are the Young's modulus and the Poisson's ratio of the linear 
elastic material, respectively; $L$ is the length of the beam, $D$ is the height of the beam, and $I$ is the
second-area moment of the beam section. The total load on the
Neumann boundary is $P=-1000$ \alejandro{lbf}. The geometry and boundary conditions, and a sample 
polygonal mesh used in this study are shown in~\fref{fig:beamproblem}.

The convergence and accuracy of the NVEM is studied for two sets of material parameters.
The first set is $E=10^7$ \alejandro{psi} and $\nu=0.3$ (compressible elasticity), and the second set
is $E=10^7$ \alejandro{psi} and $\nu=0.499999$ (nearly incompressible elasticity). 

\subsubsection{Compressible elasticity}
\begin{figure}[!bth]
	\centering
	\mbox{
	\subfigure[] {\label{fig:beamrates1_a}\includegraphics[width=0.5\linewidth]{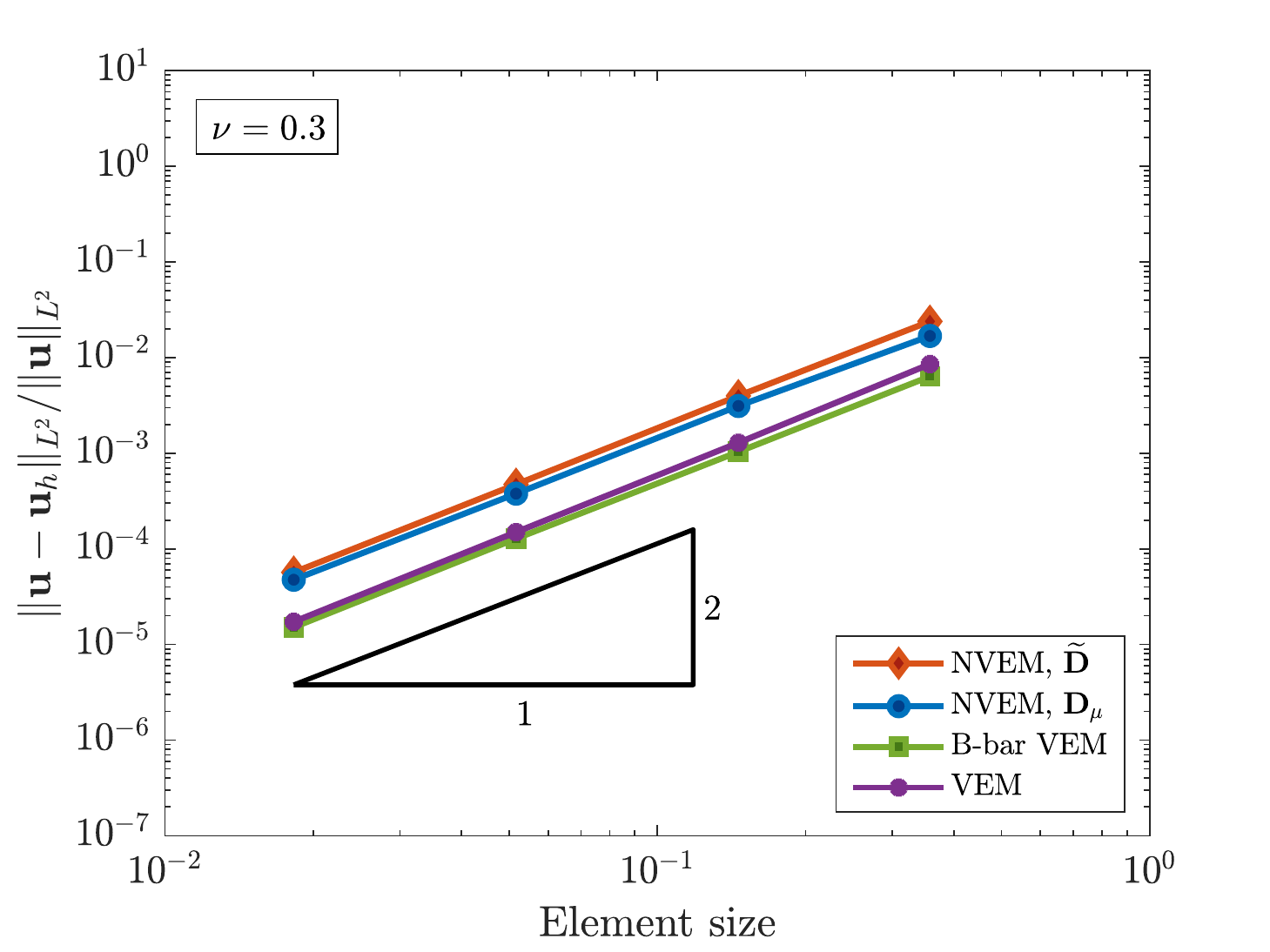}}
	\subfigure[] {\label{fig:beamrates1_b}\includegraphics[width=0.5\linewidth]{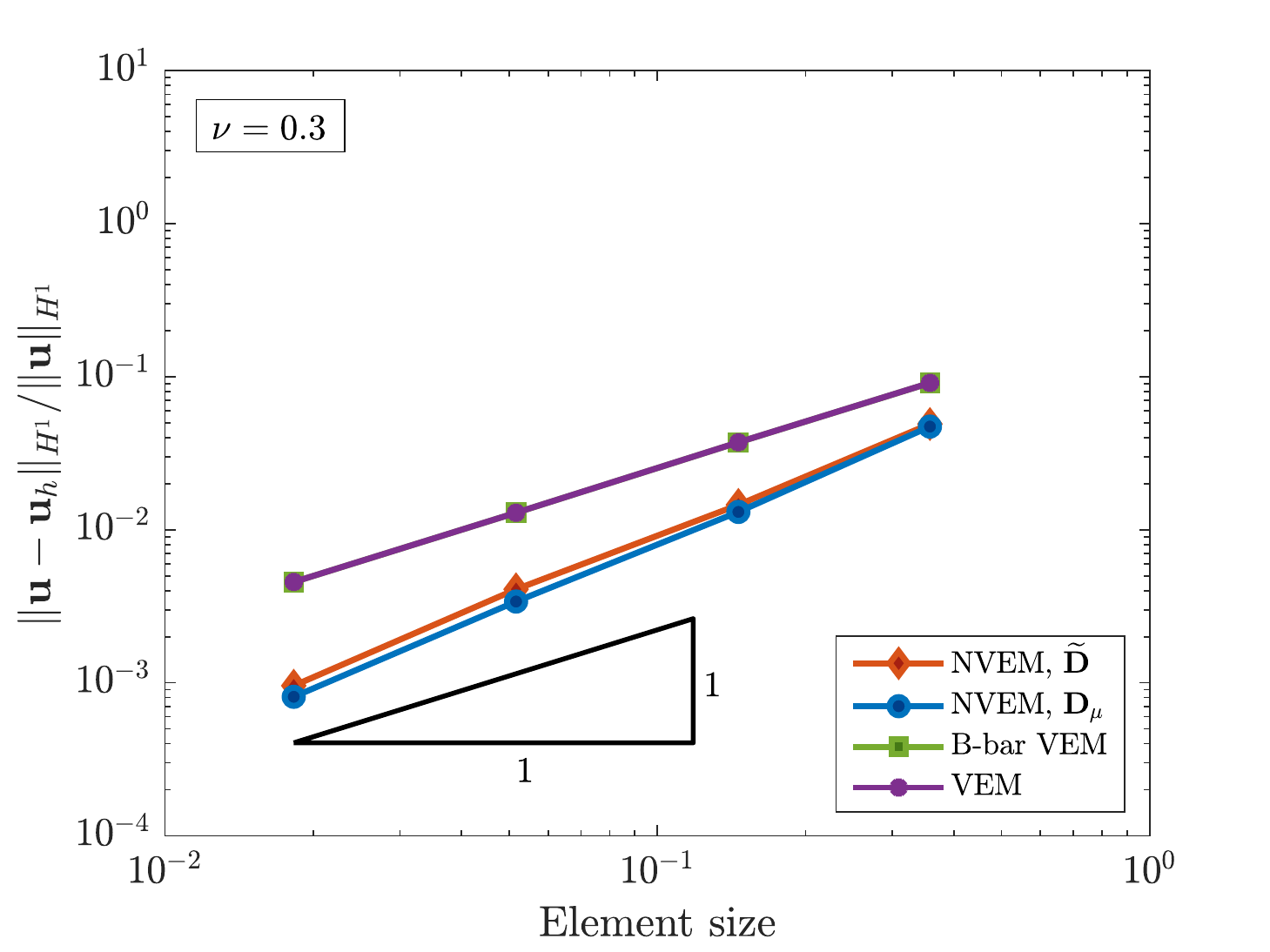}}
	}
	\subfigure[] {\label{fig:beamrates1_c}\includegraphics[width=0.5\linewidth]{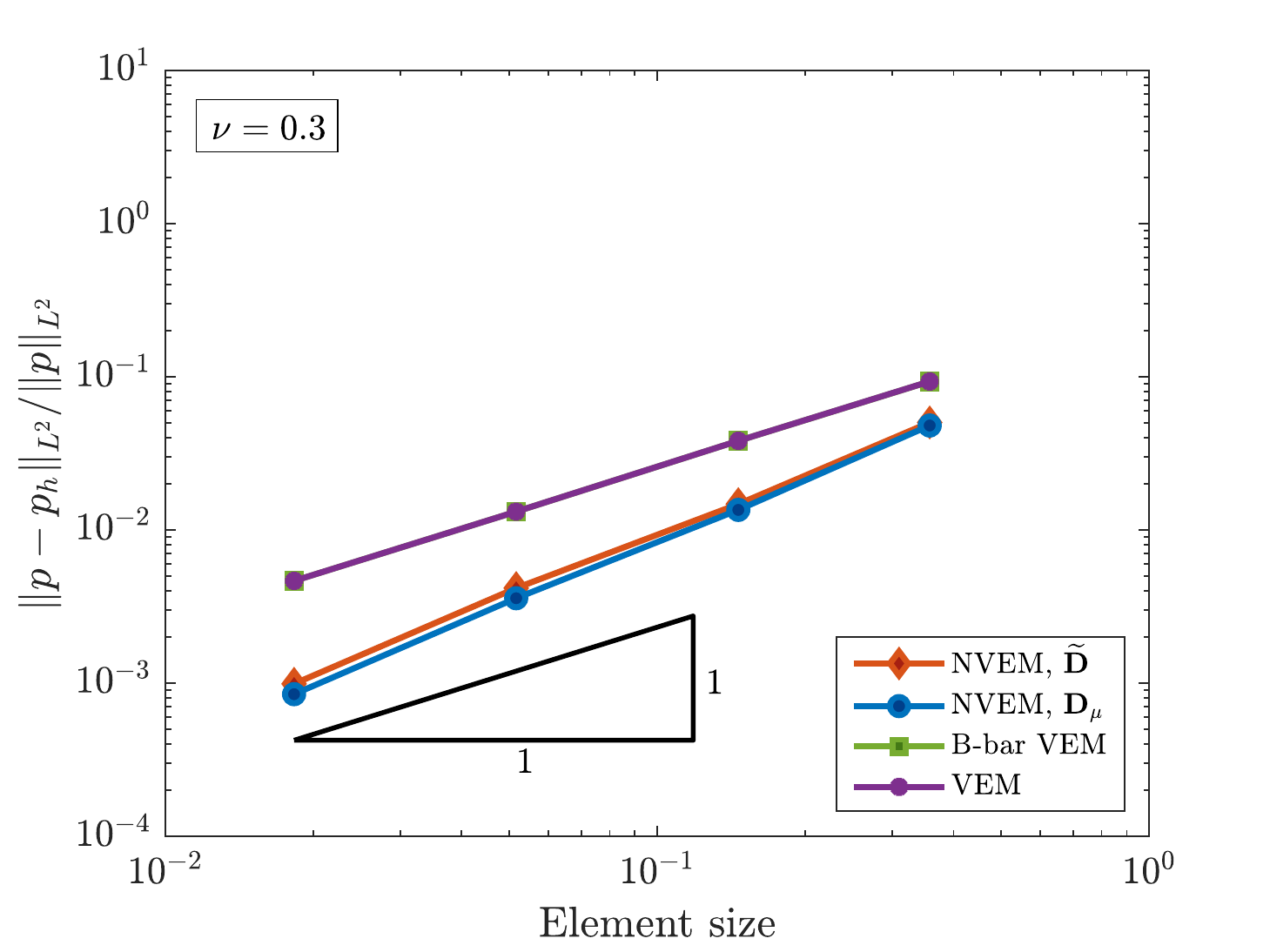}}	
	\caption{Cantilever beam problem with material parameters $E_Y=10^7$ \alejandrog{psi} and $\nu=0.3$. 
	Convergence rates in the (a) $L^2$ norm of the displacement error, (b) $H^1$ seminorm 
	of the displacement error and (c) $L^2$ norm of the pressure error for the VEM, B-bar VEM and NVEM.}
	\label{fig:beamrates1}
\end{figure}

For compressible elasticity, \fref{fig:beamrates1} presents the $L^2$ norm and 
the $H^1$ seminorm of the displacement error, and the $L^2$ norm of the pressure error, 
where it is observed that accurate solutions with optimal convergence rates of 2, 1 and 1, 
respectively, are delivered by the NVEM with $\widetilde{\vm{D}}$ and $\vm{D}_\mu$ stabilizations.
The VEM and the B-bar VEM also exhibit accurate solutions with optimal rates of convergence. In
comparison with the VEM and the B-bar VEM, the NVEM is more accurate in the
$H^1$ seminorm of the displacement error and the $L^2$ norm of the pressure error, 
whereas less accurate in the $L^2$ norm of the displacement error. 

\subsubsection{Nearly incompressible elasticity}
\begin{figure}[!bth]
	\centering
	\mbox{
	\subfigure[] {\label{fig:beamrates2_a}\includegraphics[width=0.5\linewidth]{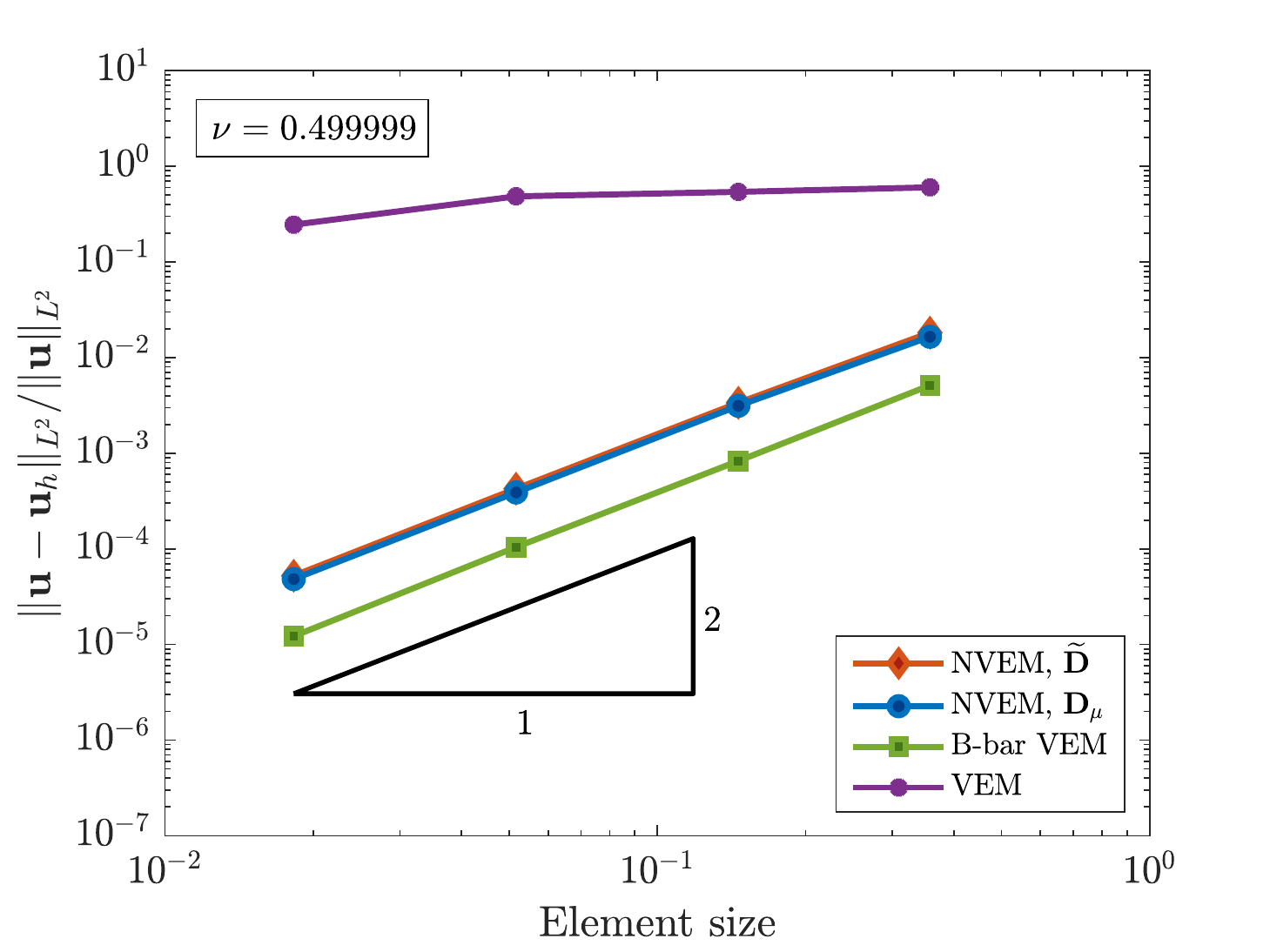}}
	\subfigure[] {\label{fig:beamrates2_b}\includegraphics[width=0.5\linewidth]{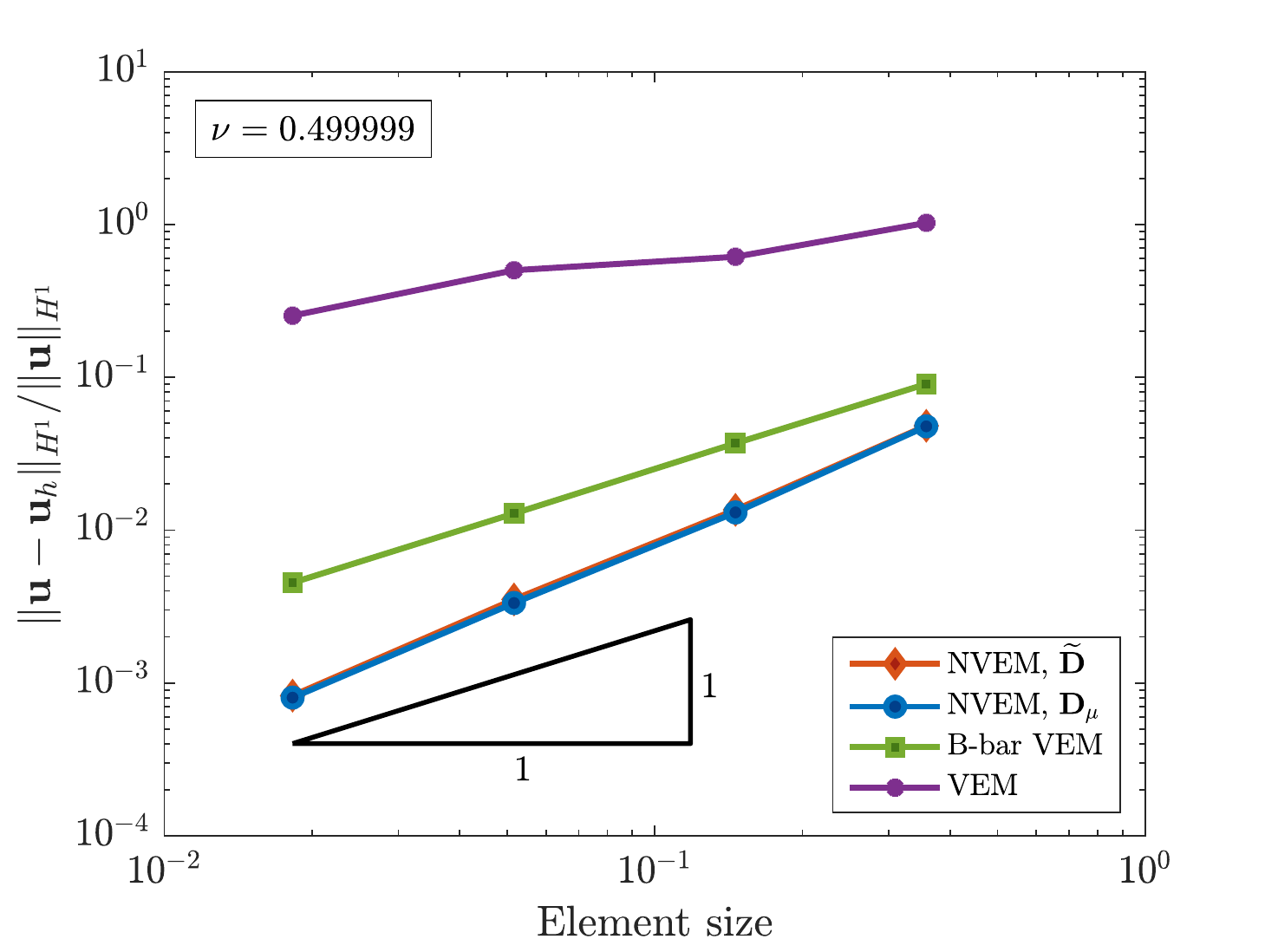}}
	}
	\subfigure[] {\label{fig:beamrates2_c}\includegraphics[width=0.5\linewidth]{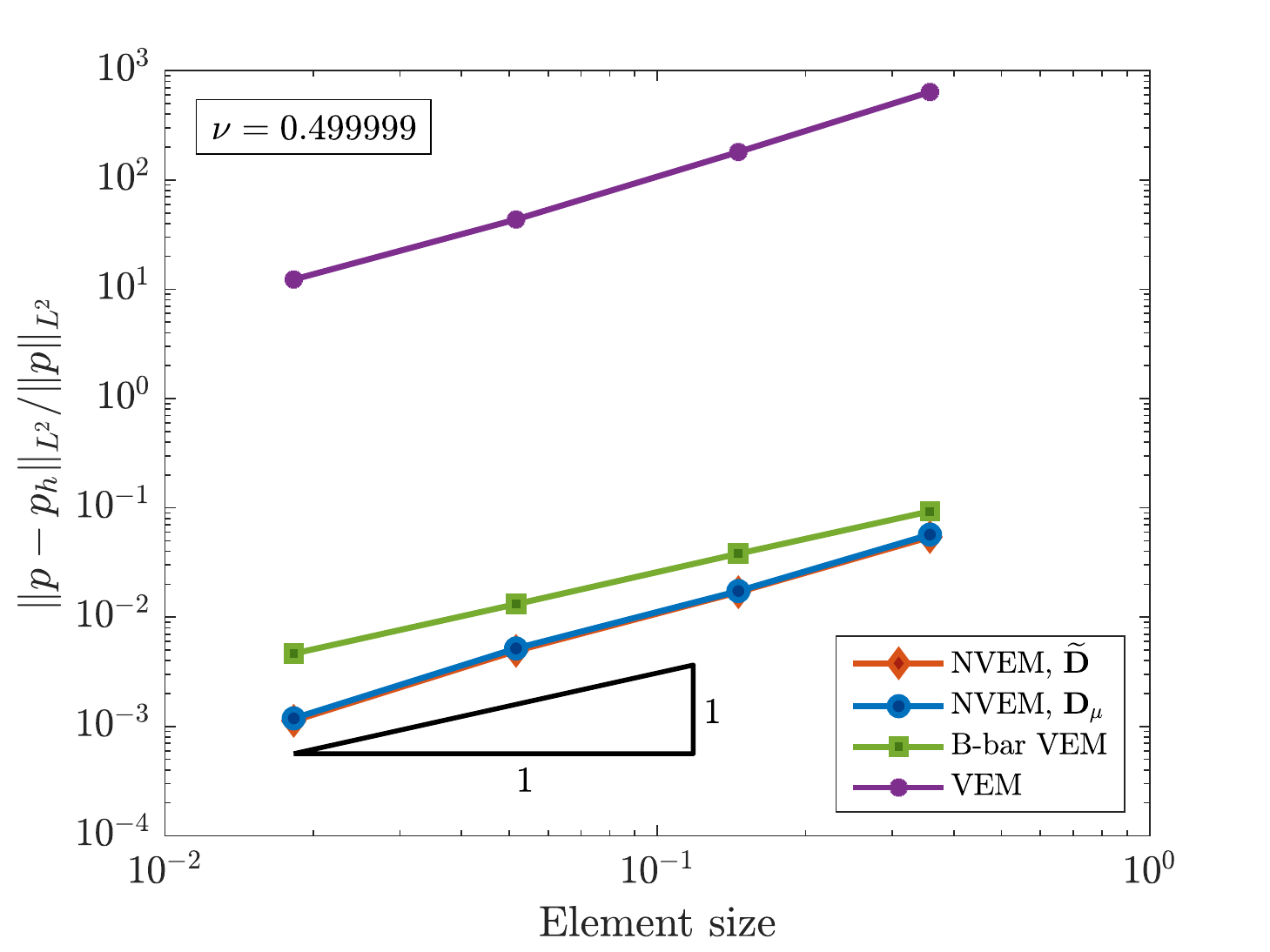}}	
	\caption{Cantilever beam problem with material parameters $E_Y=10^7$ \alejandrog{psi} and $\nu=0.499999$. 
	Convergence rates in the (a) $L^2$ norm of the displacement error, (b) $H^1$ seminorm 
	of the displacement error and (c) $L^2$ norm of the pressure error for the VEM, B-bar VEM and NVEM.}
	\label{fig:beamrates2}
\end{figure}

For nearly incompressible elasticity, the NVEM and the B-bar VEM are accurate and 
optimally convergent, as shown in~\fref{fig:beamrates2}. Like in the compressible case,
the NVEM is more accurate than the B-bar VEM in the $H^1$ seminorm of the displacement error 
and the $L^2$ norm of the pressure error, whereas less accurate in the $L^2$ norm of the 
displacement error. As expected, \fref{fig:beamrates2} also exhibits very inaccurate solutions 
for the VEM due to volumetric locking.

\subsection{Cook's membrane}
\label{sec:cook}

\begin{figure}[!bth]
	\centering
	\includegraphics[width=0.4\linewidth]{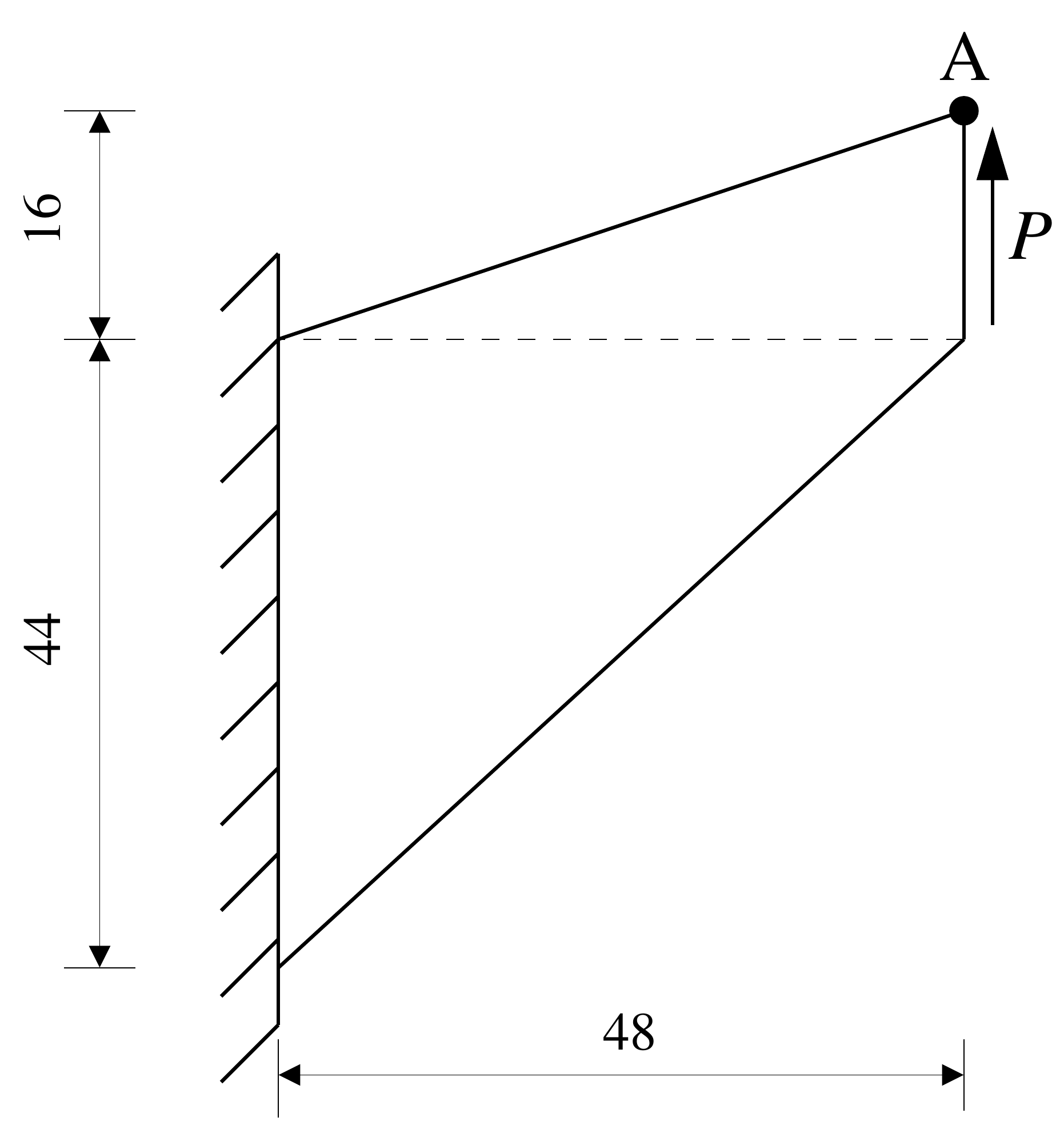}
	\caption{Geometry and boundary conditions for the Cook's membrane problem \alejandro{(dimensions in mm)}.}
	\label{fig:cookproblem}
\end{figure}

\begin{figure}[!bth]
	\centering
	\mbox{
	\subfigure[] {\label{fig:cookmesh_a}\includegraphics[width=0.36\linewidth]{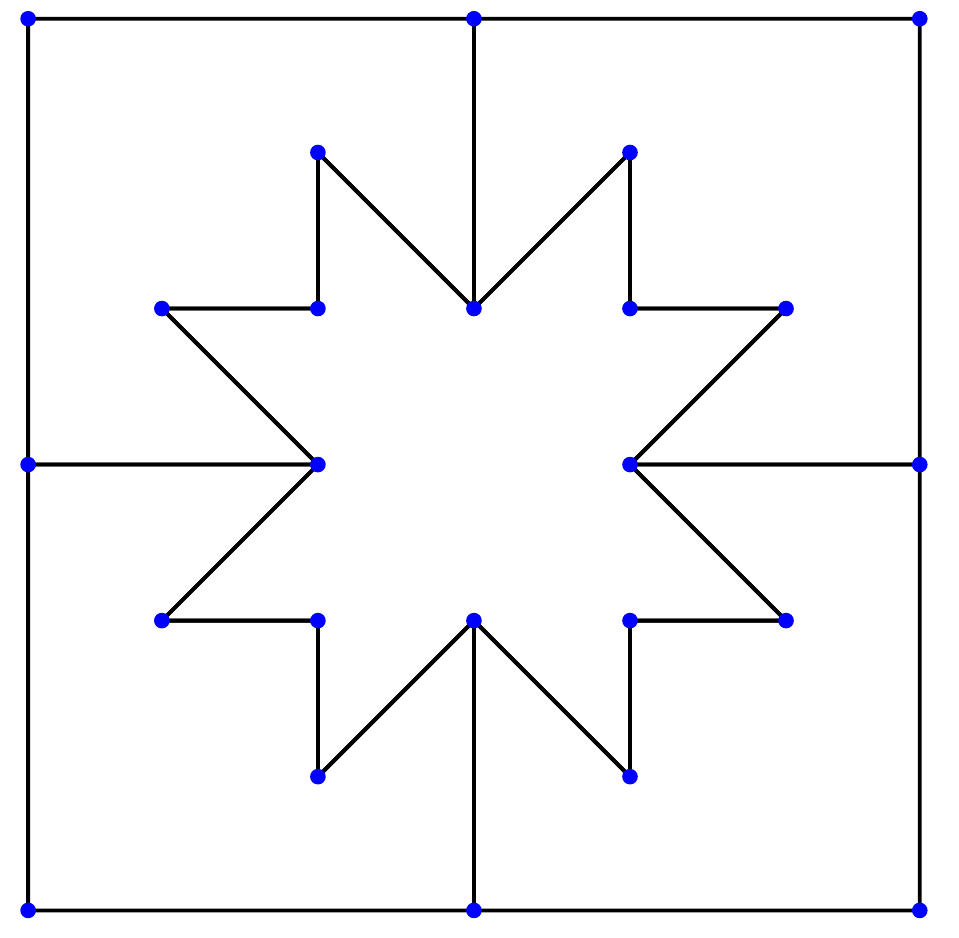}}
	\subfigure[] {\label{fig:cookmesh_b}\includegraphics[width=0.38\linewidth]{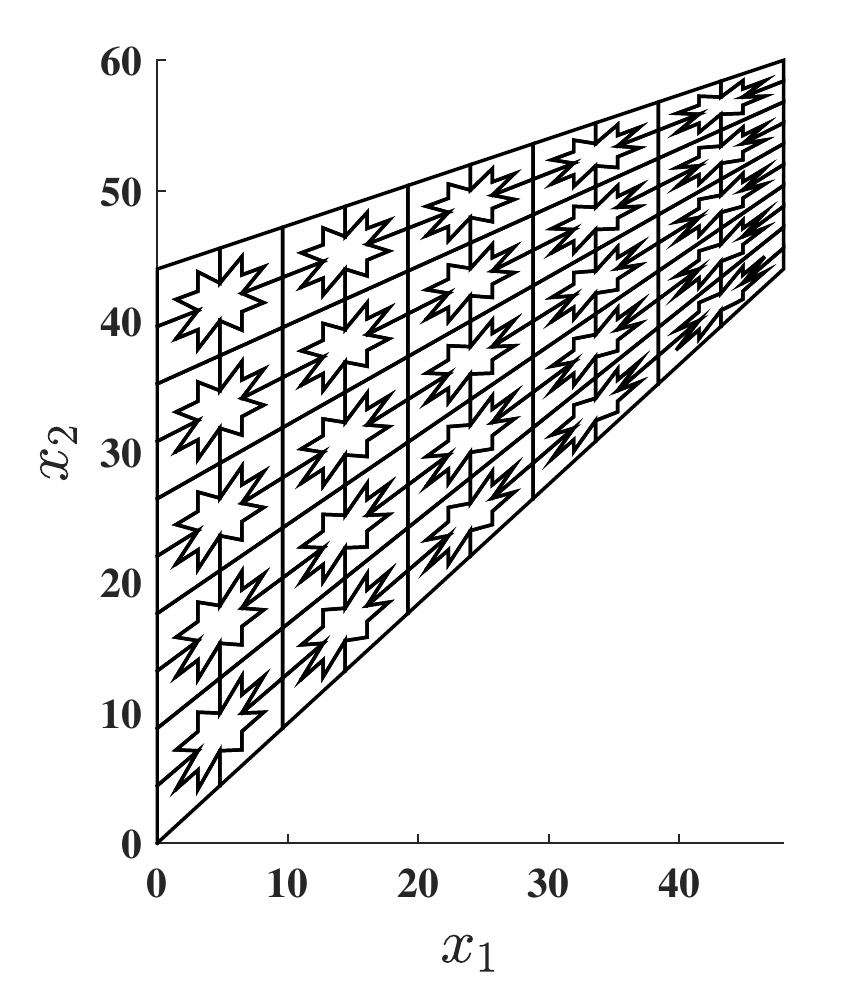}}
	}
	\caption{Cook's membrane problem. (a) Master multielement. The center polygon defines 
	an eight-point star known as `Gu\~nelve'; (b) a sample `Gu\~nelve' mesh.}
	\label{fig:cookmesh}
\end{figure}

This standard benchmark problem is suitable to test the performance of numerical 
formulations for nearly incompressible solid materials under combined bending and shear.
The geometry and boundary conditions 
for this problem are schematically shown in \fref{fig:cookproblem}. 
There, the left edge is clamped and the 
right edge is subjected to a shear load $P = 6.25$ \alejandro{N/mm}
(total shear load of $100$ \alejandro{N}). The following material parameters are
used: $E=250$ \alejandro{MPa} and $\nu=0.4999$. In this example, a non standard mesh is tested.
To this end, we define a master multielement region formed by five polygons, 
as shown in \fref{fig:cookmesh_a}, where the center polygon is an eight-point 
star known as `Gu\~nelve'. This master multielement is mapped into the
membrane geometry several times to form the mesh. A sample `Gu\~nelve' mesh is shown
in \fref{fig:cookmesh_b}.

\begin{figure}[!bth]
	\centering
	\includegraphics[width=0.6\linewidth]{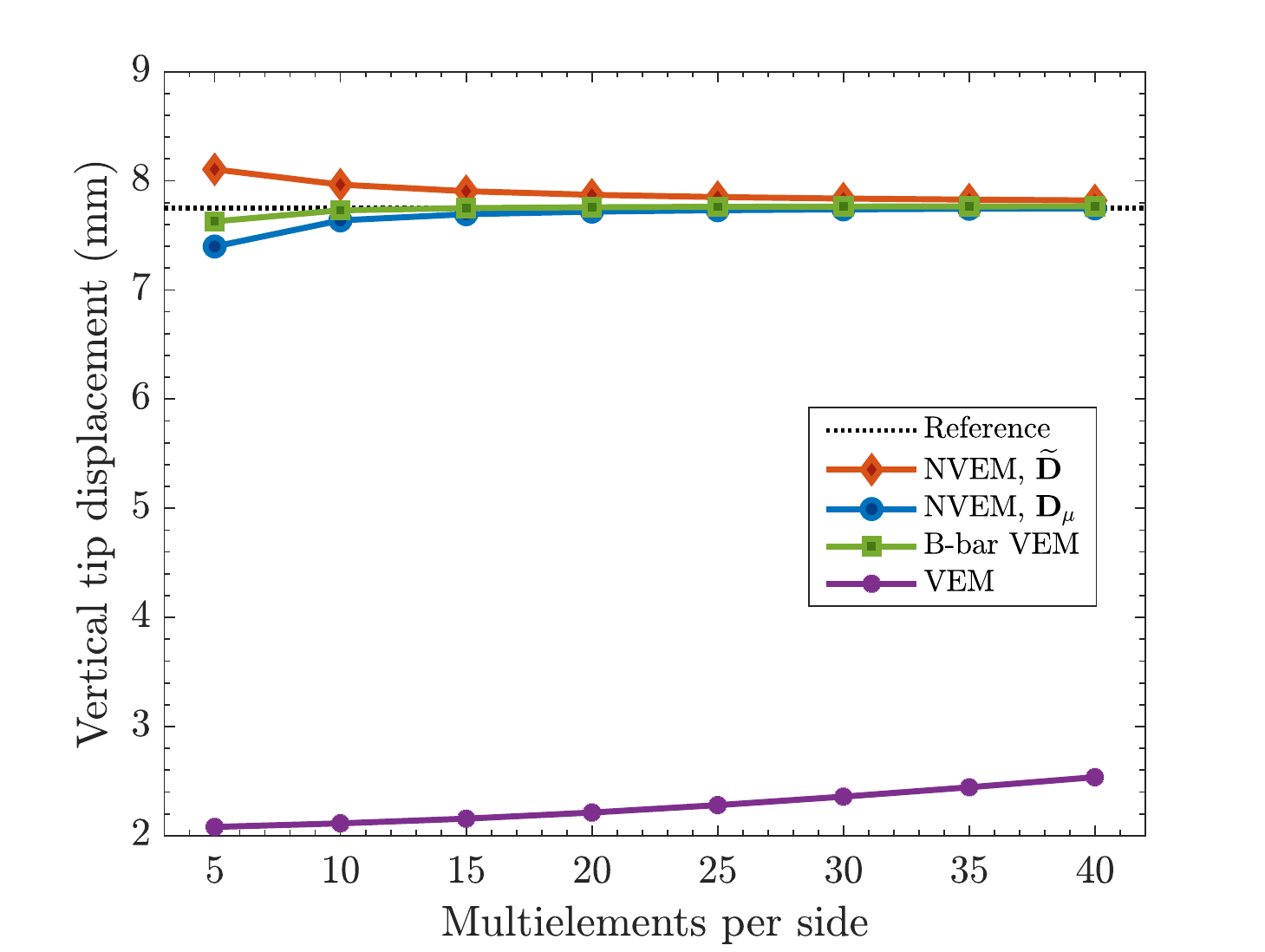}	
	\caption{Convergence of the vertical displacement at the tip of the Cook's 
	membrane (point A) upon mesh refinement for the VEM, B-bar VEM and NVEM.}

	\label{fig:cookconvergence}
\end{figure}

The convergence of the vertical tip displacement at point A upon mesh refinement is 
studied. For comparison purposes, a reference solution
for the vertical tip displacement is found using a highly refined mesh of eight-node quadrilateral
finite elements. The result of this study is summarized in \fref{fig:cookconvergence},
where it is shown that the B-bar VEM and the NVEM solutions approach the reference value 
as the mesh is refined, and, as expected, the standard VEM exhibits very poor convergence
due to volumetric locking.

\begin{figure}[!bth]
	\centering
	\mbox{
	\subfigure[] {\label{fig:cookcontouruy_a}\includegraphics[width=0.4\linewidth]{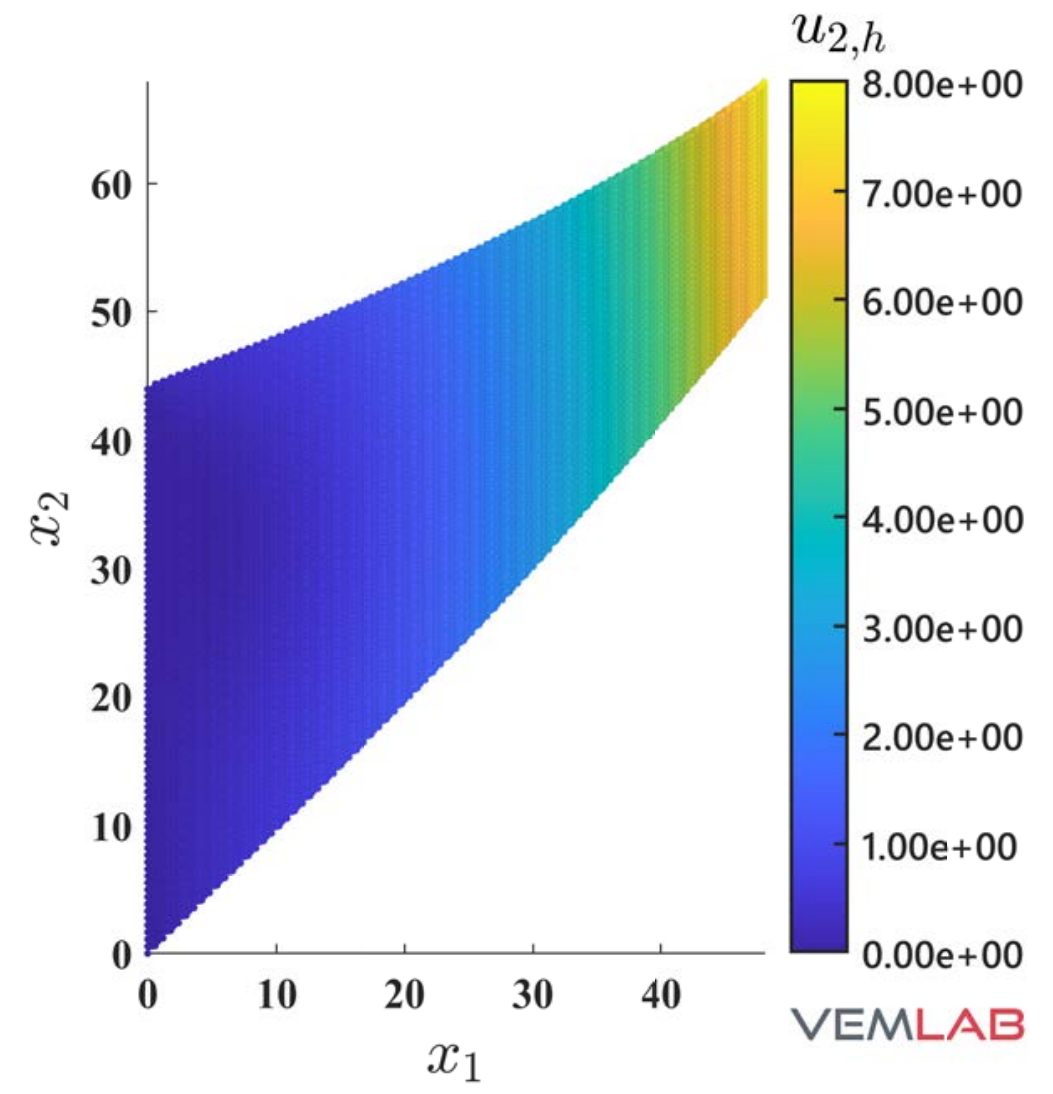}}
	\subfigure[] {\label{fig:cookcontouruy_b}\includegraphics[width=0.4\linewidth]{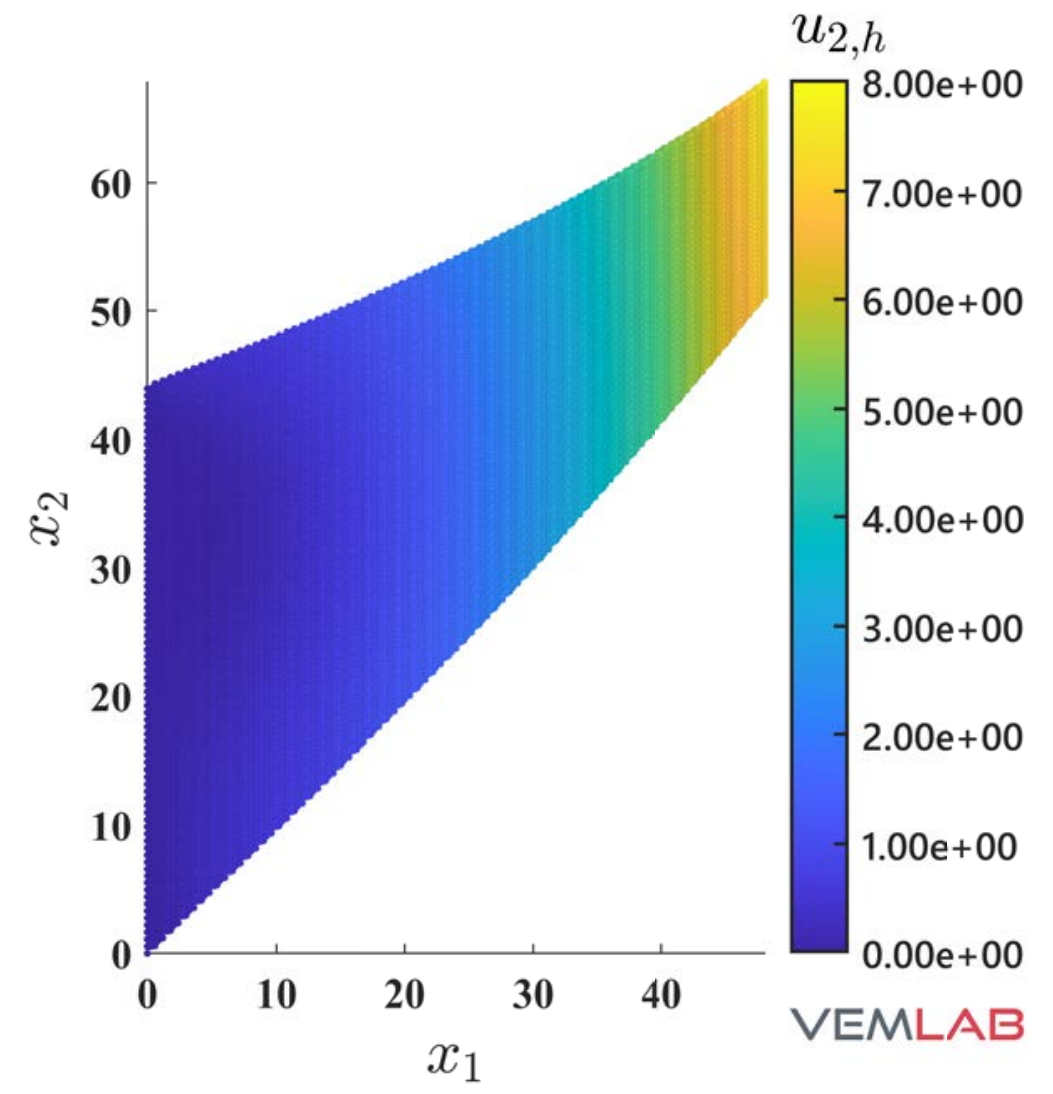}}
	}
	\subfigure[] {\label{fig:cookcontouruy_c}\includegraphics[width=0.4\linewidth]{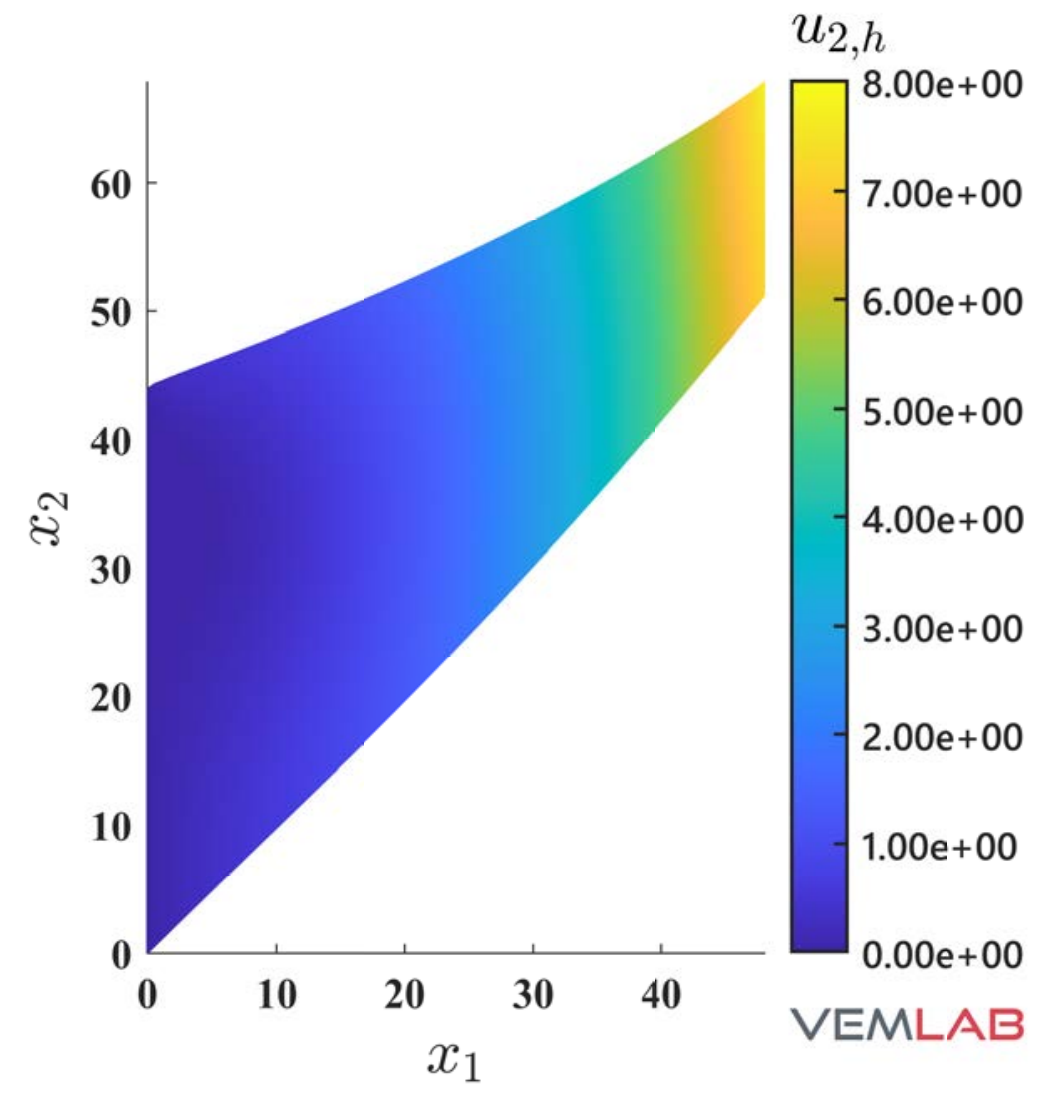}}
	\caption{Cook's membrane problem. Plots of the vertical displacement field
	solution ($u_{2,h}$) \alejandro{in mm} for the (a) NVEM ($\tilde{\vm{D}}$ stabilization), 
	(b) NVEM ($\vm{D}_{\mu}$ stabilization) and (c) B-bar VEM approaches. Plots are deformed according to $u_{2,h}$.}
	\label{fig:cookcontouruy}
\end{figure}

\begin{figure}[!bth]
	\centering
	\mbox{
	\subfigure[] {\label{fig:cookcontourp_a}\includegraphics[width=0.4\linewidth]{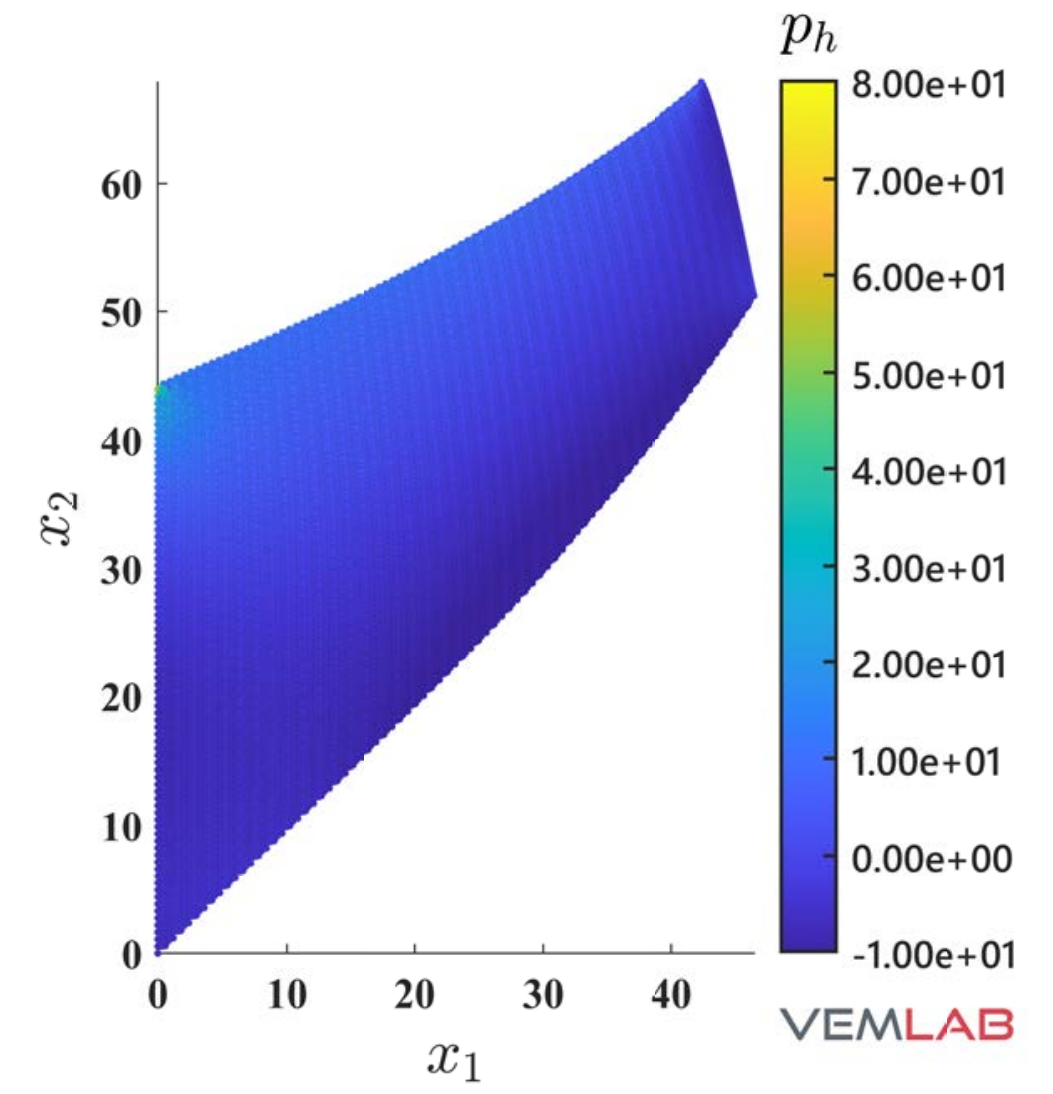}}
	\subfigure[] {\label{fig:cookcontourp_b}\includegraphics[width=0.4\linewidth]{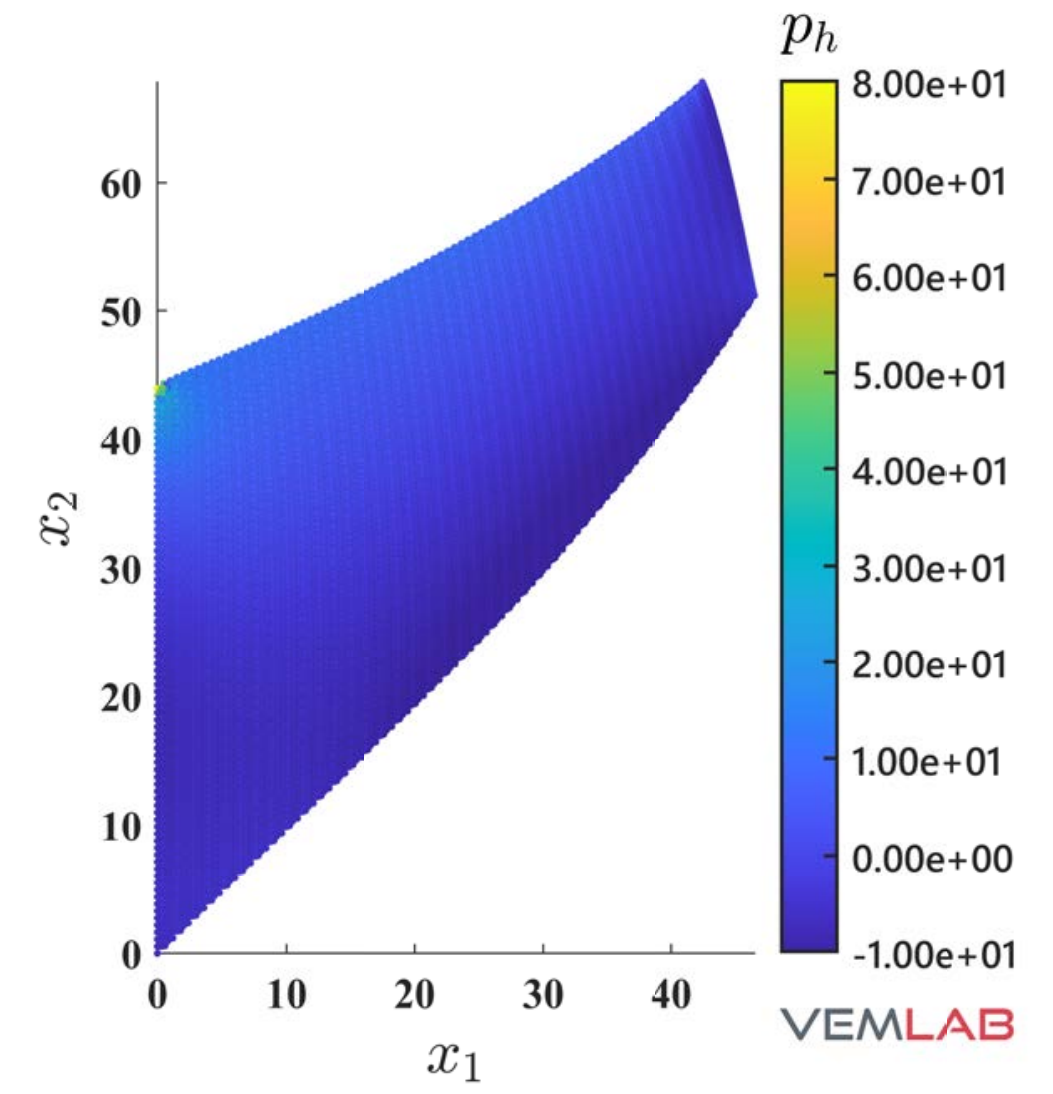}}
	}
	\subfigure[] {\label{fig:cookcontourp_c}\includegraphics[width=0.4\linewidth]{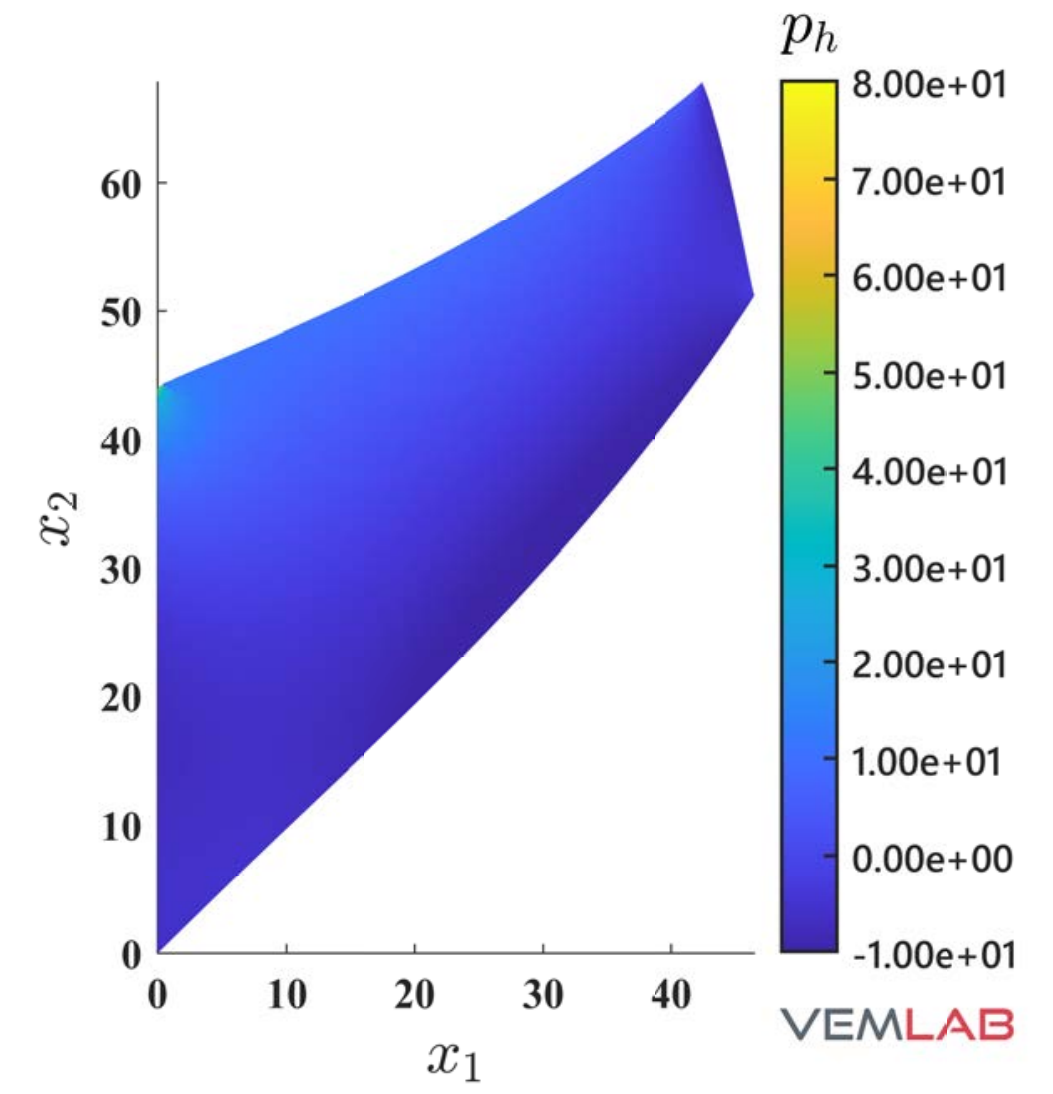}}
	\caption{Cook's membrane problem. Plots of the pressure field
	solution ($p_h$) \alejandro{in MPa} for the (a) NVEM ($\tilde{\vm{D}}$ stabilization), 
	(b) NVEM ($\vm{D}_{\mu}$ stabilization) and (c) B-bar VEM approaches. Plots are deformed according to $\|\vm{u}_h\|$.}
	\label{fig:cookcontourp}
\end{figure}

\begin{figure}[!bth]
	\centering
	\mbox{
	\subfigure[] {\label{fig:cookcontourvm_a}\includegraphics[width=0.4\linewidth]{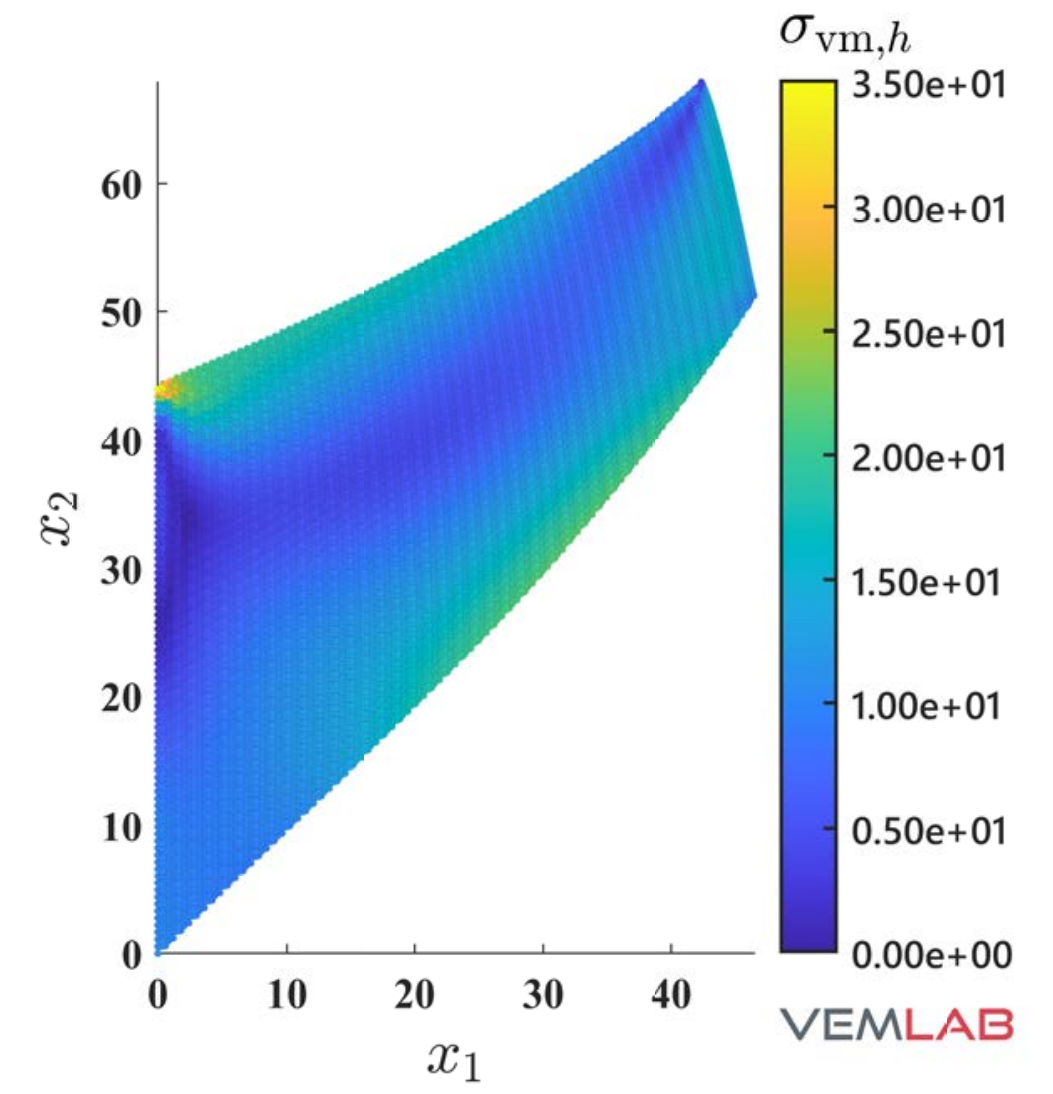}}
	\subfigure[] {\label{fig:cookcontourvm_b}\includegraphics[width=0.4\linewidth]{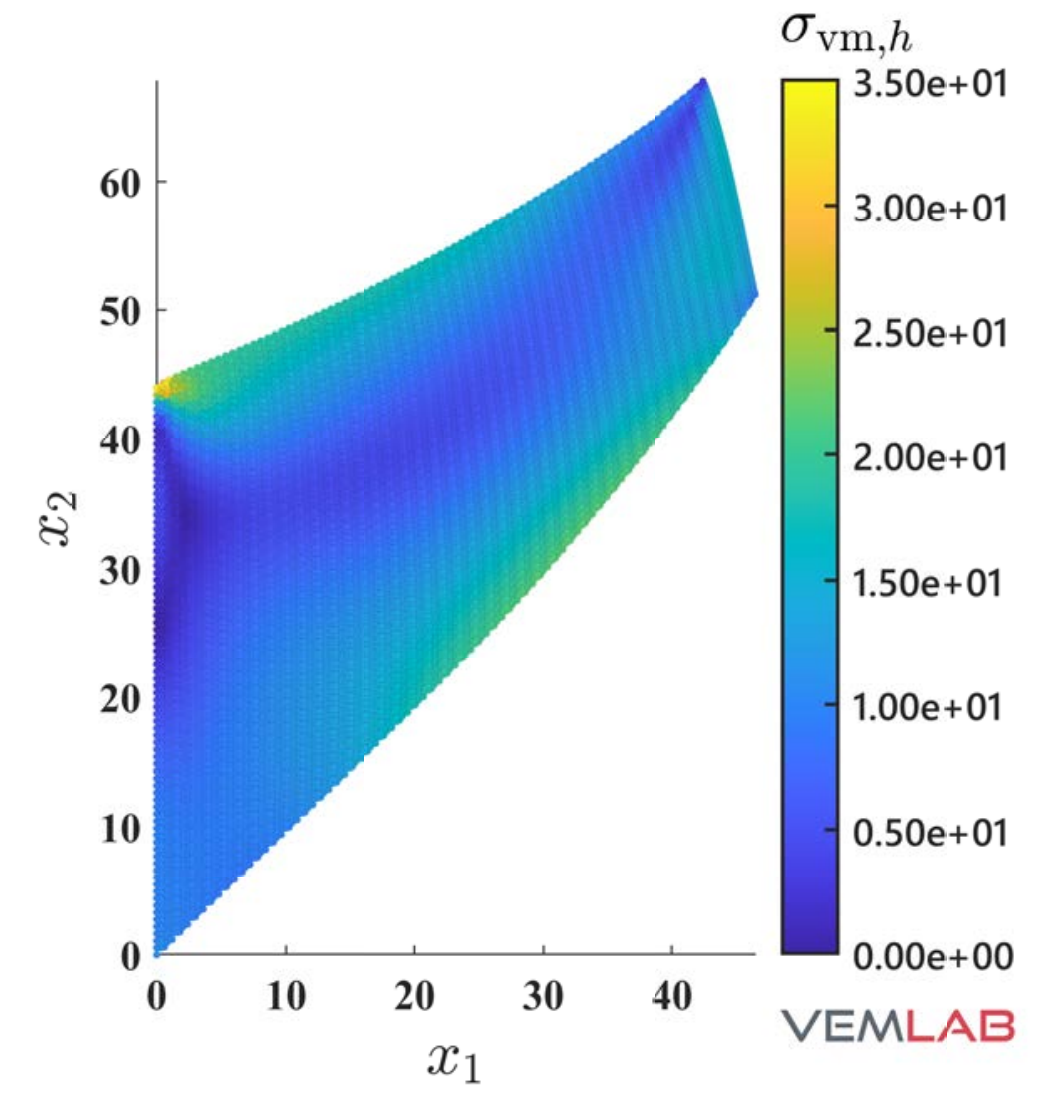}}
	}
	\subfigure[] {\label{fig:cookcontourvm_c}\includegraphics[width=0.4\linewidth]{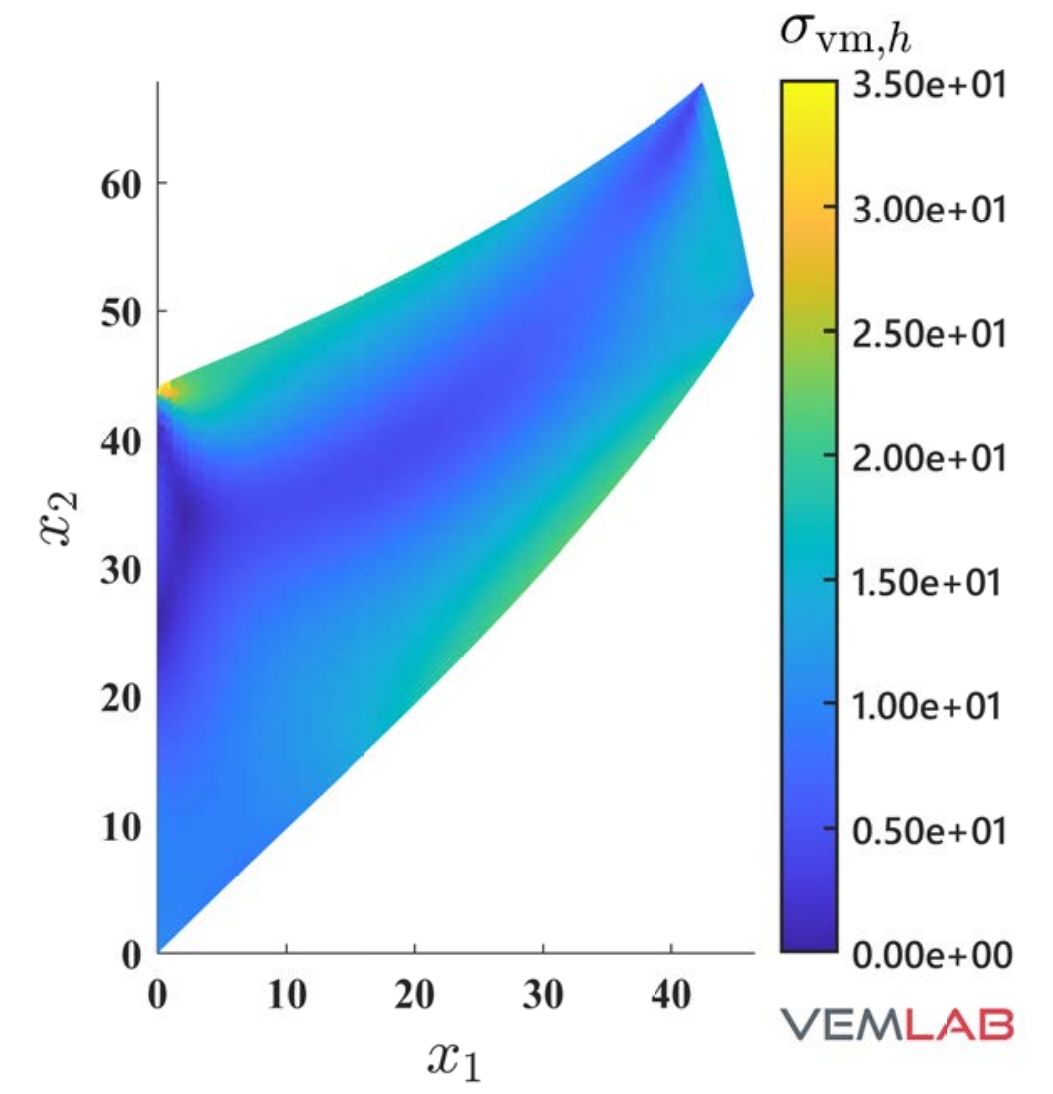}}	
	\caption{Cook's membrane problem. Plots of the von Mises stress field
	solution ($\sigma_{\textrm{vm},h}$) \alejandro{in MPa} for the (a) NVEM ($\tilde{\vm{D}}$ stabilization), 
	(b) NVEM ($\vm{D}_{\mu}$ stabilization) and (c) B-bar VEM approaches. Plots are deformed according to $\|\vm{u}_h\|$.}
	\label{fig:cookcontourvm}
\end{figure}

Plots of the vertical displacement, pressure and von Mises stress fields are presented in
Figs. \ref{fig:cookcontouruy}, \ref{fig:cookcontourp} and \ref{fig:cookcontourvm}, respectively.
Scatter plots are used for the NVEM as in this approach the field variables are known at the nodes. 
As expected, a very good agreement between the NVEM and the B-bar VEM solutions is observed in these
plots.

\subsection{Infinite plate with a circular hole}
\label{sec:platehole}

\begin{figure}[!bth]
	\centering
	\mbox{
	\subfigure[] {\label{fig:plateholeproblem_a}\includegraphics[width=0.38\linewidth]{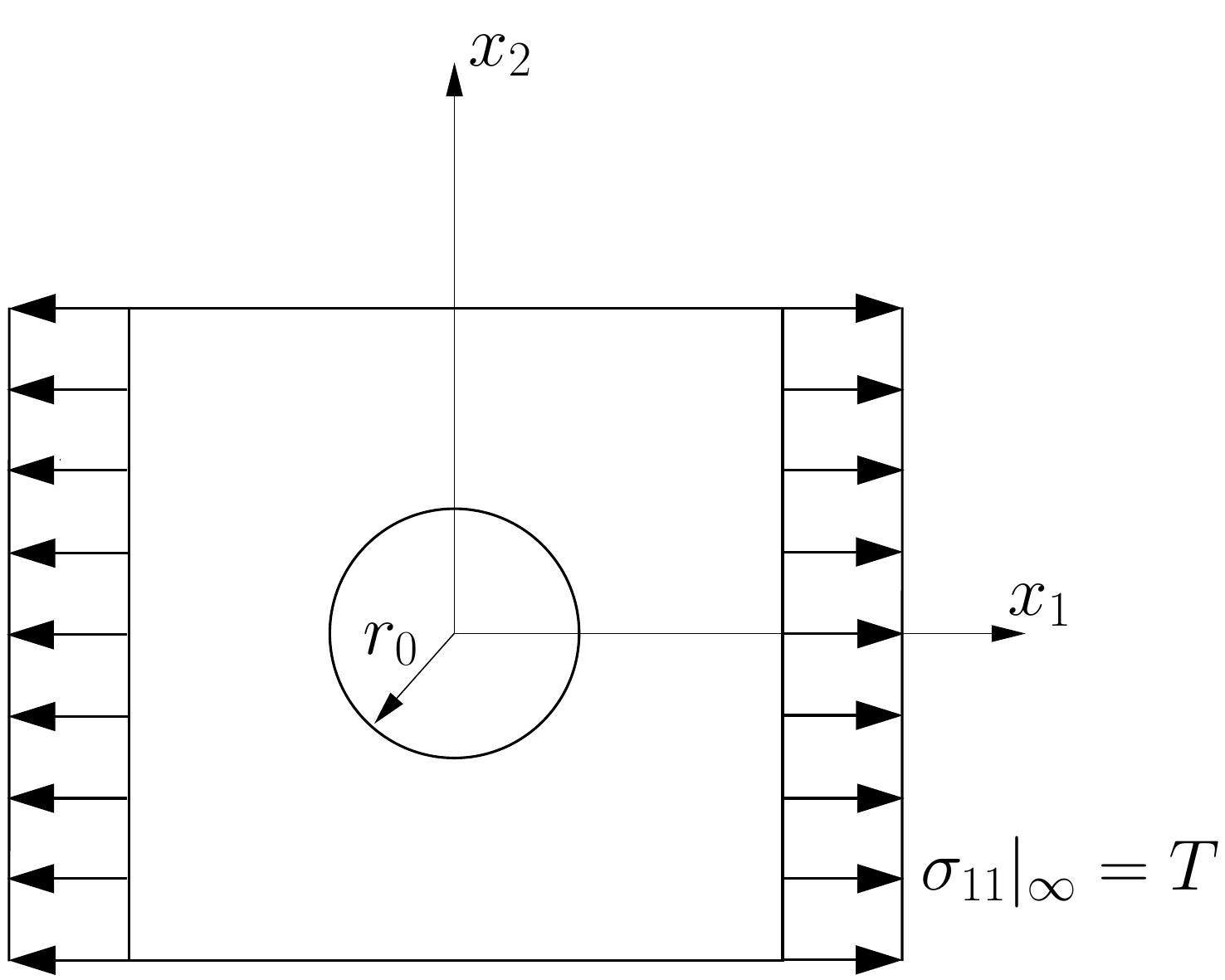}}
	\subfigure[] {\label{fig:plateholeproblem_b}\includegraphics[width=0.38\linewidth]{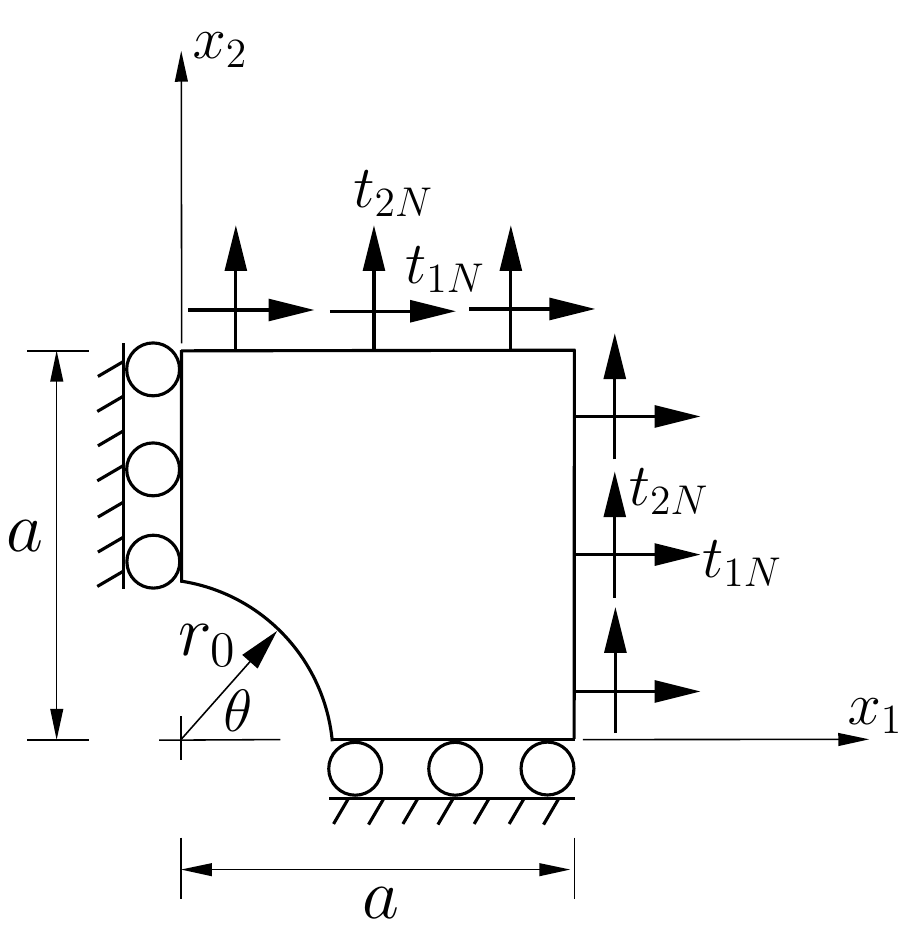}}	
	}
	\caption{Geometry and boundary conditions for the infinite plate with a circular hole. 
	(a) Infinite domain and (b) a quarter of the domain.}
	\label{fig:plateholeproblem}
\end{figure}

This example is devoted to study an infinite plate with a circular hole that is loaded 
at infinity with the following tractions:
$\sigma_{11}=T$ and $\sigma_{22}=\sigma_{12}=0$ (see \fref{fig:plateholeproblem_a}). 
Due to the symmetry of the geometry and boundary conditions, a quarter of
the domain is used as the domain of analysis (see \fref{fig:plateholeproblem_b}). Plane stress and plane 
strain conditions are considered. For plane stress condition, the material parameters are set to
$E_\mathrm{Y}=10^3$ \alejandro{psi} and $\nu = 0.3$ (compressible elasticity), and for plane strain condition, they are set to
$E_\mathrm{Y}=10^3$ \alejandro{psi} and $\nu = 0.499999$ (nearly incompressible elasticity).
The exact solution is given by~\cite{Timoshenko-Goodier:1970}
\begin{equation*}
\vm{u}=\left[
\begin{array}{c}
\frac{T}{4G} \bigg( \frac{\kappa+1}{2}r\cos{\theta}+
 \frac{r_0^2}{r}\Big( (\kappa+1)\cos{\theta}+\cos{3\theta}\Big)-\frac{r_0^4}{r^3}\cos{3\theta} \bigg)\\
 \frac{T}{4G} \bigg( \frac{\kappa-3}{2}r\sin{\theta}+
 \frac{r_0^2}{r}\Big( (\kappa-1)\sin{\theta}+\sin{3\theta}\Big)-\frac{r_0^4}{r^3}\sin{3\theta} \bigg)
\end{array}
\right],
\end{equation*}
where $G=E_\mathrm{Y}/(2(1+\nu))$ and $\kappa = (3-\nu)/(1+\nu)$. The exact stress field is
\begin{equation*}
\left[\begin{array}{c}
\sigma_{11}\\
\sigma_{22}\\
\sigma_{12}
\end{array}
\right]
=\left[\begin{array}{c}
T\bigg(1-\frac{r_0^2}{r^2}\left( \frac{3}{2}\cos{2\theta}
+\cos{4\theta}\right)+\frac{3r_0^4}{2r^4}\cos{4\theta} \bigg)\\
-T\bigg(\frac{r_0^2}{r^2}\left( \frac{1}{2}\cos{2\theta}-\cos{4\theta}\right)
+\frac{3r_0^4}{2r^4}\cos{4\theta} \bigg)\\
-T\bigg(\frac{r_0^2}{r^2}\left( \frac{1}{2}\sin{2\theta}
+\sin{4\theta}\right)-\frac{3r_0^4}{2r^4}\sin{4\theta}\bigg)
\end{array}
\right],
\end{equation*}
where $r$ is the radial distance from the center $(x_1=0,x_2=0)$ to a
point $(x_1,x_2)$ in the domain of analysis. In the computations, the
following data are used: $T=100$ \alejandro{psi}, $r_0=1$ \alejandro{in}
and $a=5$ \alejandro{in}. The Dirichlet 
boundary conditions \alejandrog{in inches are imposed} on the domain 
of analysis (\fref{fig:plateholeproblem_b}) as follows: $u_{1D} = 0$
on the left side and $u_{2D} = 0$ on the bottom side. 
The Neumann boundary conditions are prescribed using the exact 
stresses, as follows:
$\vm{t}_N=\bigl[\; t_{1N} \; \; t_{2N}\; \bigr]^\transpose=\bigl[\;\sigma_{12} \;\; \sigma_{22}\;\bigr]^\transpose$ on the top side and
$\vm{t}_N=\bigl[\; t_{1N} \; \; t_{2N}\; \bigr]^\transpose=\bigl[\;\sigma_{11} \;\; \sigma_{12}\;\bigr]^\transpose$ on the right side.

\begin{figure}[!bth]
	\centering
	\mbox{
	\subfigure[] {\label{fig:plateholemeshes_a}\includegraphics[width=0.38\linewidth]{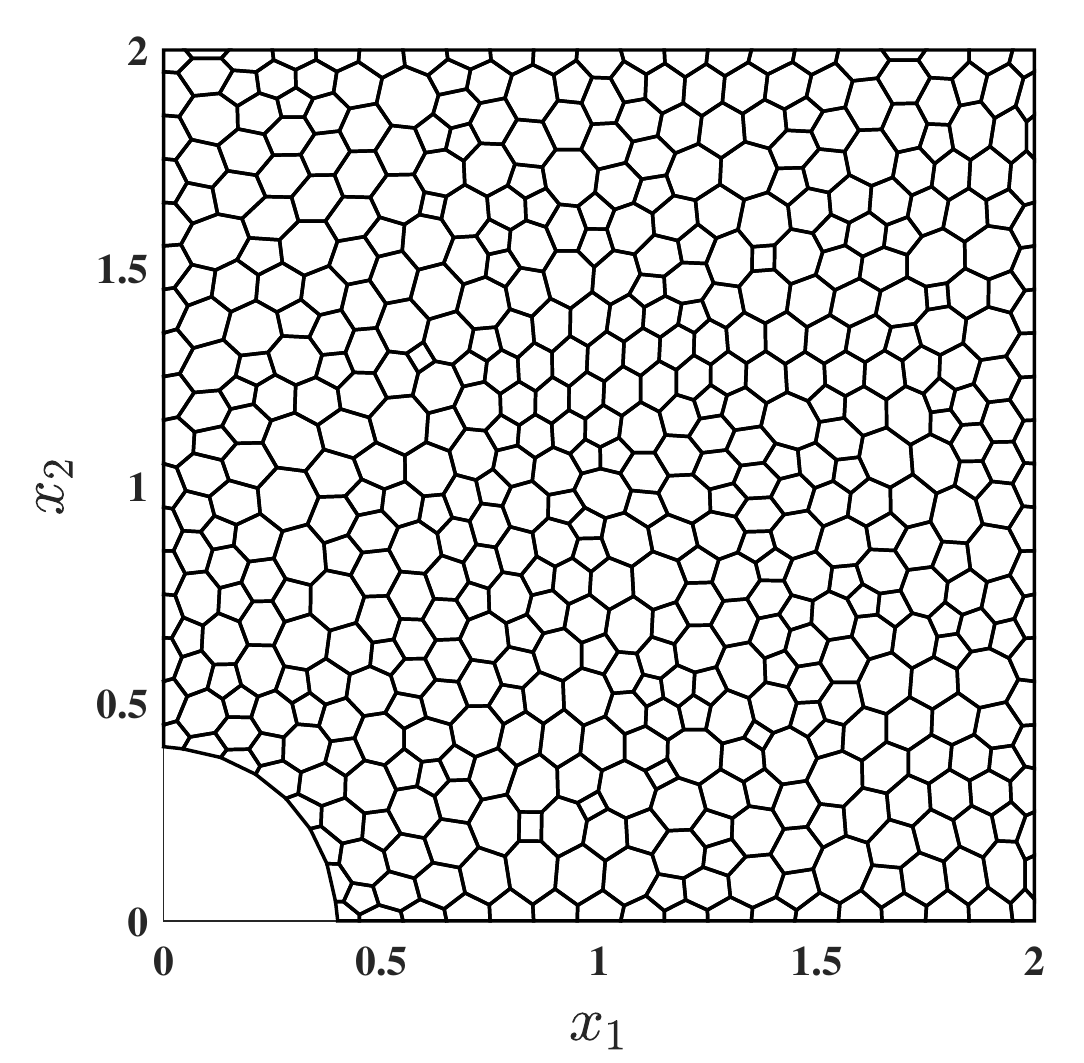}}
	\subfigure[] {\label{fig:plateholemeshes_b}\includegraphics[width=0.38\linewidth]{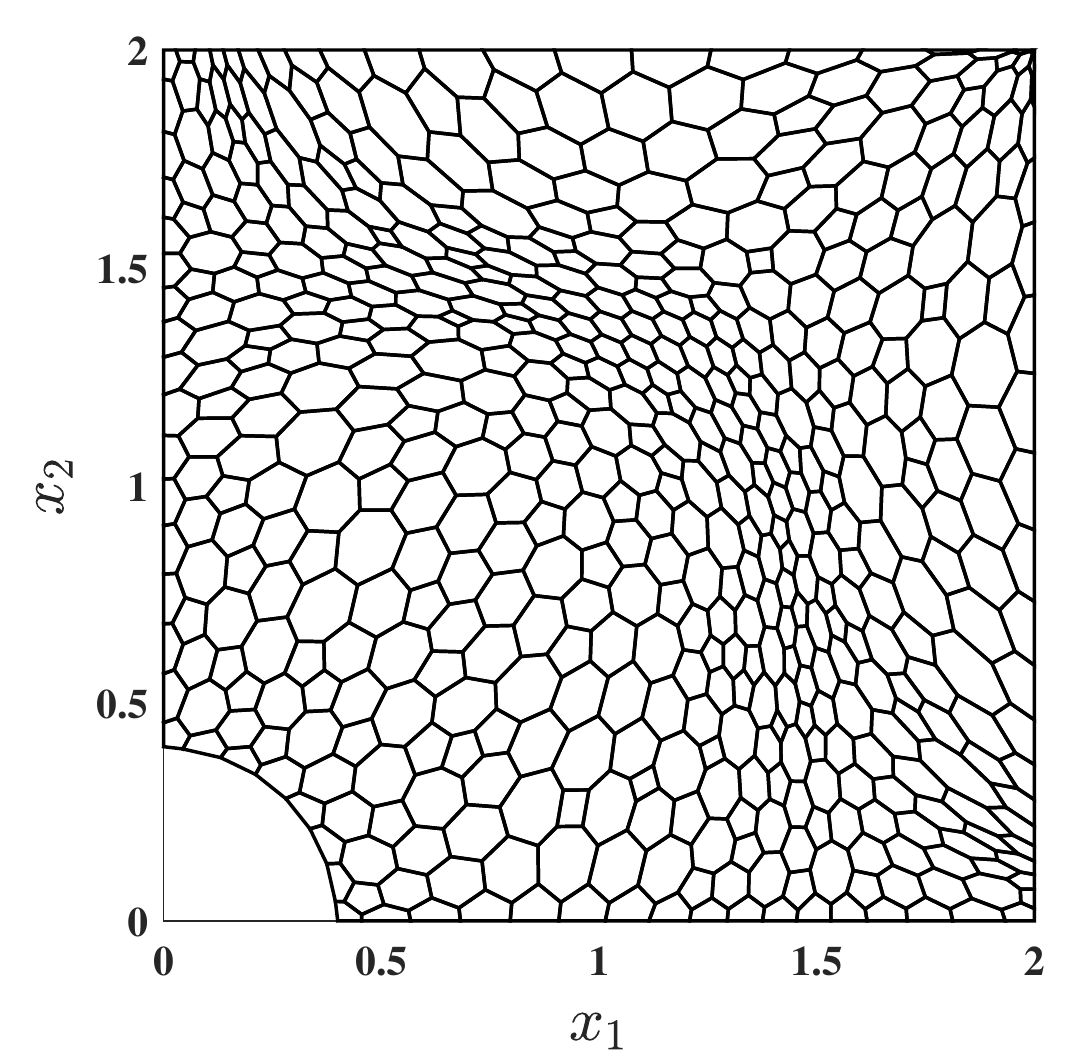}}	
	}
	\mbox{
	\subfigure[] {\label{fig:plateholemeshes_c}\includegraphics[width=0.38\linewidth]{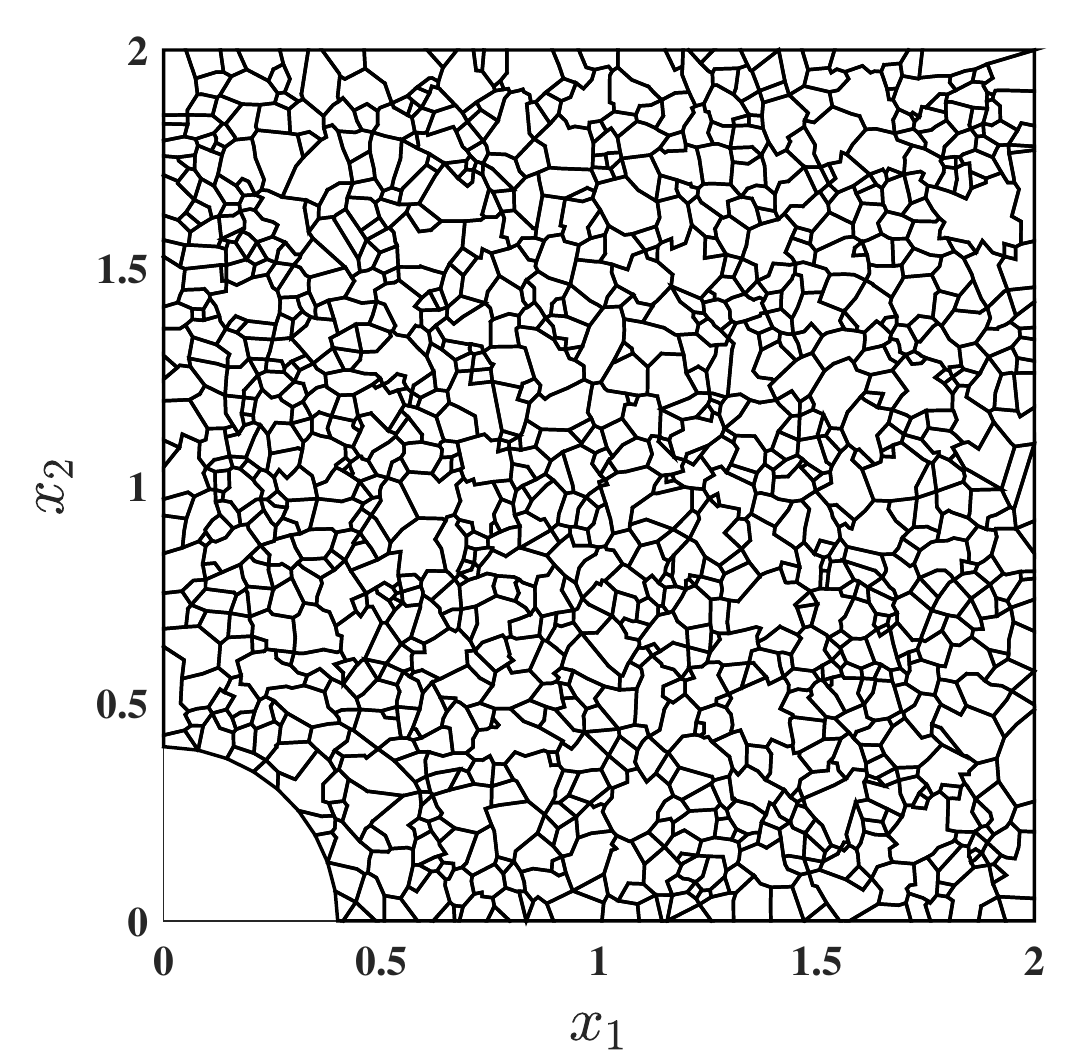}}	
	}
	\caption{Sample meshes for a quarter of the infinite plate with a circular hole. 
	(a) Regular, (b) distorted and (c) random meshes.}
	\label{fig:plateholemeshes}
\end{figure}

In this example, regular, distorted and random meshes are used. For each of them, 
\fref{fig:plateholemeshes} shows a sample mesh. The random meshes are built using 
Polylla, a polygonal mesh generator~\cite{Salinas-Hitschfeld-Ortiz-Si:2022}. Polylla generates 
meshes from any input triangulation (in this case, a Delaunay triangulation) using  
the longest-edge propagation path~\cite{Rivara:1997} and terminal-edge region concepts.

\subsubsection{Compressible elasticity}
\begin{figure}[!bth]
	\centering
	\mbox{
	\subfigure[] {\label{fig:platerates1_a}\includegraphics[width=0.5\linewidth]{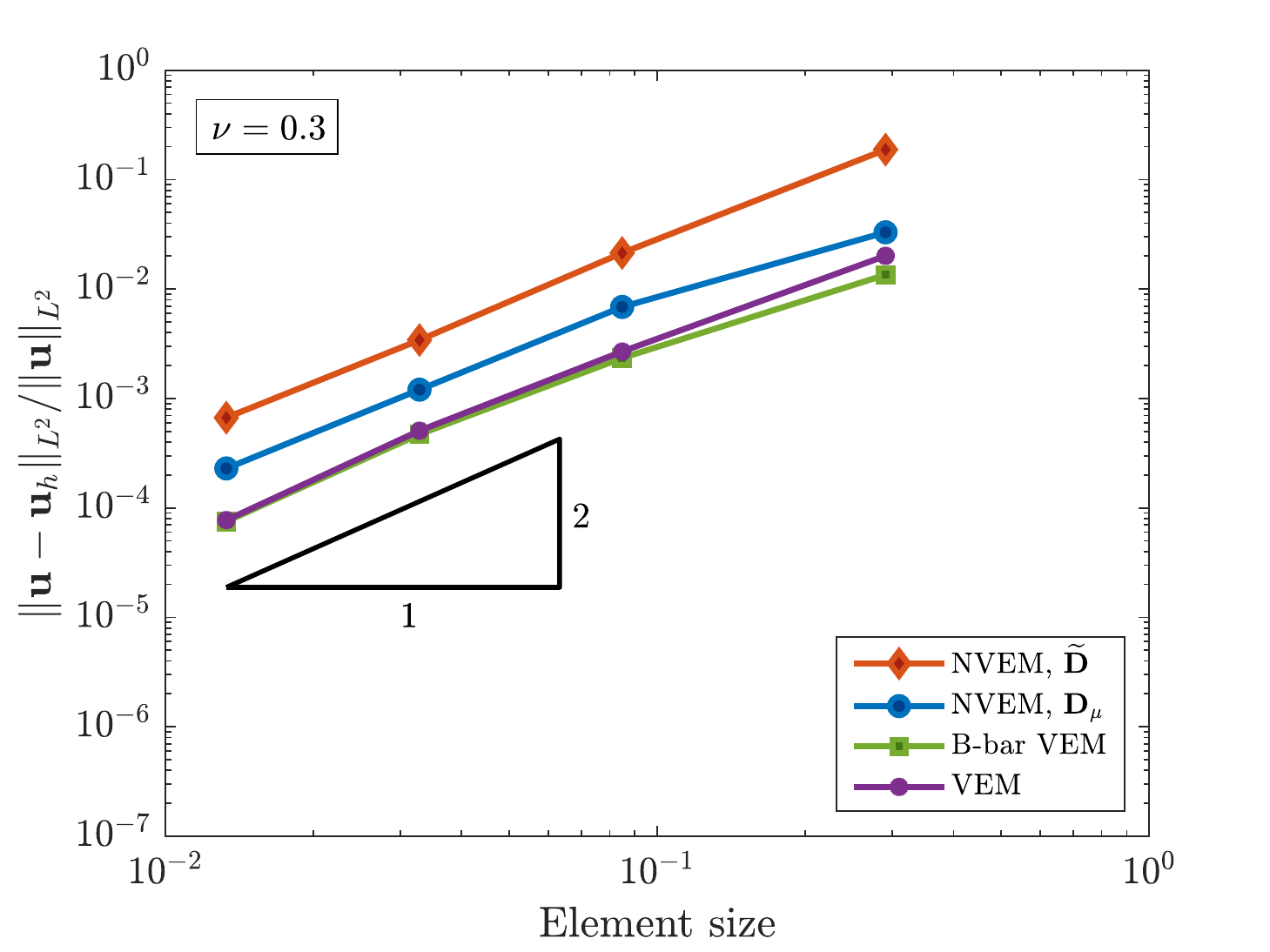}}
	\subfigure[] {\label{fig:platerates1_b}\includegraphics[width=0.5\linewidth]{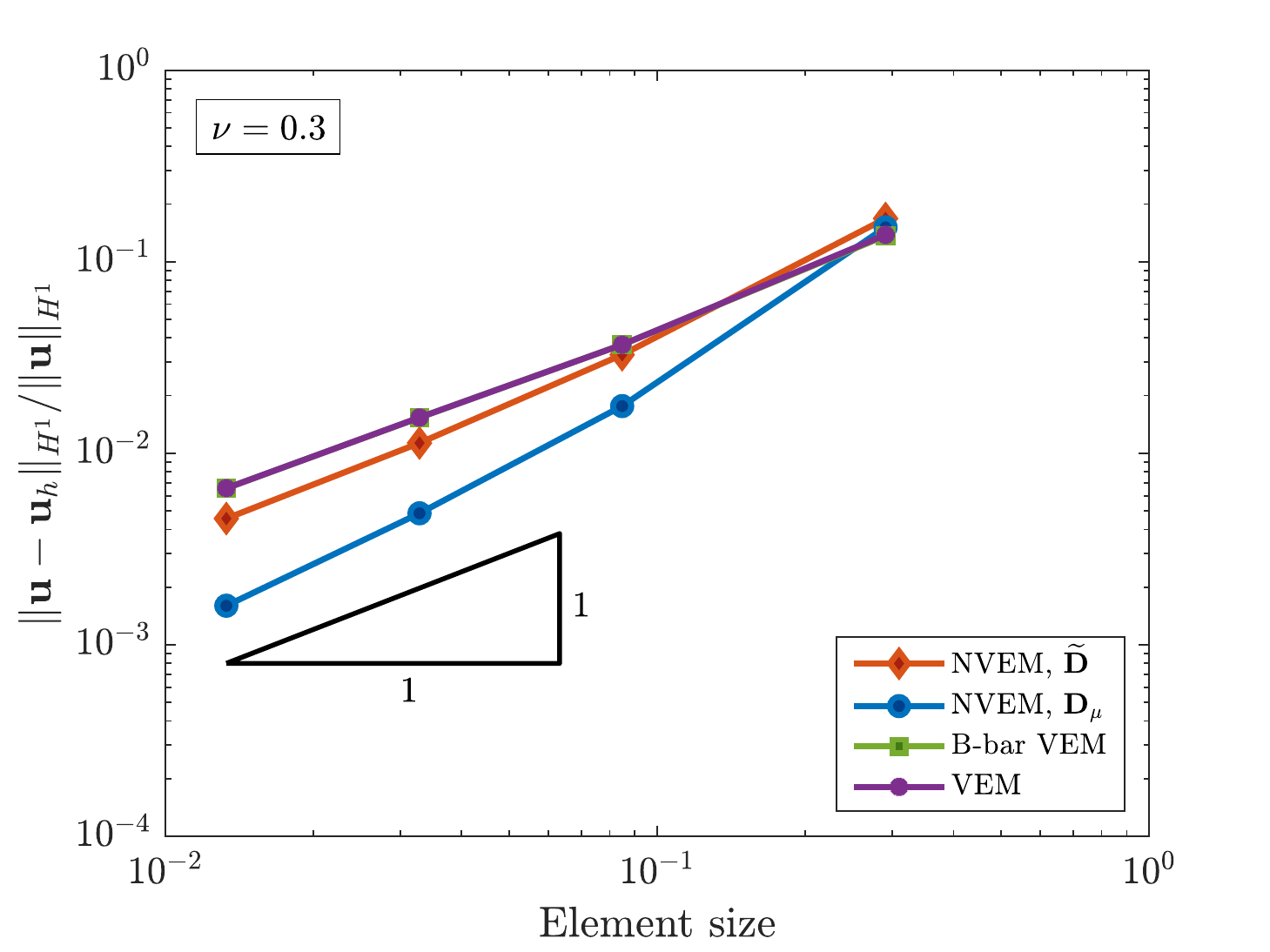}}
	}
	\subfigure[] {\label{fig:platerates1_c}\includegraphics[width=0.5\linewidth]{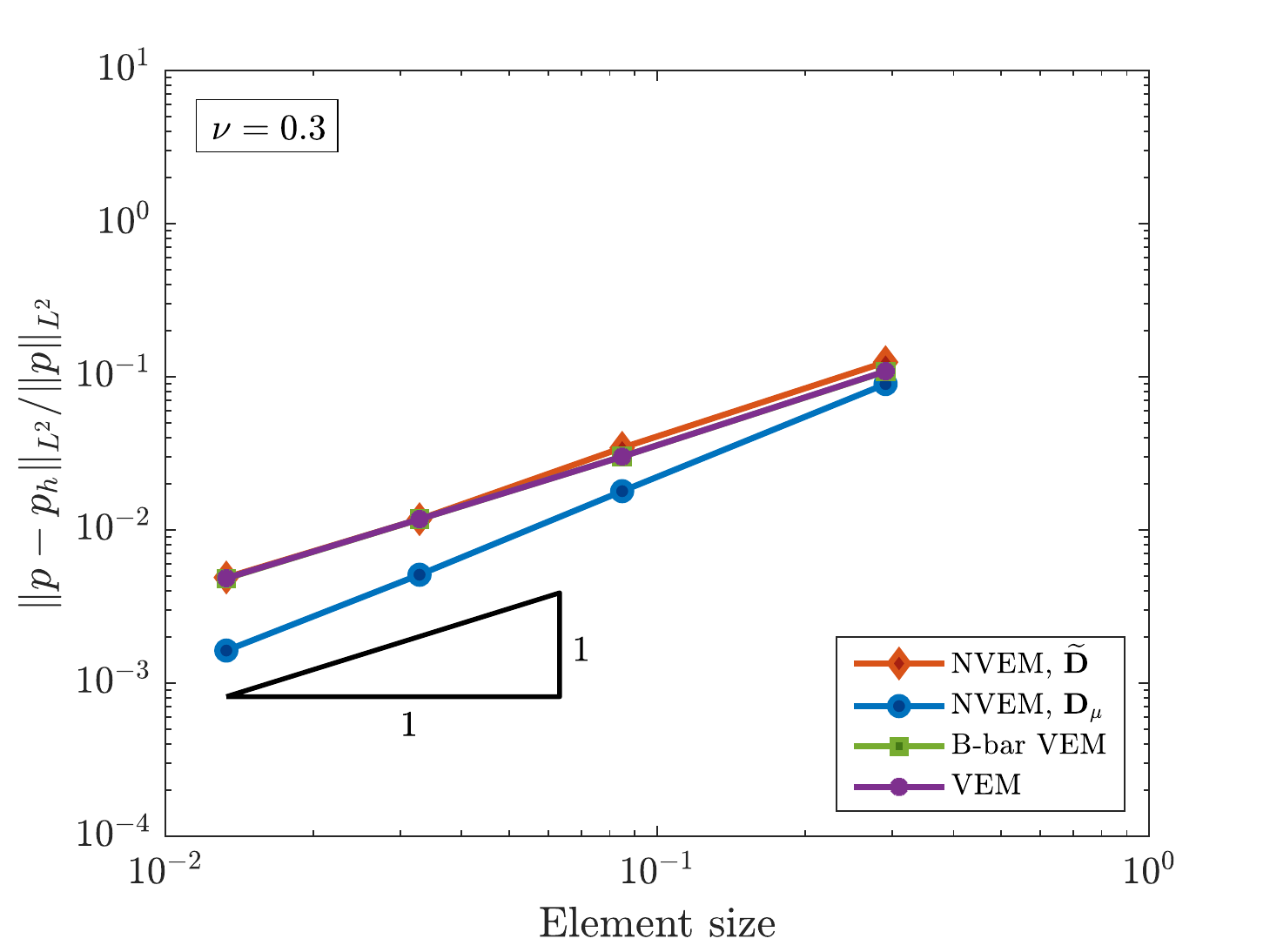}}	
	\caption{Infinite plate with a circular hole problem on regular meshes. Material parameters $E_Y=10^3$ \alejandrog{psi} and $\nu=0.3$
	for plane stress condition. Convergence rates in the (a) $L^2$ norm of the displacement error, (b) $H^1$ seminorm 
	of the displacement error and (c) $L^2$ norm of the pressure error for the VEM, B-bar VEM and NVEM.}
	\label{fig:platerates1}
\end{figure}

In terms of accuracy and convergence upon mesh refinement, \fref{fig:platerates1} reveals 
that for the compressible case on regular meshes all the methods are accurate and converge
optimally in the $L^2$ norm and $H^1$ seminorm of the displacement error, and the $L^2$ 
norm of the pressure error (rates of 2, 1 and 1, respectively). In the $L^2$ norm of
the displacement error, the VEM and the B-bar VEM are more accurate than the NVEM.
In the $H^1$ seminorm of the displacement error and the $L^2$ norm of the pressure error,
the VEM, the B-bar VEM and the NVEM with $\widetilde{\vm{D}}$ stabilization behave similar, while
the NVEM with $\vm{D}_\mu$ stabilization delivers the most accurate solution. 

\subsubsection{Nearly incompressible elasticity}
\begin{figure}[!bth]
	\centering
	\mbox{
	\subfigure[] {\label{fig:platerates2_a}\includegraphics[width=0.5\linewidth]{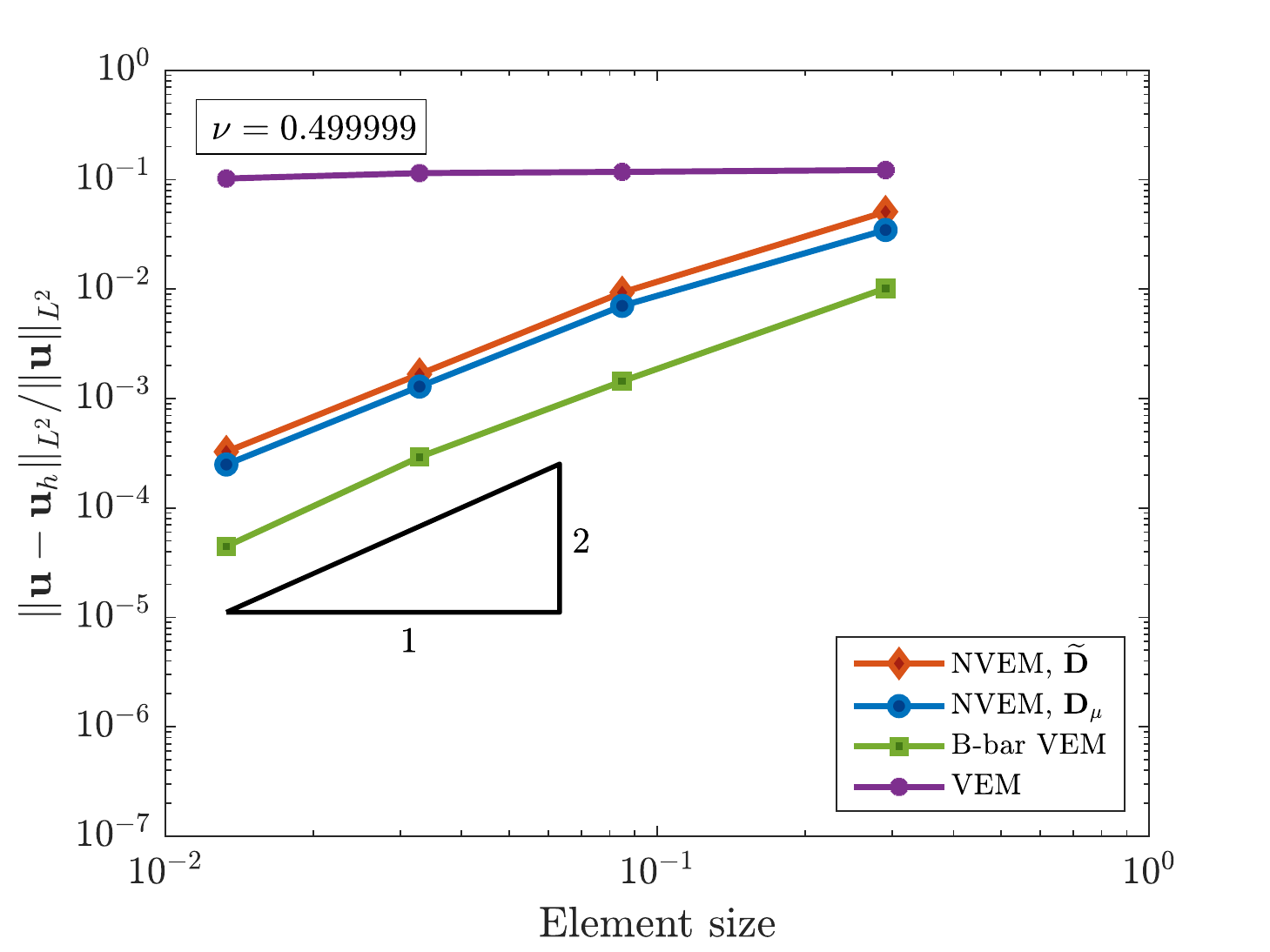}}
	\subfigure[] {\label{fig:platerates2_b}\includegraphics[width=0.5\linewidth]{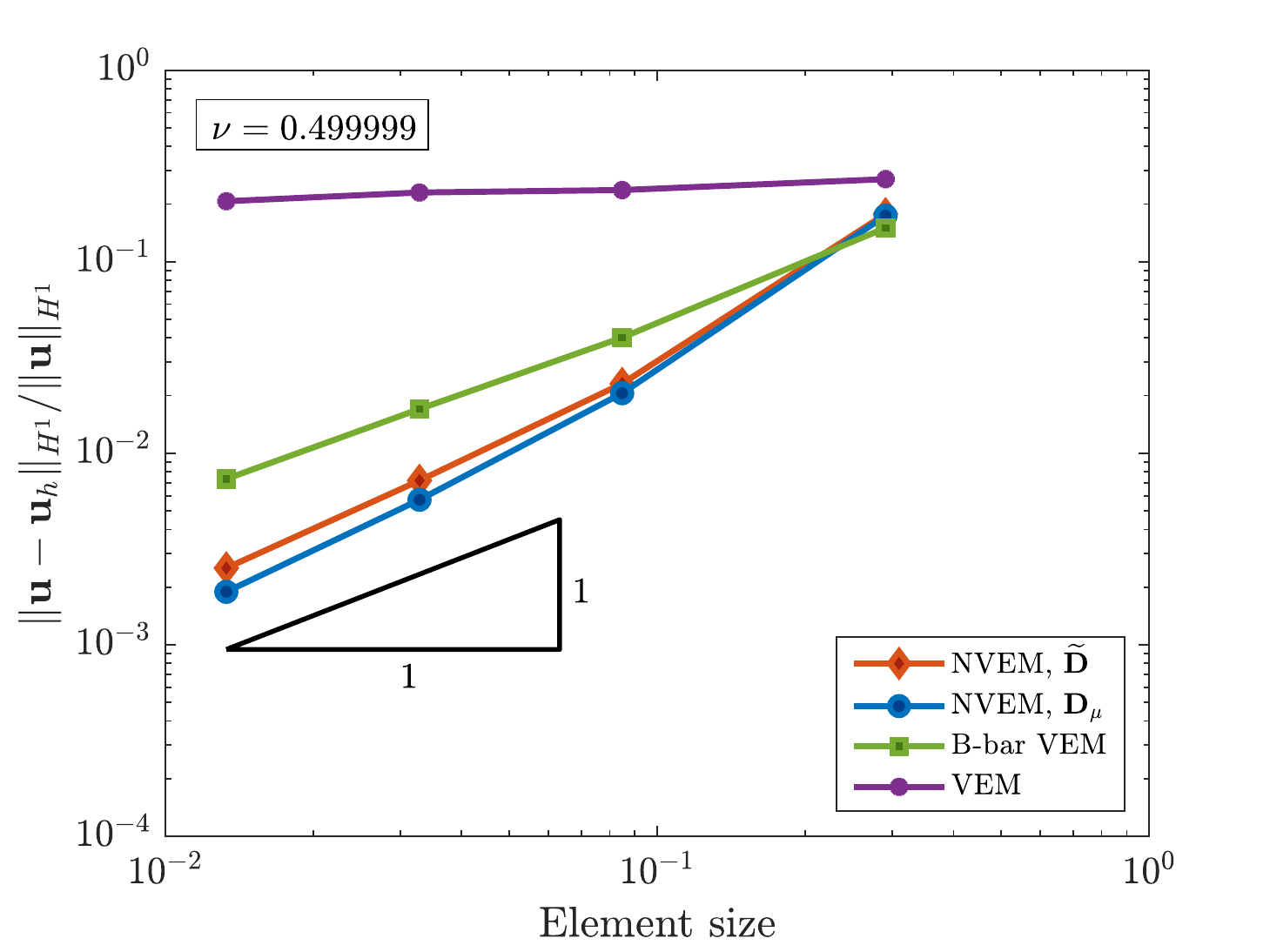}}
	}
	\subfigure[] {\label{fig:platerates2_c}\includegraphics[width=0.5\linewidth]{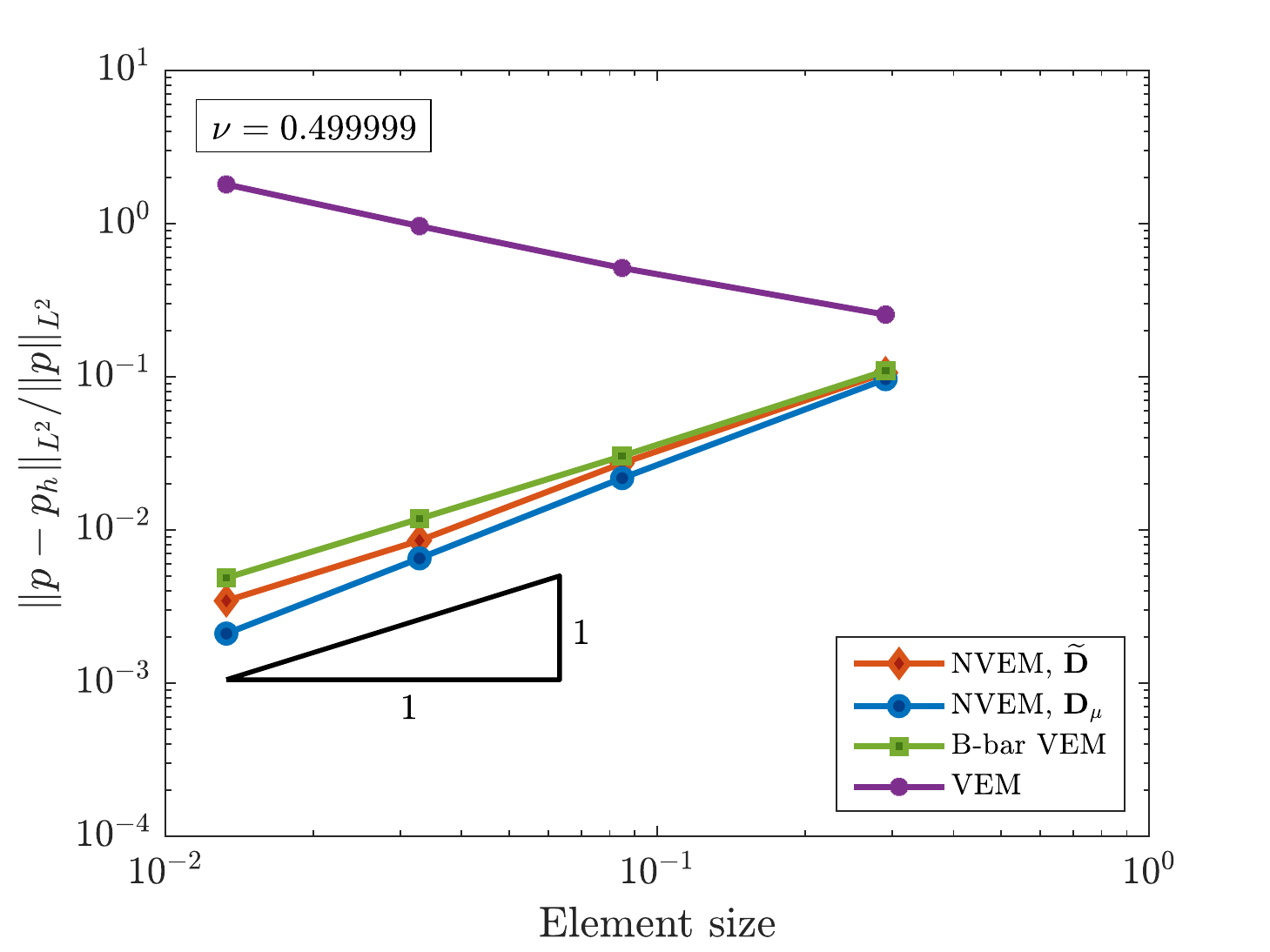}}	
	\caption{Infinite plate with a circular hole problem (regular meshes). Material parameters $E_Y=10^3$ \alejandrog{psi} and $\nu=0.499999$
	for plane strain condition. Convergence rates in the (a) $L^2$ norm of the displacement error, (b) $H^1$ seminorm 
	of the displacement error and (c) $L^2$ norm of the pressure error for the VEM, B-bar VEM and NVEM.}
	\label{fig:platerates2}
\end{figure}

\begin{figure}[!bth]
	\centering
	\mbox{
	\subfigure[] {\label{fig:platerates3_a}\includegraphics[width=0.5\linewidth]{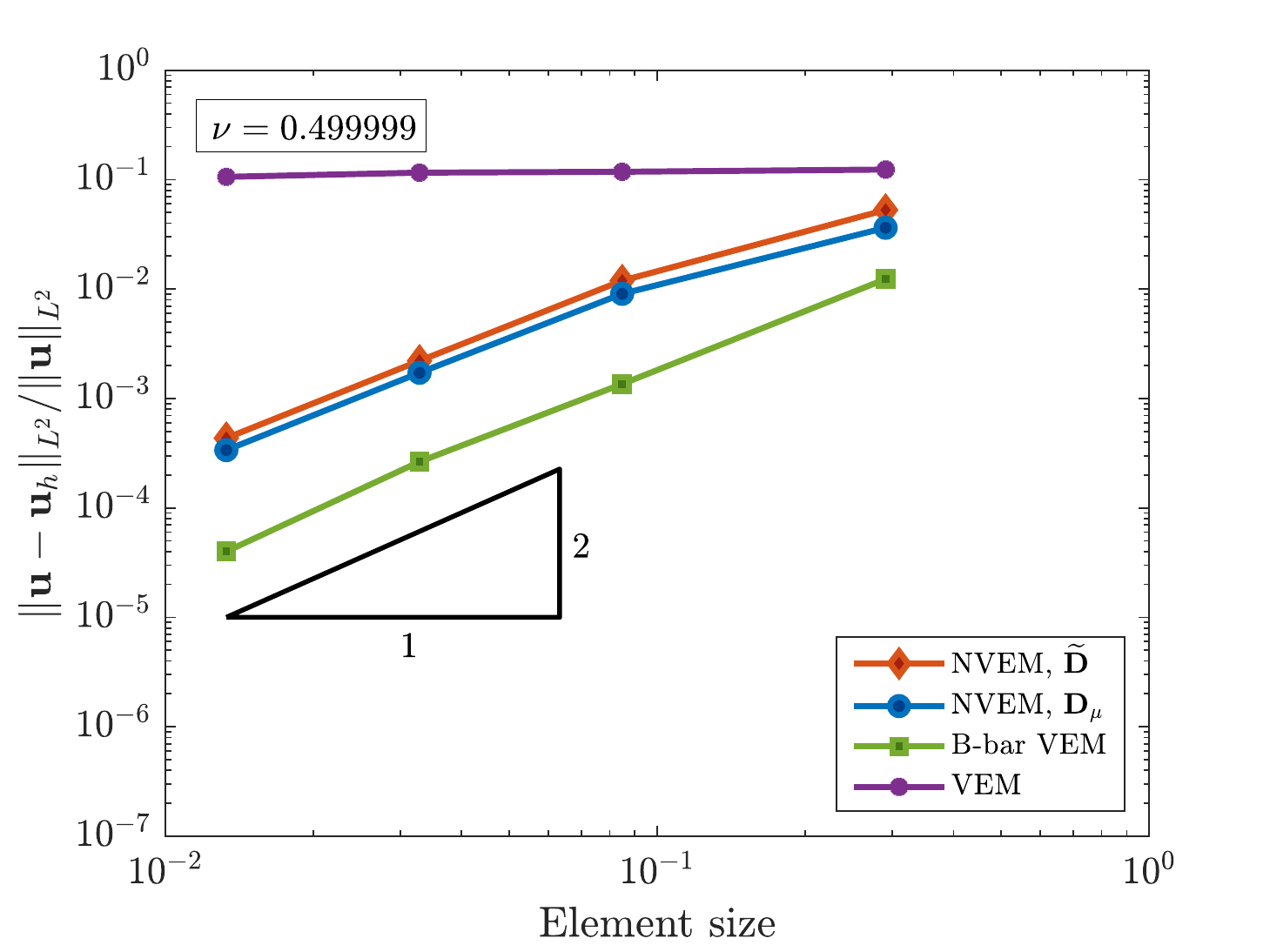}}
	\subfigure[] {\label{fig:platerates3_b}\includegraphics[width=0.5\linewidth]{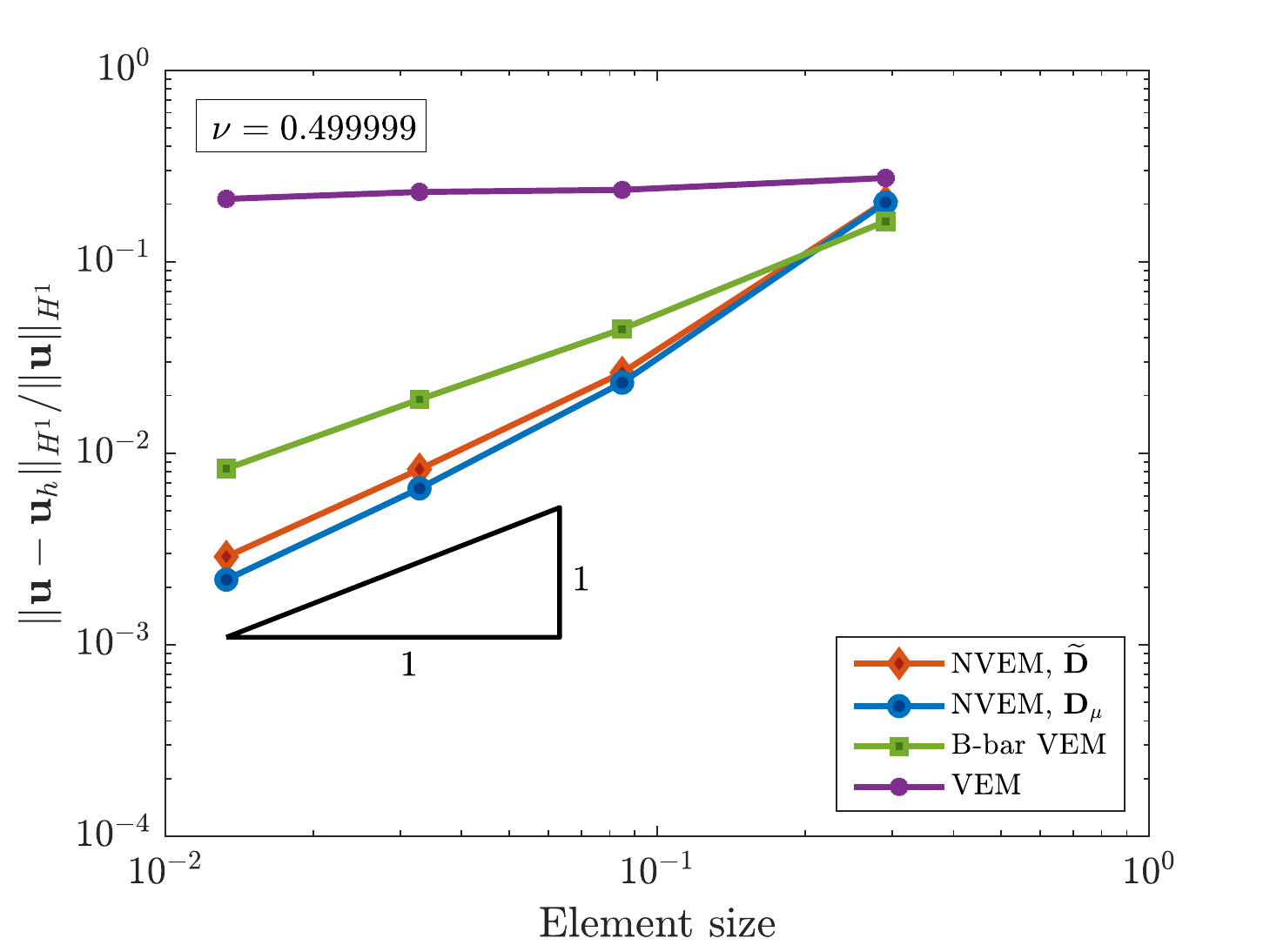}}
	}
	\subfigure[] {\label{fig:platerates3_c}\includegraphics[width=0.5\linewidth]{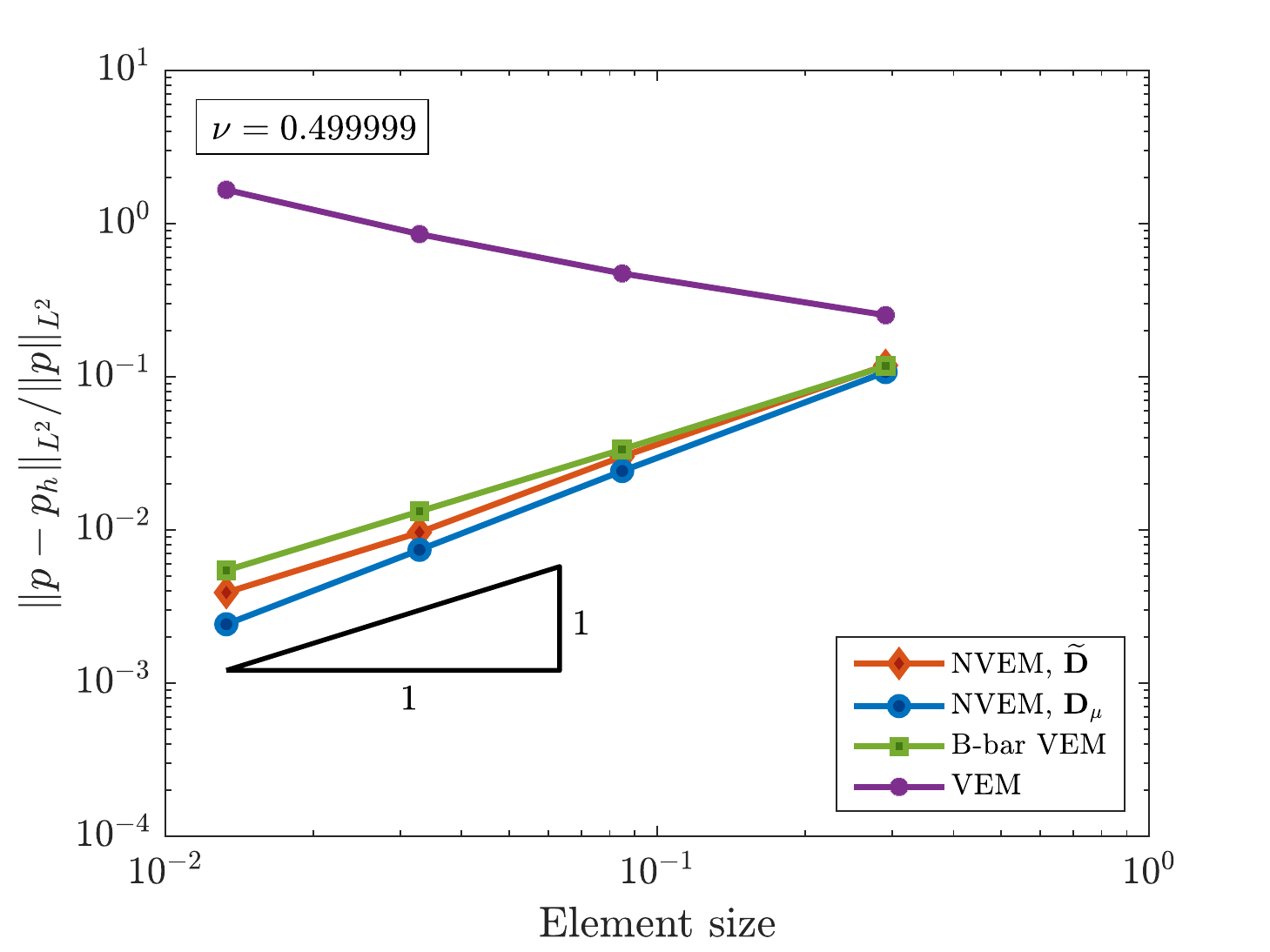}}	
	\caption{Infinite plate with a circular hole problem (distorted meshes). Material parameters $E_Y=10^3$ \alejandrog{psi} and $\nu=0.499999$
	for plane strain condition. Convergence rates in the (a) $L^2$ norm of the displacement error, (b) $H^1$ seminorm 
	of the displacement error and (c) $L^2$ norm of the pressure error for the VEM, B-bar VEM and NVEM.}
	\label{fig:platerates3}
\end{figure}

For the nearly incompressible case, regular, distorted and random meshes are considered.  
On regular and distorted meshes, the B-bar VEM and the NVEM are accurate and optimally 
convergent in the three norms (see Figs. \ref{fig:platerates2} and \ref{fig:platerates3}),
whereas, as expected, the VEM is severely inaccurate due to volumetric locking.
In comparison with the B-bar VEM, Figs. \ref{fig:platerates2} and \ref{fig:platerates3} 
show that the NVEM is slightly more accurate in the $H^1$ seminorm of the displacement 
error and the $L^2$ norm of the pressure error, whereas less accurate in 
the $L^2$ norm of the displacement error. We also observe that between the two 
stabilized NVEM, the NVEM with $\vm{D}_\mu$ stabilization is slightly more accurate. 

\begin{figure}[!bth]
	\centering
	\mbox{
	\subfigure[] {\label{fig:platerates4_a}\includegraphics[width=0.5\linewidth]{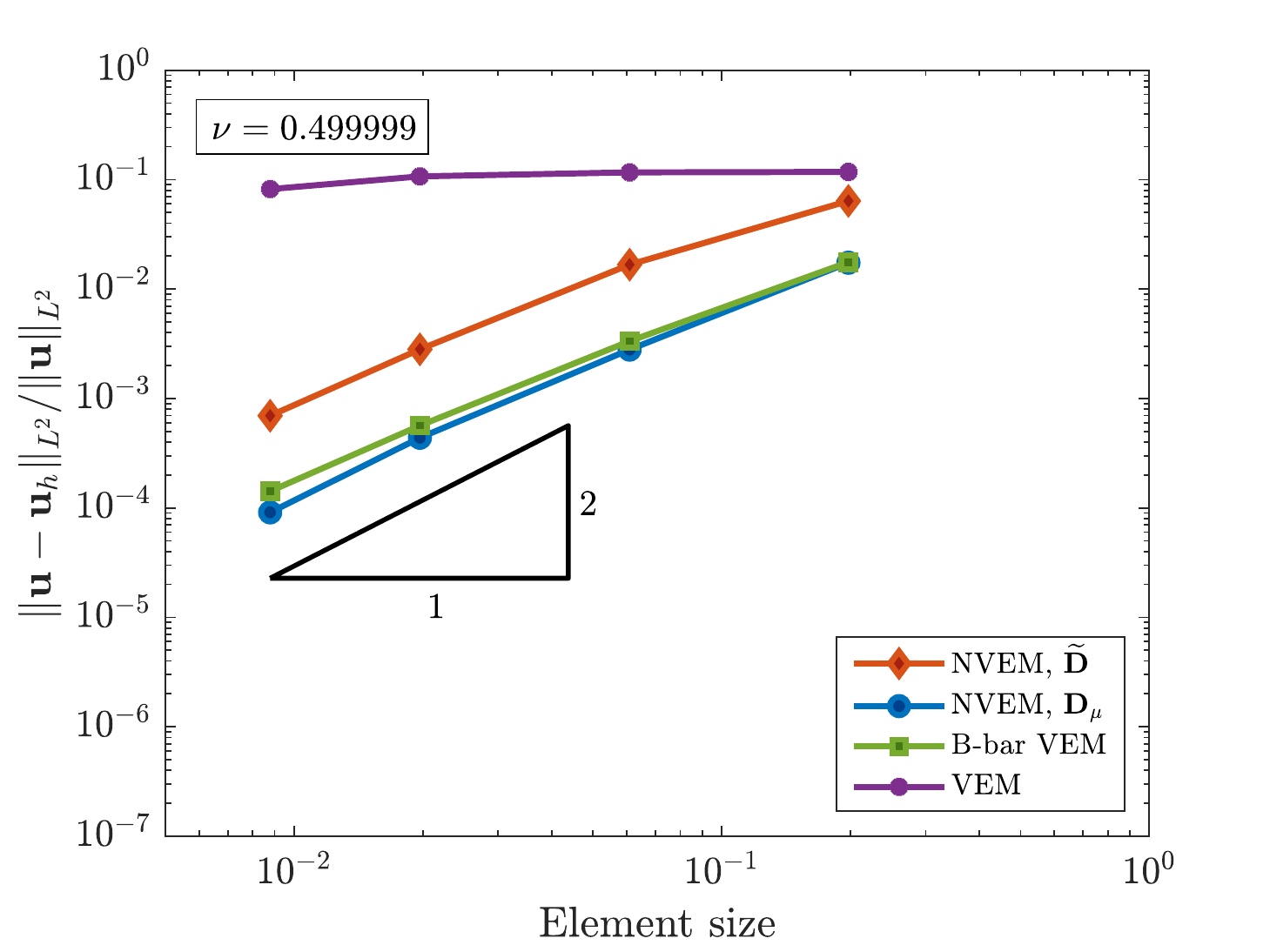}}
	\subfigure[] {\label{fig:platerates4_b}\includegraphics[width=0.5\linewidth]{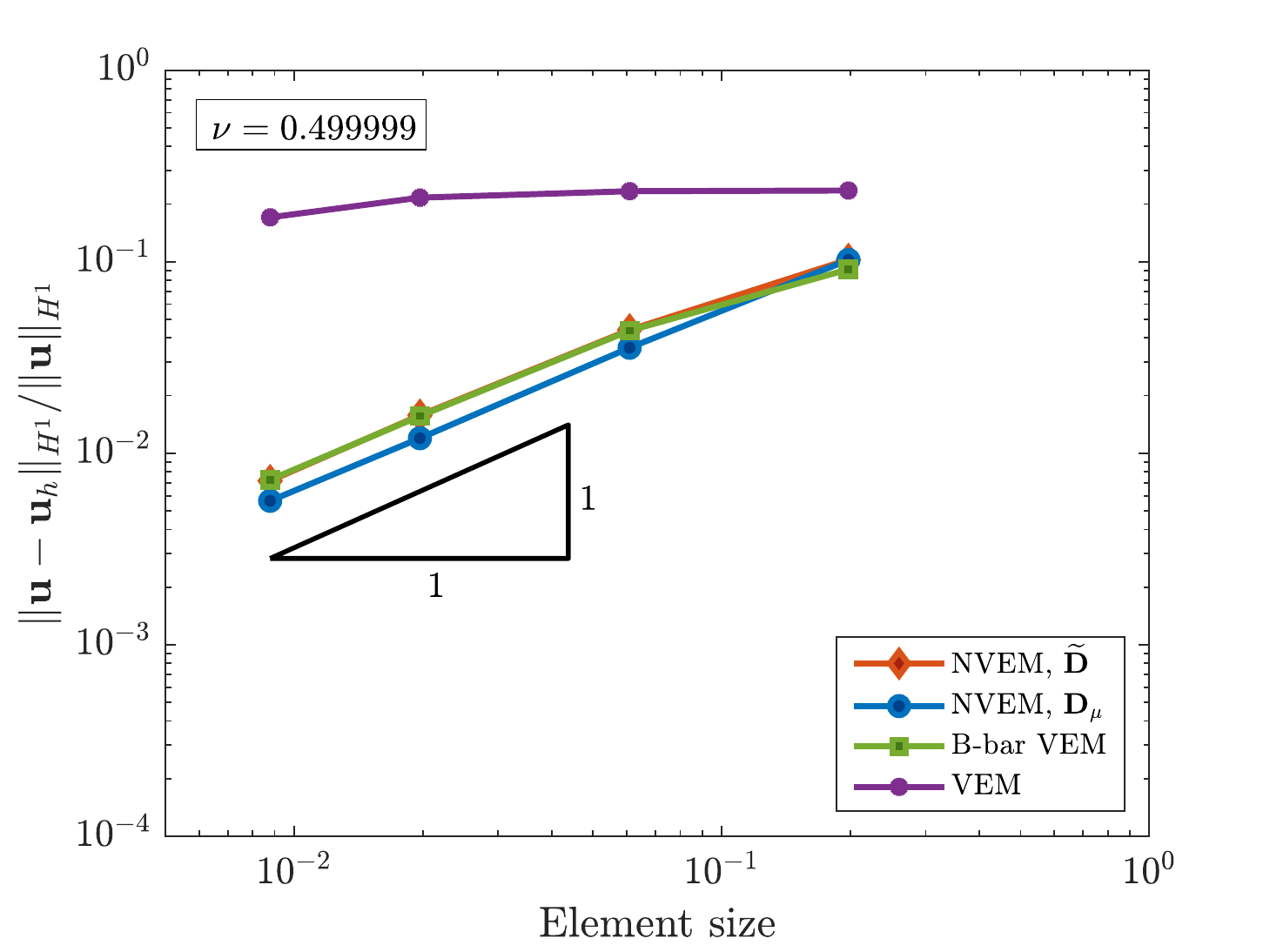}}
	}
	\subfigure[] {\label{fig:platerates4_c}\includegraphics[width=0.5\linewidth]{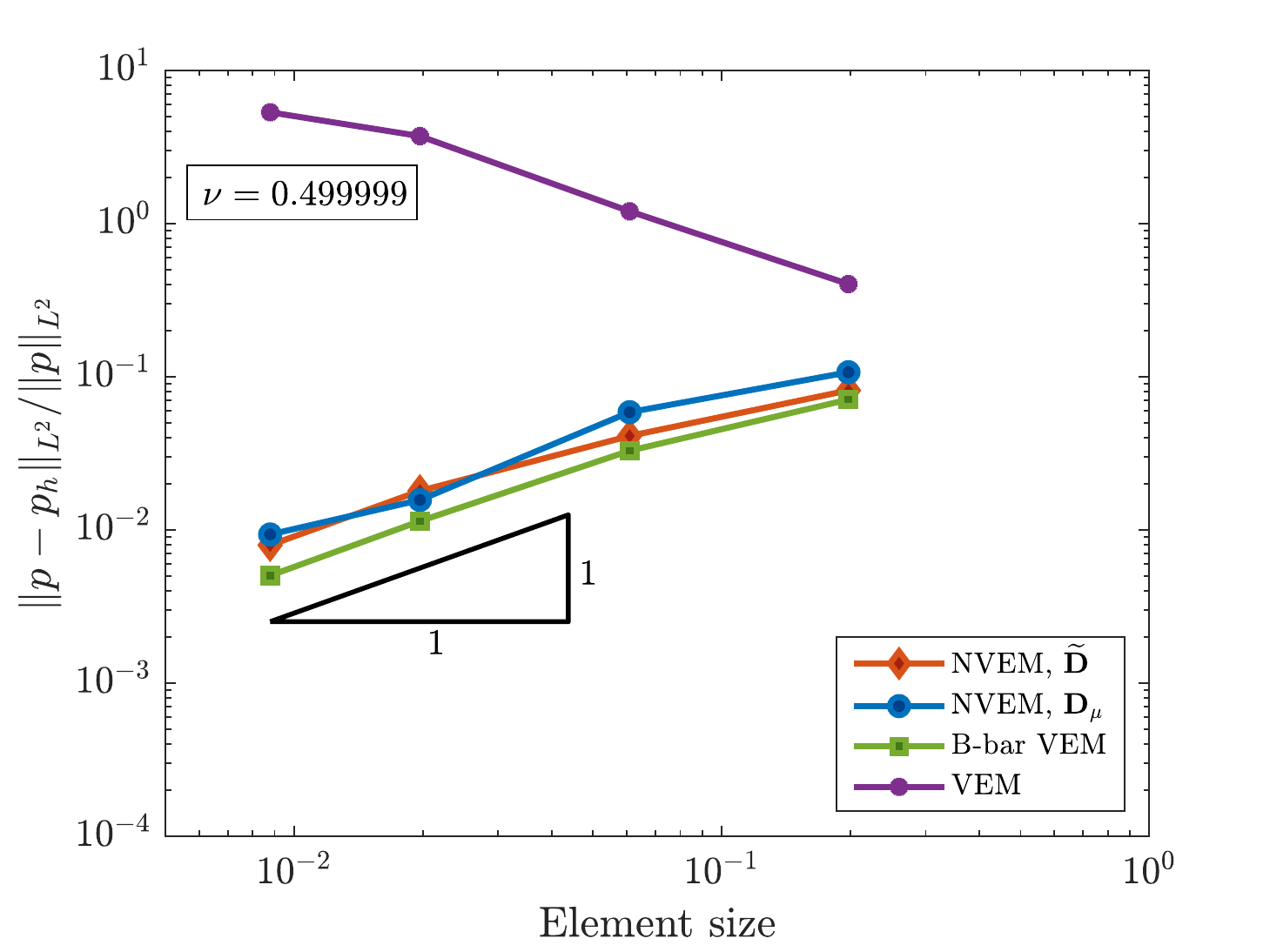}}	
	\caption{Infinite plate with a circular hole problem (random meshes). Material parameters $E_Y=10^3$ \alejandrog{psi} and $\nu=0.499999$
	for plane strain condition. Convergence rates in the (a) $L^2$ norm of the displacement error, (b) $H^1$ seminorm 
	of the displacement error and (c) $L^2$ norm of the pressure error for the VEM, B-bar VEM and NVEM.}
	\label{fig:platerates4}
\end{figure}

On random meshes, the B-bar VEM and the NVEM are accurate and optimally 
convergent in the three norms (\fref{fig:platerates4}), whereas the VEM is severely 
inaccurate due to volumetric locking. \fref{fig:platerates4_a} reveals that 
in the $L^2$ norm of the displacement error the NVEM with $\widetilde{\vm{D}}$ 
stabilization is less accurate than the B-bar VEM,
whereas slightly more accurate when $\vm{D}_\mu$ stabilization is used. Regarding the
$H^1$ seminorm of the displacement error, \fref{fig:platerates4_b} shows that the
B-bar VEM and the NVEM with $\widetilde{\vm{D}}$ stabilization behave similar, while the
NVEM with $\vm{D}_\mu$ stabilization is the most accurate. With respect to the 
$L^2$ norm of the pressure error, the B-bar VEM is slightly more accurate than 
the NVEM approach (\fref{fig:platerates4_c}).

\begin{figure}[!bth]
	\centering
	\mbox{
	\subfigure[] {\label{fig:platecontouru_a}\includegraphics[width=0.43\linewidth]{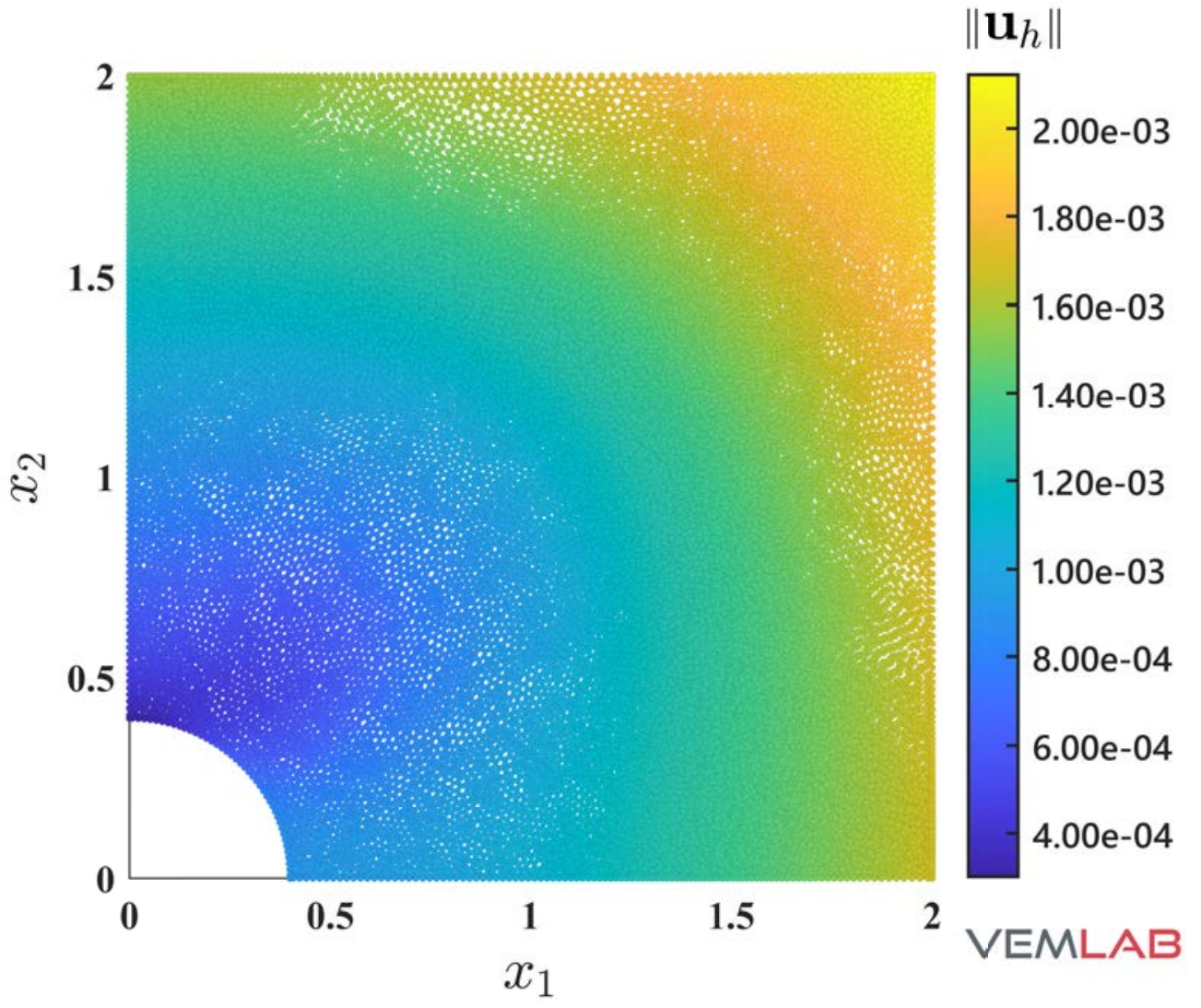}}
	\subfigure[] {\label{fig:platecontouru_b}\includegraphics[width=0.43\linewidth]{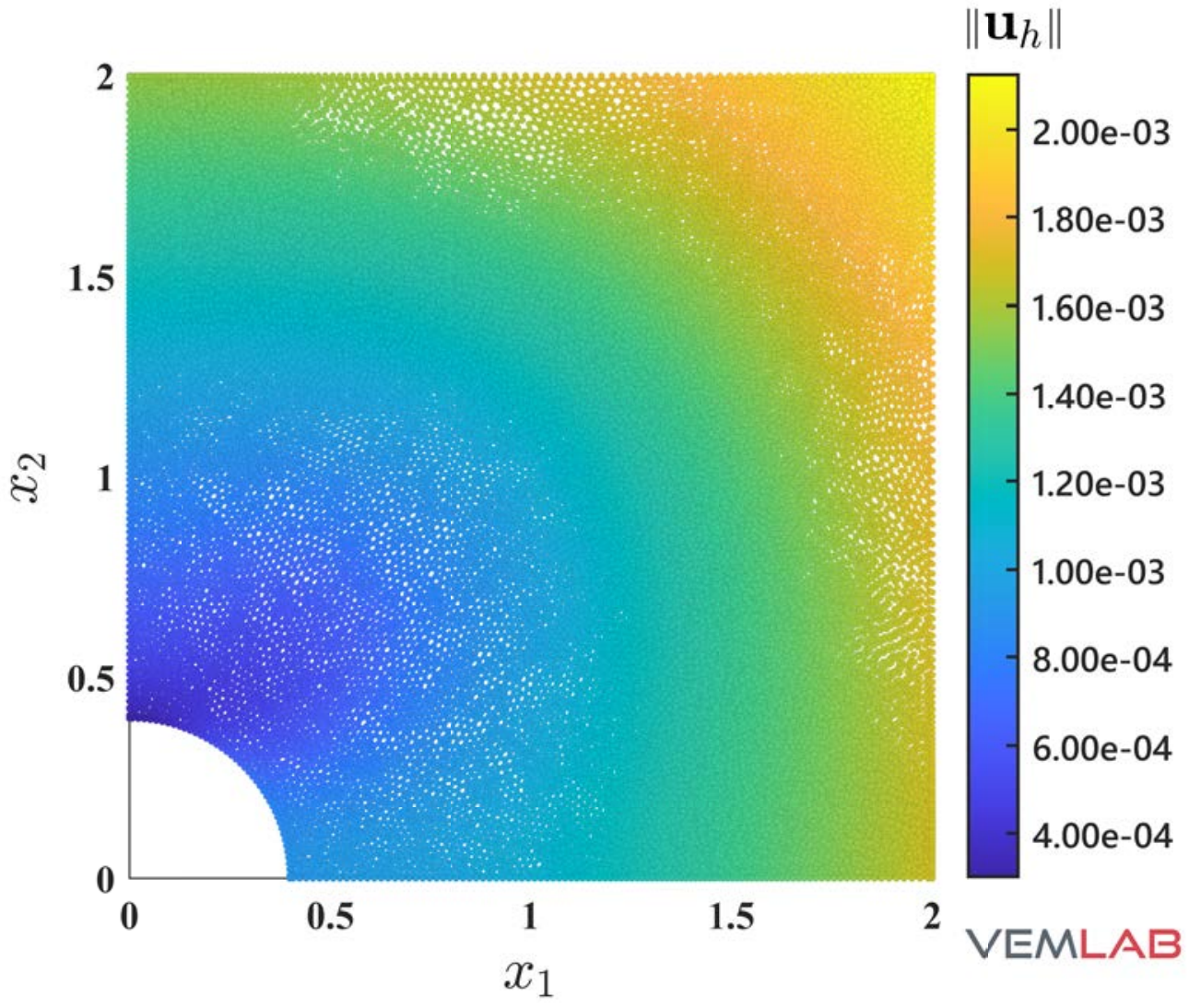}}
	}
	\subfigure[] {\label{fig:platecontouru_c}\includegraphics[width=0.43\linewidth]{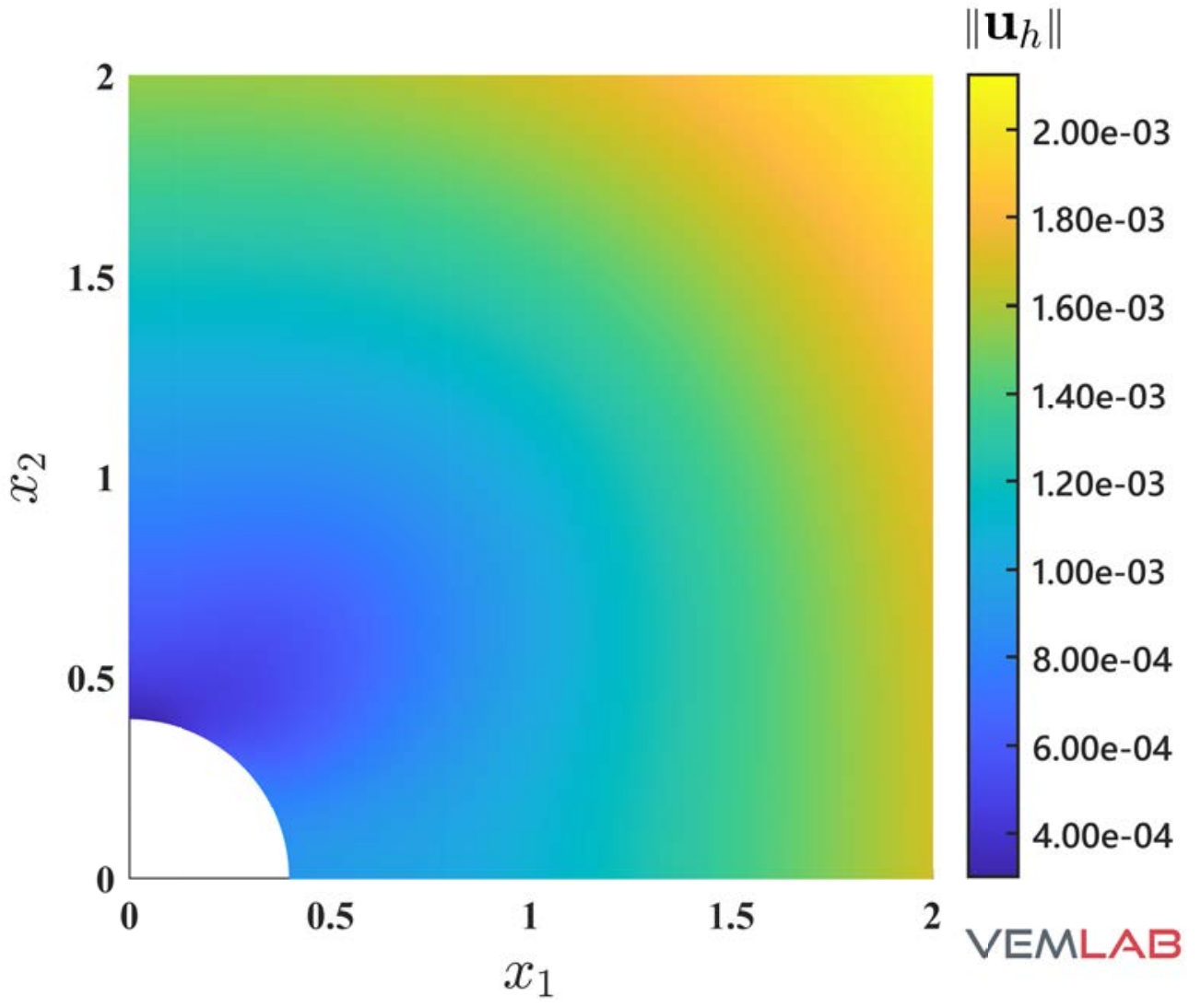}}
	\caption{Infinite plate with a circular hole problem. 
	Material parameters $E_Y=10^3$ \alejandro{psi} and $\nu=0.499999$ for plane strain condition. 
	Plots of the norm of the displacement field
	solution ($\|\vm{u}_h\|$) \alejandro{in inches} for the (a) NVEM ($\tilde{\vm{D}}$ stabilization), 
	(b) NVEM ($\vm{D}_{\mu}$ stabilization) and (c) B-bar VEM approaches.}
	\label{fig:platecontouru}
\end{figure}

\begin{figure}[!bth]
	\centering
	\mbox{
	\subfigure[] {\label{fig:platecontourp_a}\includegraphics[width=0.43\linewidth]{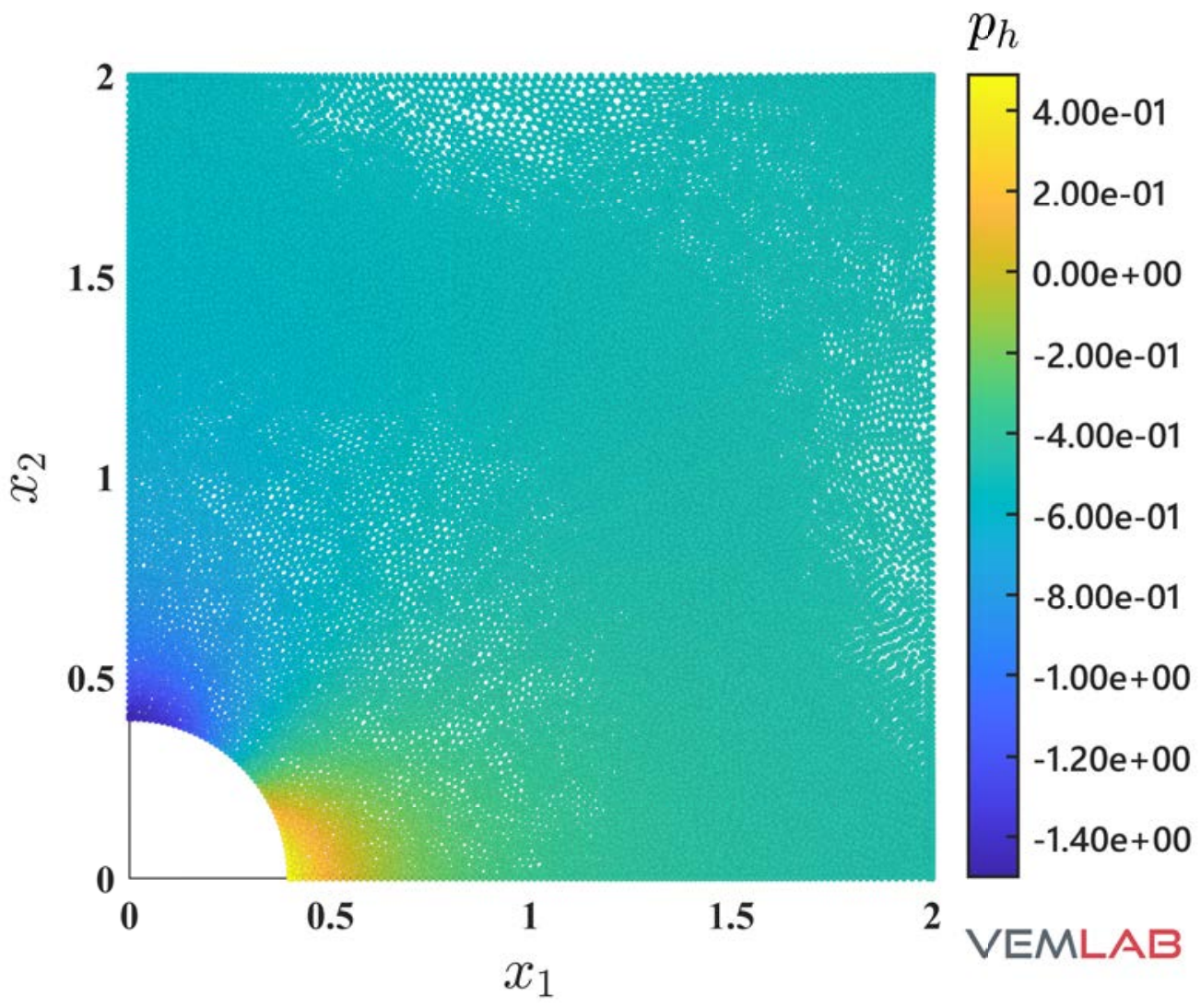}}
	\subfigure[] {\label{fig:platecontourp_b}\includegraphics[width=0.43\linewidth]{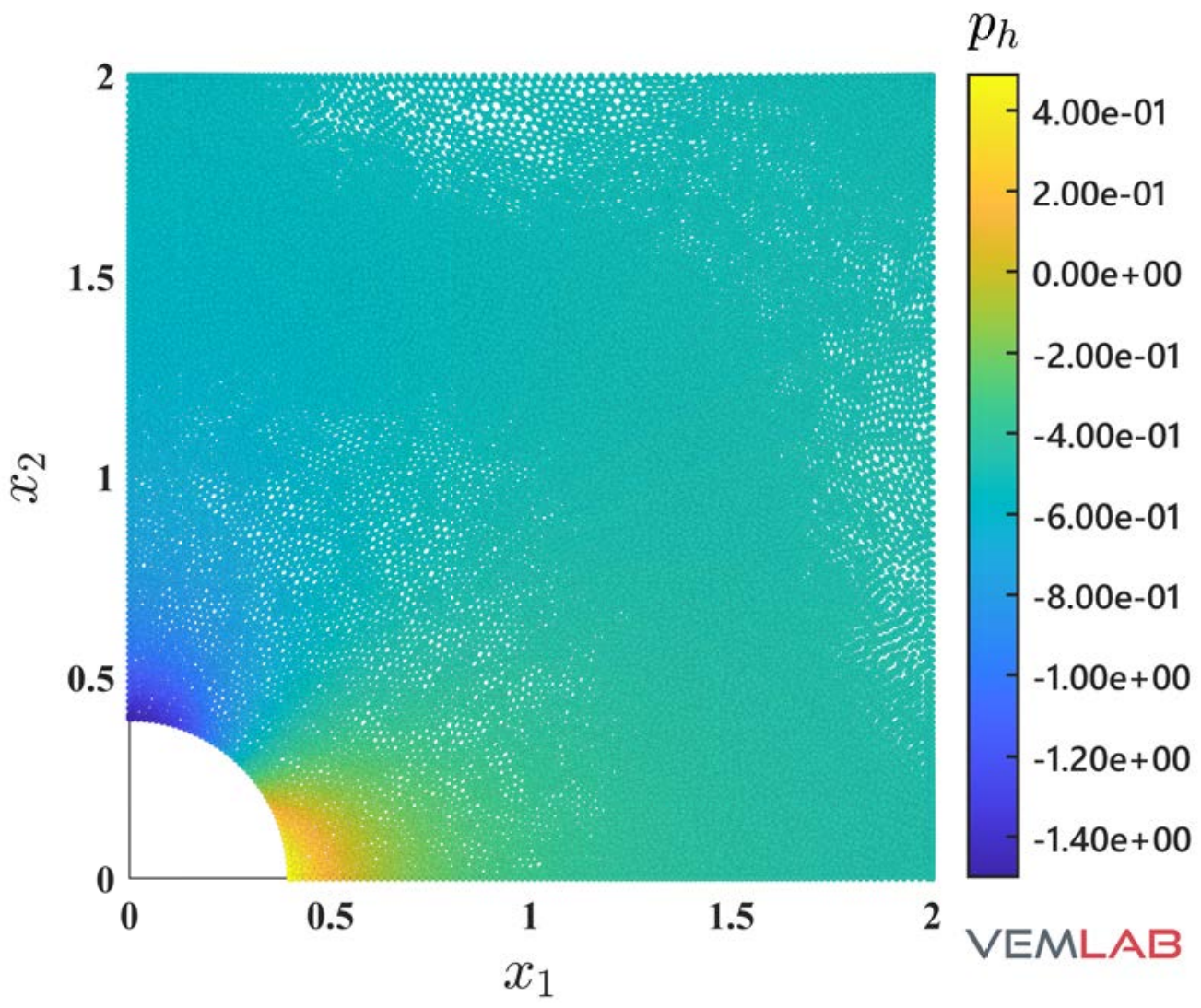}}
	}
	\subfigure[] {\label{fig:platecontourp_c}\includegraphics[width=0.43\linewidth]{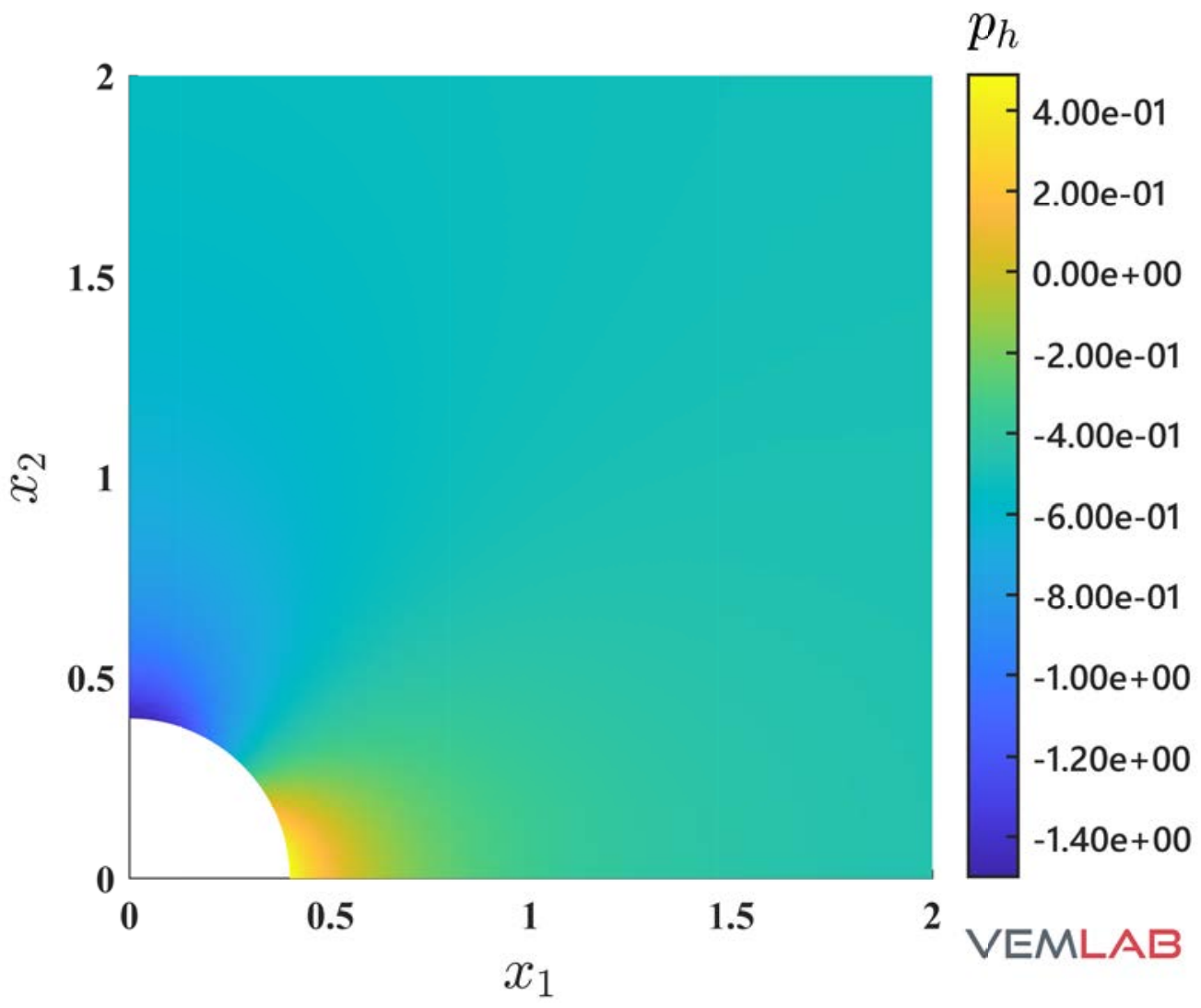}}
	\caption{Infinite plate with a circular hole problem. 
	Material parameters $E_Y=10^3$ \alejandro{psi} and $\nu=0.499999$ for plane strain condition. 
	Plots of the pressure field
	solution ($p_h$) \alejandro{in psi} for the (a) NVEM ($\tilde{\vm{D}}$ stabilization), 
	(b) NVEM ($\vm{D}_{\mu}$ stabilization) and (c) B-bar VEM approaches.}
	\label{fig:platecontourp}
\end{figure}

\begin{figure}[!bth]
	\centering
	\mbox{
	\subfigure[] {\label{fig:platecontourvm_a}\includegraphics[width=0.43\linewidth]{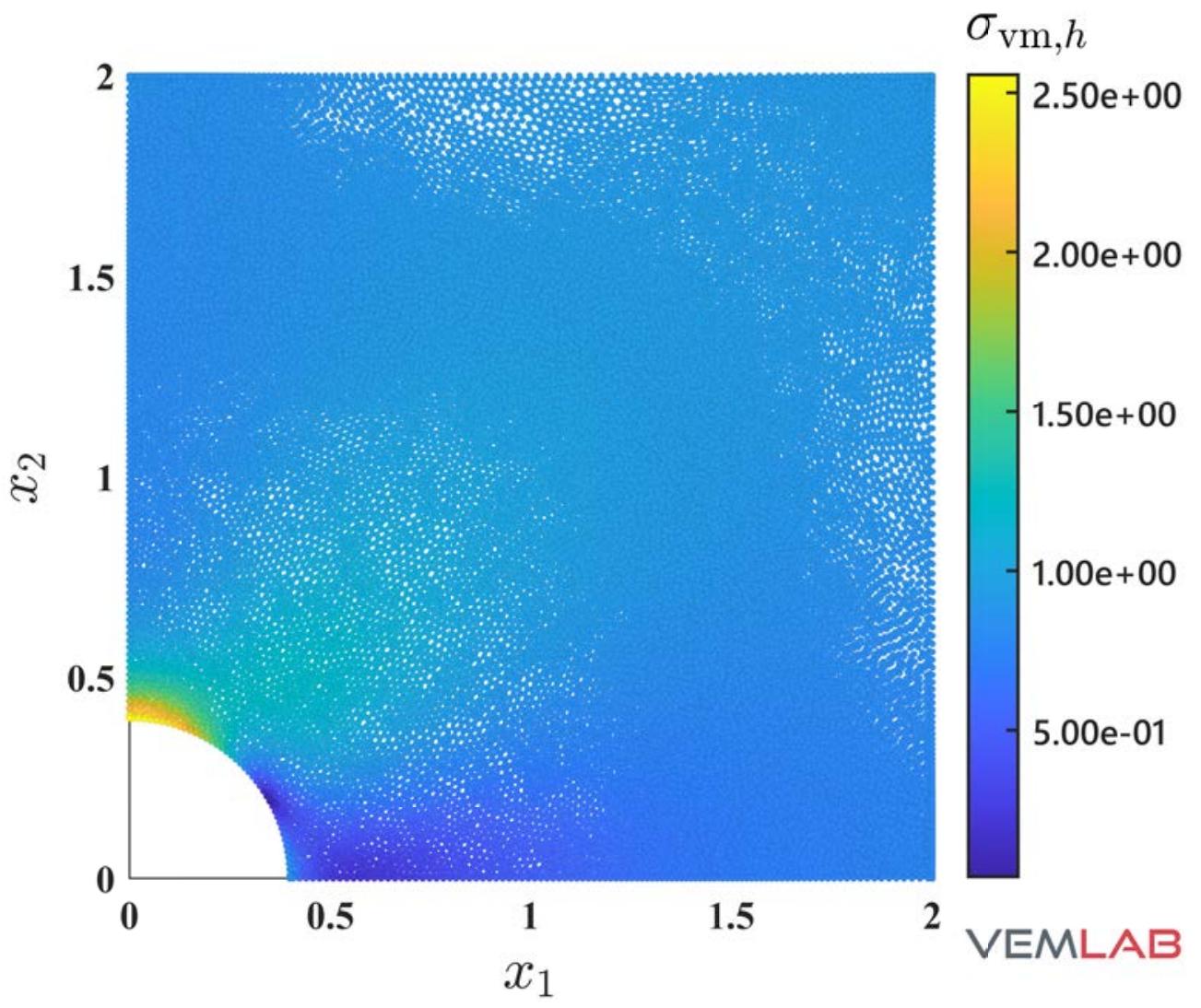}}
	\subfigure[] {\label{fig:platecontourvm_b}\includegraphics[width=0.43\linewidth]{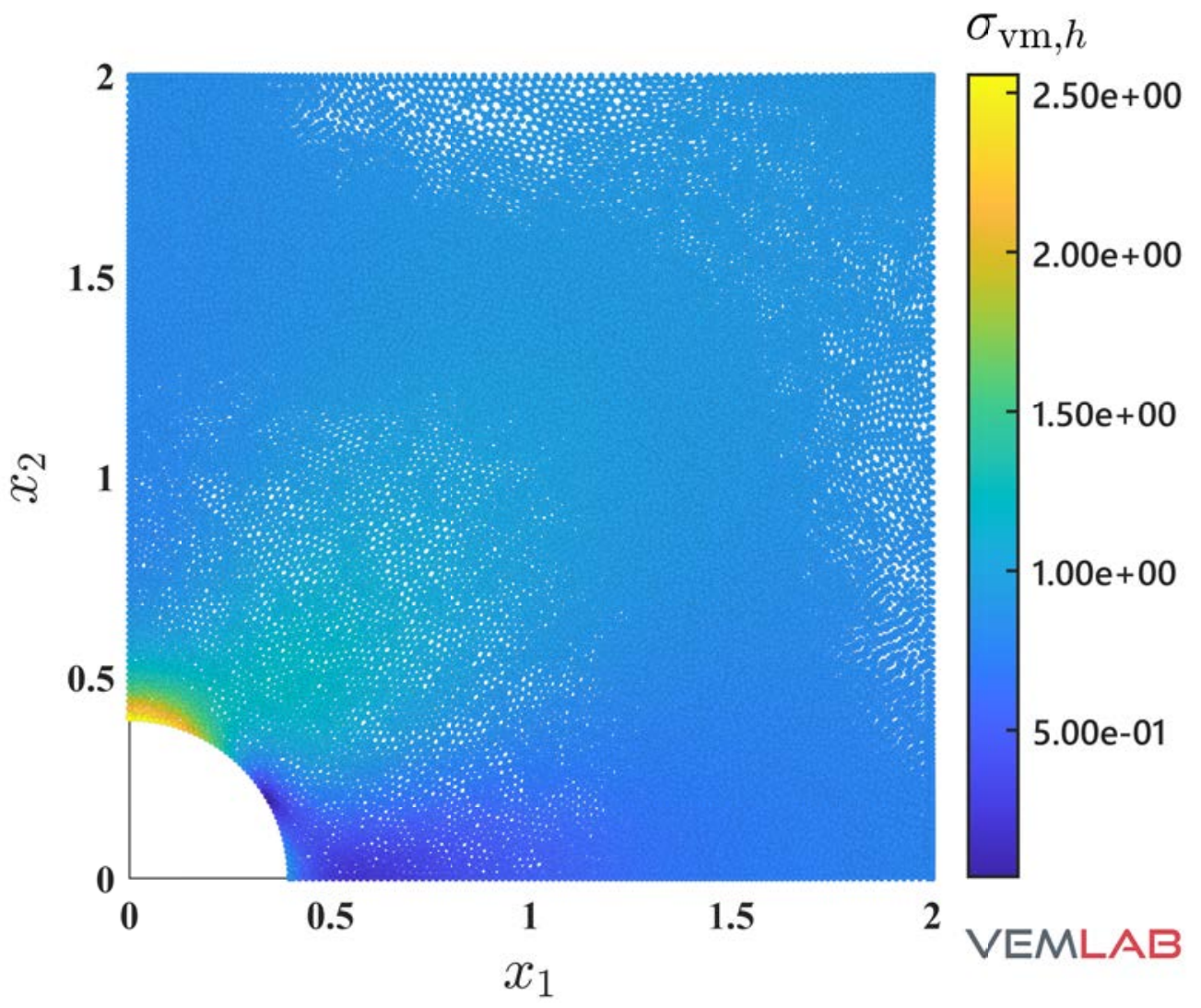}}
	}
	\subfigure[] {\label{fig:platecontourvm_c}\includegraphics[width=0.43\linewidth]{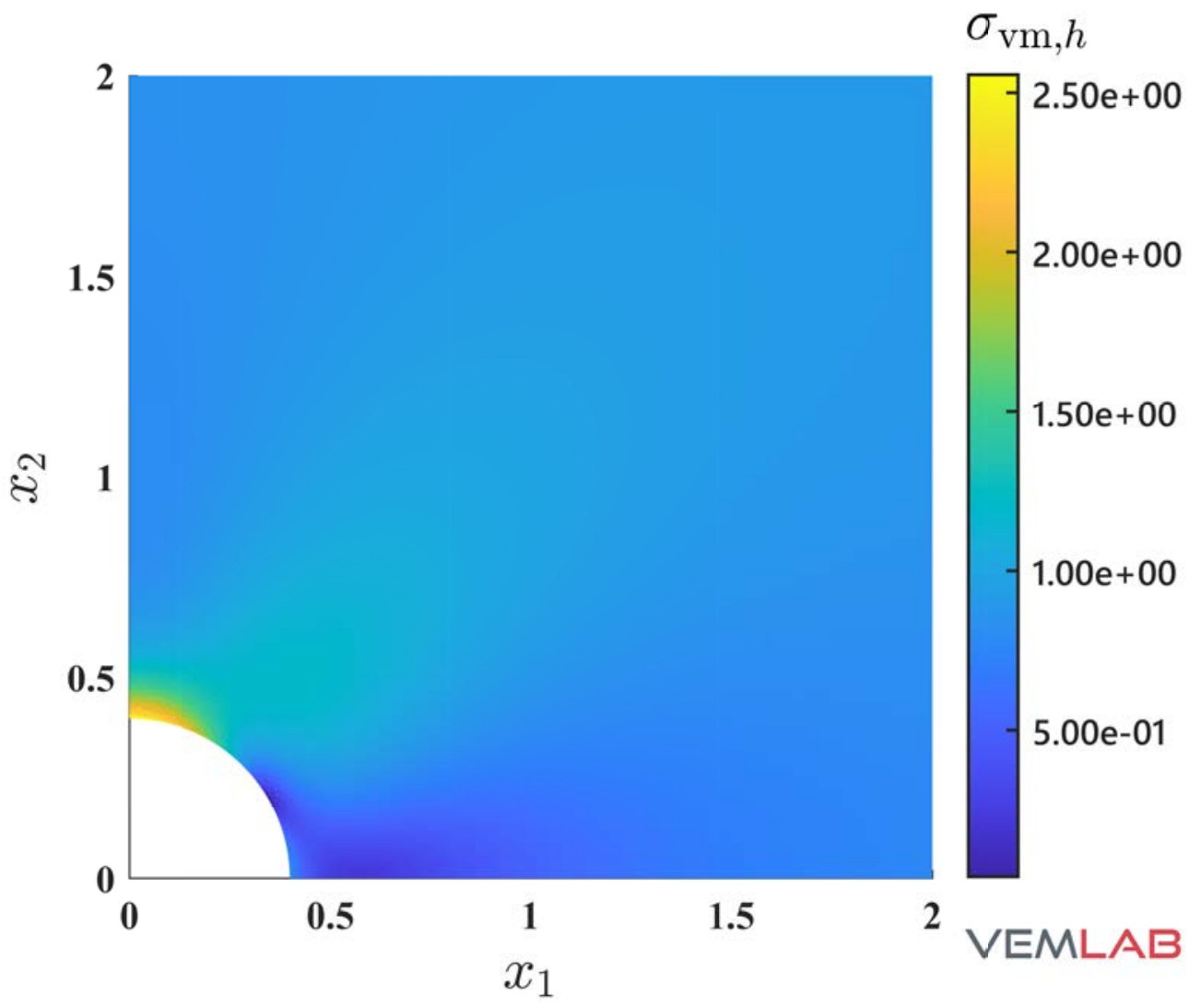}}
	\caption{Infinite plate with a circular hole problem. 
	Material parameters $E_Y=10^3$ \alejandro{psi} and $\nu=0.499999$ for plane strain condition. 
	Plots of the von Mises stress field
	solution ($\sigma_{\textrm{vm},h}$) \alejandro{in psi} for the (a) NVEM ($\tilde{\vm{D}}$ stabilization), 
	(b) NVEM ($\vm{D}_{\mu}$ stabilization) and (c) B-bar VEM approaches.}
	\label{fig:platecontourvm}
\end{figure}

Finally, plots of the displacement, pressure and von Mises stress fields are presented in
Figs. \ref{fig:platecontouru}, \ref{fig:platecontourp} and \ref{fig:platecontourvm}, respectively.
Scatter plots are used for the NVEM as in this approach the field variables are known at the nodes. 
A very good agreement between the NVEM and the B-bar VEM solutions is found in these plots.

\subsection{Dynamic response of a cantilever beam}
\label{sec:dynamicbeam}
\alejandro{
The last numerical example aims to assess the NVEM approach in an elastodynamic
analysis of a cantilever beam subjected to a point force that is applied at the tip of the beam. 
Three cases are considered for the tip load: a constant impact load, a variable impact load 
and a harmonic load. The body force due to the acceleration of gravity is present at all times in the analysis.
The beam geometry and boundary conditions are shown in~\fref{fig:elastobeam_problem}. 
Four meshes of increasing density are considered, which
are shown in~\fref{fig:elastobeam_meshes}. For each of the three cases considered, a reference
solution is obtained with a highly refined mesh of eight-node quadrilateral finite elements (FEM, Q8)
consisting of 62850 degrees of freedom.

\begin{figure}[!bth]
        \centering
        \includegraphics[width=1.0\linewidth]{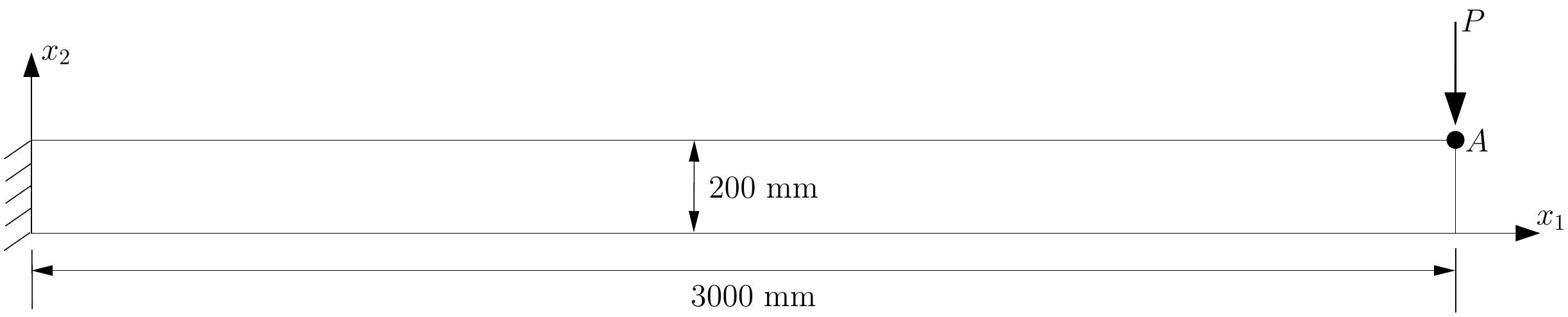}
        \caption{Geometry and boundary conditions for the dynamic response of 
        a cantilever beam problem.}
        \label{fig:elastobeam_problem}
\end{figure}
\begin{figure}[!bth]
        \centering
        \subfigure[] {\label{fig:elastobeam_mesh_a}\includegraphics[width=0.9\linewidth]{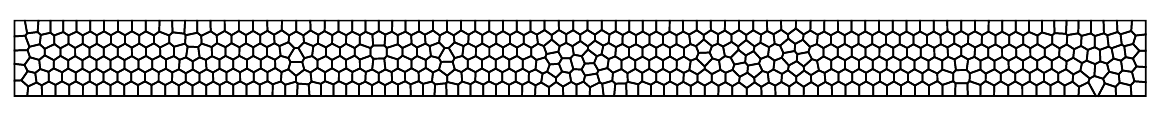}}
        \subfigure[] {\label{fig:elastobeam_mesh_b}\includegraphics[width=0.9\linewidth]{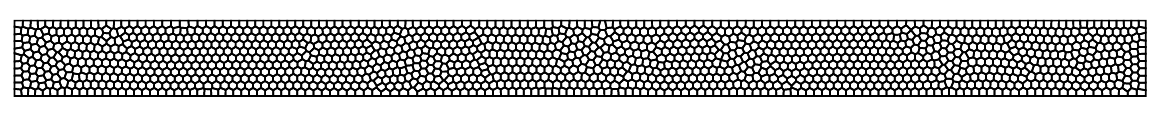}}   
        \subfigure[] {\label{fig:elastobeam_mesh_c}\includegraphics[width=0.9\linewidth]{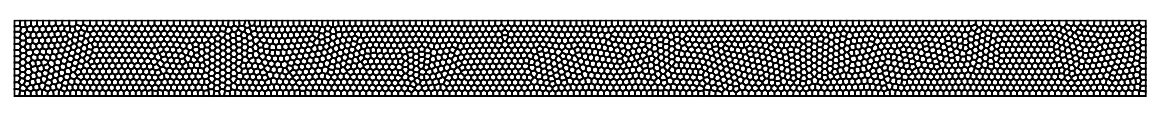}}
        \subfigure[] {\label{fig:elastobeam_mesh_d}\includegraphics[width=0.9\linewidth]{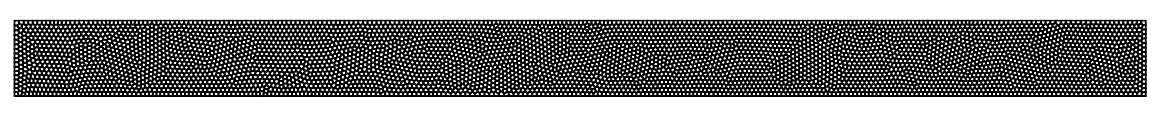}}
        \caption{Meshes for the dynamic response of a cantilever beam. (a) 500 polygonal elements (2004 degrees of freedom), 
                 (b) 1500 polygonal elements (6004 degrees of freedom), (c) 2500 polygonal elements (10002 degrees of freedom) and (d) 5000 polygonal elements (20000 degrees of freedom).}
        \label{fig:elastobeam_meshes}
\end{figure}

The semidiscrete equation of motion of the elastodynamic problem is given by~\cite{hughes:2000}
\begin{equation*}
\vm{M}\ddot{\vm{d}}+\vm{C}\dot{\vm{d}}+\vm{K}\vm{d}=\vm{F},
\end{equation*}
where $\vm{M}$ is the global mass matrix, $\vm{C}$ is the global viscous damping matrix, 
$\vm{K}$ is the global stiffness matrix, $\vm{F}$ is the global vector of external forces, and
$\ddot{\vm{d}}$, $\dot{\vm{d}}$ and $\vm{d}$ are the nodal accelerations, nodal 
velocities and nodal displacements, respectively. In the NVEM, $\vm{K}$ and $\vm{F}$ 
are formed, respectively, by assembling the element matrices $\vm{K}_I$ and $\vm{f}_I$ that were 
developed in the preceding sections. 
The element mass matrix for the NVEM can be constructed following 
the construction of the element mass matrix
for the VEM as given in Ref.~\cite{BeiraodaVeiga-Brezzi-Marini-Russo:2014}, which, 
in the notation used in the present article, is:
\begin{equation*}
(\vm{M}_E)_{ab} = \rho\,\left(\int_E\Pi\phi_a\Pi\phi_b\,\diffx + |E|\left[(\vm{I}-\vm{P})^\transpose (\vm{I}-\vm{P})\right]_{ab}\right).
\end{equation*}
Notice that applying the averaging operator $\pi_I$ to $\Pi\phi_a$, 
that is $\pi_I[\Pi\phi_a]=(\Pi\phi_a)_I$, gives a constant matrix representing the
evaluation of $\Pi\phi_a$ at node $I$. Hence, the NVEM representation of the 
element mass matrix is written as follows:
\begin{equation*}
(\vm{M}_I)_{ab} = \rho\,|I|\Bigl((\Pi\phi_a)_I (\Pi\phi_b)_I + \left[(\vm{I}-\vm{P})_I^\transpose (\vm{I}-\vm{P})_I\right]_{ab}\Bigr).
\end{equation*}
The damping matrix is estimated using the Rayleigh formula~\cite{hughes:2000}:
\begin{equation*}
\vm{C}=q_1\vm{M}+q_2\vm{K},
\end{equation*}
where $q_1$ and $q_2$ are constant input parameters.
The Hilber-Hughes-Taylor $\alpha$-method~\cite{hughes:2000} is used to solve the semidiscrete 
equation of motion, which is an implicit, second-order accurate, unconditionally
stable algorithm with high-frequency numerical dissipation. The $\alpha$-method
leads to the trapezoidal rule for time stepping when the parameter $\alpha$ is set
to zero. 

The following input data is used for all the simulations performed in this
example: $\widetilde{\vm{D}}$ stabilization (results are not affected if $\vm{D}_\mu$
stabilization is chosen instead), plane stress assumption, beam thickness 20 mm, 
Young's modulus $E_\mathrm{Y}$ = 200000 MPa,
Poisson's ratio $\nu=0.3$,
mass density $\rho=7.85\times 10^{-9}$ t/$\mathrm{mm}^3$, acceleration of gravity 
$g=9800$ $\mathrm{mm}/\mathrm{s}^2$, Rayleigh damping parameters $q_1=0.025$ and
$q_2=0.001$, $\alpha$-method parameter $\alpha=-0.1$, time interval in seconds $[0,\,1.6]$,
time step $\Delta t=0.001$ s --- the time step is chosen according to the usual rule that
the period associated with the lowest natural frequency is at least 30 times the time 
step size. Using the Euler-Bernoulli theory to estimate the lowest natural 
frequency gives 18.12 hz; hence, the step size is chosen such that
$\Delta t < (1/18.12)/30 = 0.0018$ s. 

\subsubsection{Constant impact load}
The beam is initially undeformed and at rest. A constant impact load $P=-10000$ N is
applied at time $t=0$ s at the tip of the beam (point $A$) and is assumed to remain 
fixed throughout the time interval of analysis. \fref{fig:response_constant_impact_load} 
presents the response curves at point $A$ for the meshes shown in \fref{fig:elastobeam_meshes}.
As expected, it is observed that with mesh refinement the accuracy of the dynamic response 
increases. In particular, the finer mesh (mesh (d)) matches very well with the 
reference finite element dynamic response.

\begin{figure}[!bth]
        \centering
        \includegraphics[width=0.55\linewidth]{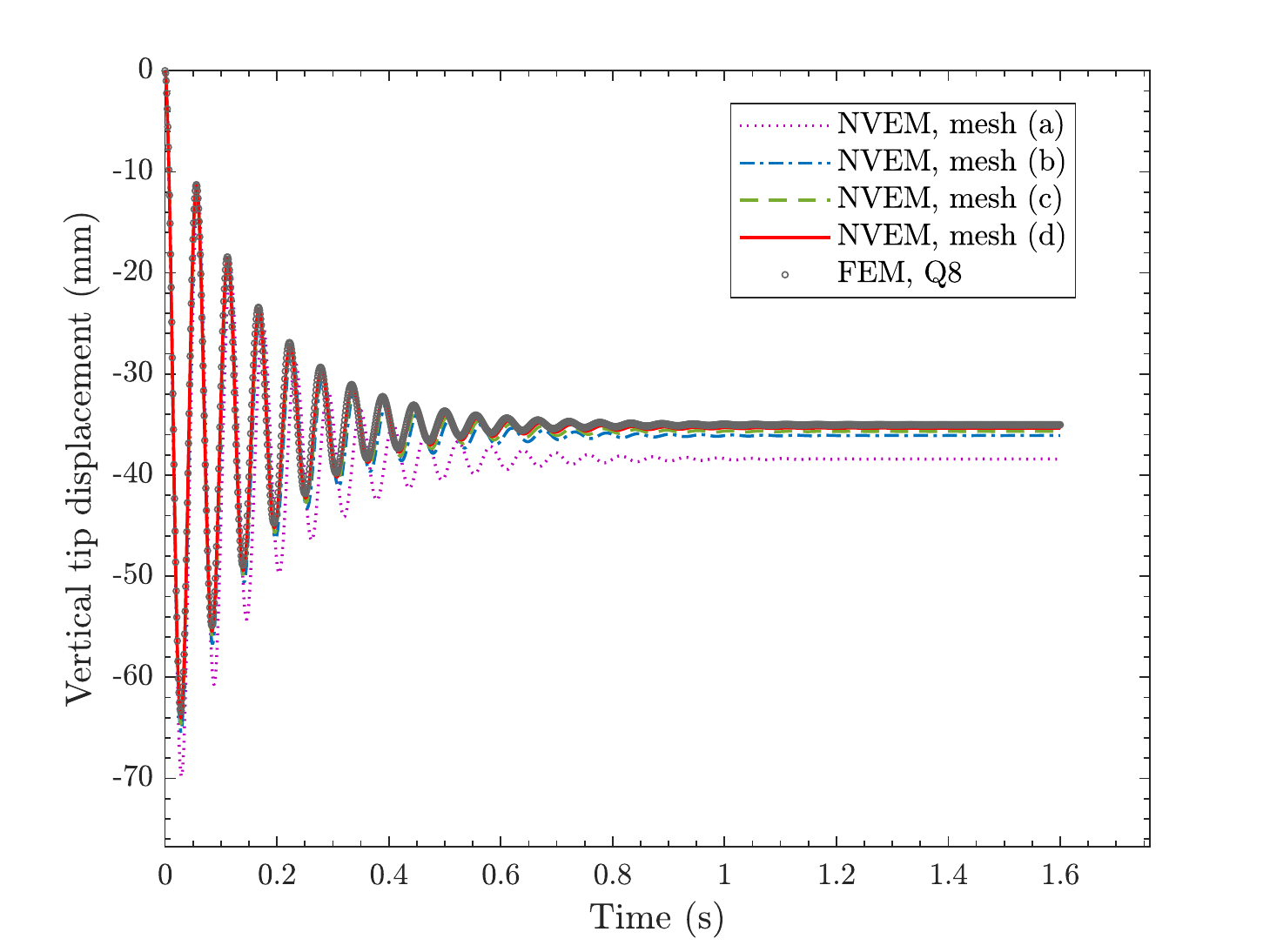}       
        \caption{Dynamic response of the cantilever beam subjected to a constant impact tip load.
        Comparison among the response curves at the tip of beam for the NVEM 
        using four meshes of increasing density and a reference
        response curve obtained with 8-node quadrilateral finite elements (FEM, Q8).}
        \label{fig:response_constant_impact_load}
\end{figure}

\subsubsection{Variable impact load}
The beam is initially undeformed and at rest. The impact load at the tip of 
the beam (point $A$) is $P=f(t)F$, where $F=-10000$ N and $f(t)$ is defined as
\begin{equation*}
f(t)=\left\{ \begin{array}{cl}
         -500t+1, &  0 \leq t \leq 0.002\\
        0, & t > 0.002\end{array} \right.
\end{equation*}
\fref{fig:response_variable_impact_load} presents the response curves at point $A$ 
for the meshes shown in \fref{fig:elastobeam_meshes}. It is observed 
that the accuracy of the dynamic response increases with mesh refinement and
that the finer mesh (mesh (d)) matches very well with the reference finite element dynamic 
response.

\begin{figure}[!bth]
        \centering
        \includegraphics[width=0.55\linewidth]{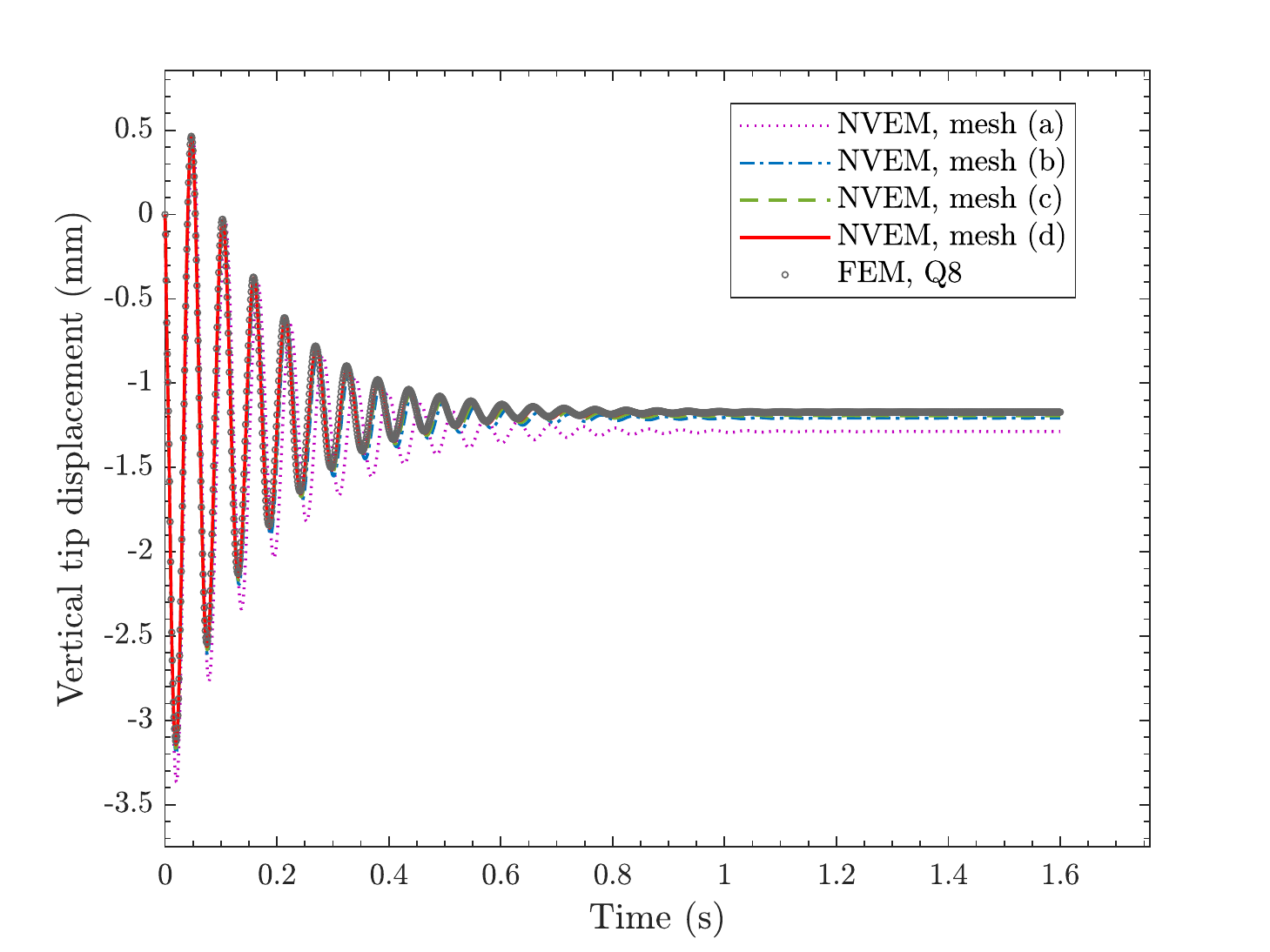}       
        \caption{Dynamic response of the cantilever beam subjected to a variable impact tip load.
        Comparison among the response curves at the tip of beam for the NVEM 
        using four meshes of increasing density and a reference
        response curve obtained with 8-node quadrilateral finite elements (FEM, Q8).}
        \label{fig:response_variable_impact_load}
\end{figure}

\subsubsection{Harmonic load}
The beam is initially undeformed and at rest. 
A harmonic load defined as $P=-10000 \cos{(\omega_f t)}$ N, where $\omega_f=8$ rad/s,
is applied at point $A$. \fref{fig:response_harmonic_load} depicts the response curves at point $A$ 
for the meshes shown in \fref{fig:elastobeam_meshes}. As in the two preceding load cases,
the results show that the accuracy of the dynamic response increases with mesh refinement.
Among the meshes considered, the finer mesh (mesh (d)) provides the best agreement 
with the reference finite element dynamic response.
}

\begin{figure}[!bth]
        \centering
        \includegraphics[width=0.55\linewidth]{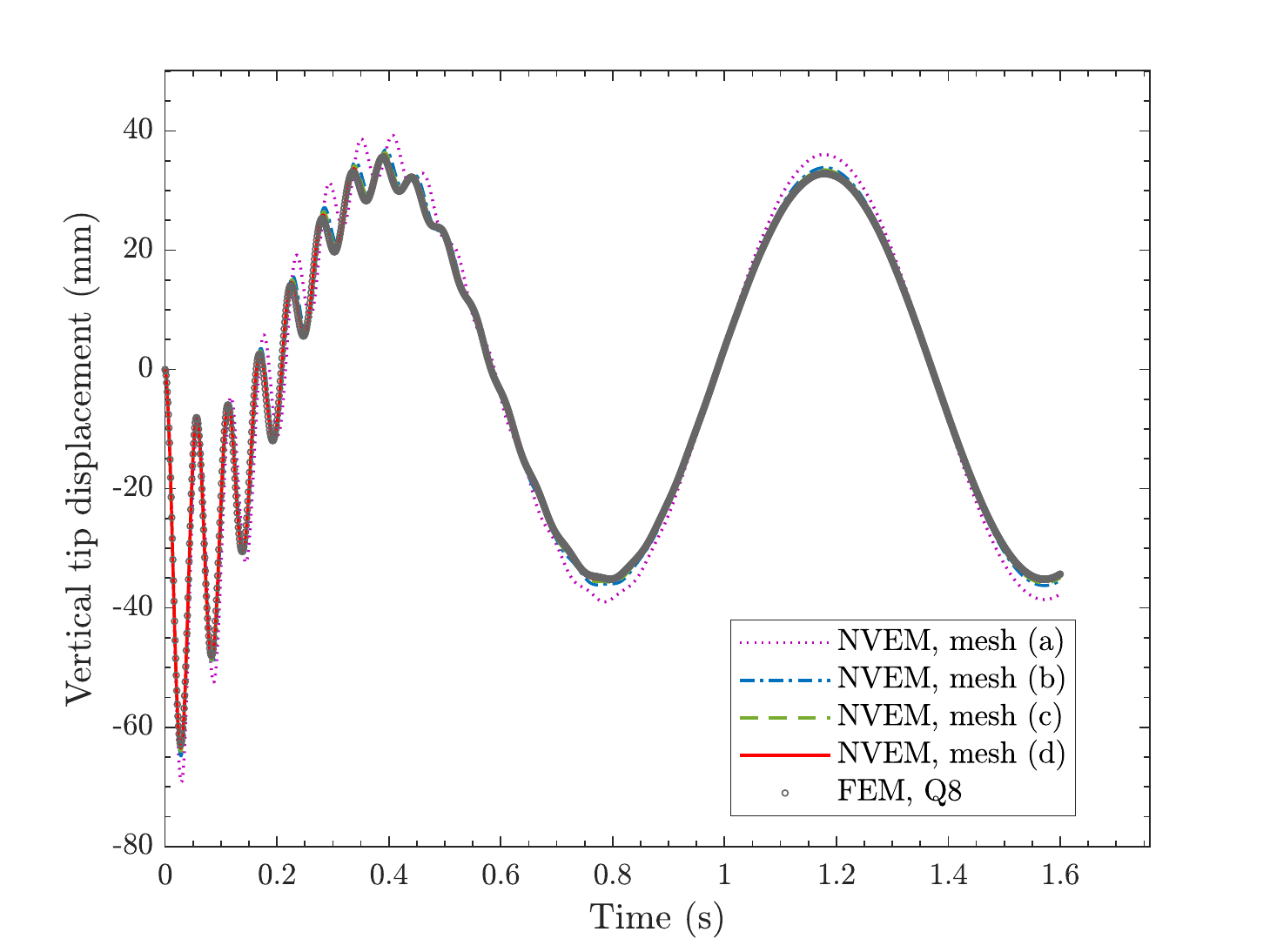}       
        \caption{Dynamic response of the cantilever beam subjected to a harmonic tip load.
        Comparison among the response curves at the tip of beam for the NVEM 
        using four meshes of increasing density and a reference
        response curve obtained with 8-node quadrilateral finite elements (FEM, Q8).}
        \label{fig:response_harmonic_load}
\end{figure}

\section{Concluding remarks}
\label{sec:conclusions}

In this paper, we proposed a combined nodal integration and virtual element method,
wherein the strain is averaged at the nodes from the strain of surrounding virtual elements. 
For the strain averaging procedure, a nodal averaging operator was constructed using a generalization to
virtual elements of the node-based uniform strain approach for 
finite elements~\cite{Dohrmann-Heinstein-Jung-Key-Witkowski:2000}.
This gave rise to the node-based uniform strain virtual element method (NVEM). 
In line with nodal integration techniques, two stabilizations were proposed for 
the NVEM using modified constitutive matrices. The main distinction 
between the NVEM and existing VEM approaches for compressible and 
nearly incompressible elasticity, such as the B-Bar VEM~\cite{Park-Chi-Paulino:2021}, 
is that the stresses and strains in the NVEM are nodal variables just like 
displacements, whereas in the existing methods the stresses and strains are element variables.

Several examples in compressible and nearly incompressible elasticity were conducted
using regular, distorted and random discretizations in two dimensions. The NVEM showed
accurate, convergent and stable solutions in the $L^2$ norm and $H^1$ seminorm of the displacement
error in all the numerical tests. \alejandro{An elastodynamic numerical experiment was also conducted
for the NVEM. The dynamic response of the NVEM was found to be in good agreement
with reference solutions obtained with 8-node quadrilateral finite elements.}
Even though the NVEM and the B-Bar VEM 
behave similar in terms of accuracy and convergence, we remark that all 
the field variables in the NVEM are related to the nodes, 
which would obviously ease the tracking of the state and history-dependent variables 
in the nonlinear regime. 
\alejandro{Also, as in meshfree methods and 
with the state and history-dependent variables stored at the nodes, the NVEM
would be suitable for large deformations settings with remeshing if
a remeshing strategy were designed such that the connectivity
is revaluated while retaining the current set of nodes. This would
obviate the remapping of state and history-dependent variables 
between the old and new meshes making the method very attractive.}
The aforementioned features provide
room for further development of the NVEM. In the short term perspective, 
our current work on the NVEM is focused on its extension to 
materially-nonlinear-only formulation. In the medium 
term, we plan to develop its extension to \alejandro{large deformations with remeshing}.

\section*{Acknowledgements}
This work was performed under the auspices of the Chilean National Fund for Scientific 
and Technological Development (FONDECYT) through grants 
ANID FONDECYT No. 1221325 (A.O-B and R.S-V) and ANID FONDECYT No. 1211484 (N.H-K).
S.S-F gratefully acknowledges the research support of the 
Chilean National Agency of Research and Development (ANID)
through grant ANID doctoral scholarship No. 21202379.


\clearpage

\appendix
\section{List of main symbols}
\label{appA}

\begin{flushleft}
\begin{tabular}{ l l }

	$\vm{b}$  & body force vector\\ 
	$\widehat{\vm{b}}$  & element average of the body force vector\\ 	
	$\widehat{\vm{b}}_I$  & nodal average of $\widehat{\vm{b}}$\\ 		
	$\vm{D}$   & constitutive matrix for the linear elastic material \\
	$\widetilde{\vm{D}}$   & modified constitutive matrix for the linear elastic material \\
	$\vm{D}_\mu$   & modified constitutive matrix for the linear elastic material \\		
 	$E$ & element or domain of the element \\
	$e$ & edge of an element \\ 
  $|E|$  &  area of the element  \\
  $|e|$  & length of an element edge  \\	
  $E_\textrm{Y}$ & Young's modulus of the linear elastic material \\
  $\vm{f}_E^b$ & VEM element force vector associated with $\vm{b}$ \\
  $\vm{f}_I^b$ & NVEM nodal force vector associated with $\vm{b}$ \\  
  $\vm{f}_e^t$ & VEM element force vector associated with $\vm{t}_N$ \\
  $\vm{f}_I^t$ & NVEM nodal force vector associated with $\vm{t}_N$ \\  
 	$I$ & a node of the mesh \\  
 	$|I|$ & representative area of node $I$ \\  	
 	$\vm{K}_E$ & VEM element stiffness matrix \\
 	$\vm{K}_I$ & NVEM nodal stiffness matrix \\ 	
	$\vm{K}_E^\cons$ & consistency part of the VEM element stiffness matrix \\ 
	$\vm{K}_I^\cons$ & consistency part of the NVEM nodal stiffness matrix \\ 
        $\vm{K}_E^\stab$  &  stability part of the VEM element stiffness matrix \\
        $\vm{K}_I^\stab$  &  stability part of the NVEM nodal stiffness matrix \\    
        $\vm{M}_I$ & NVEM nodal mass matrix \\            
        $\vm{n}$  &  unit outward normal to the element boundary \\     
        $\vm{n}_a$  &  unit outward normal to the $a$-th edge of the element \\                 
        $\vm{n}_\Gamma$  & unit outward normal to the domain boundary \\          			
\end{tabular}
\end{flushleft}

\begin{flushleft}
\begin{tabular}{ l l }      
        $N_E^V$  & number of edges/nodes of an element \\                       
        $N_e^V$  & number of nodes of an element edge \\        
        $\vm{t}_N$ & Neumann boundary condition \\       
        $\widehat{\vm{t}}_N$  & edge average of the Neumann boundary condition\\        
        $\widehat{\vm{t}}_{N,I}$  & nodal average of $\widehat{\vm{t}}_N$\\     
        $\vm{u}$ &  displacement vector  \\     
        $\vm{u}_D$ &  Dirichlet boundary condition  \\          
        $\vm{u}_h$ &  discrete trial displacement vector  \\    
        $\vm{u}_a$ &  nodal displacement vector  \\
	$\bar{\vm{u}}$ & mean of the values that $\vm{u}$ takes over the vertices of an element\\
	$\vm{v}$ &  test displacement vector \\	
	$\vm{v}_h$ &  discrete test displacement vector \\	
	$\vm{x}$ & position vector\\	
	$\vm{x}_a$ & nodal coordinates vector\\	
	$\bar{\vm{x}}$ & mean of the values that $\vm{x}$ takes over the vertices of an element\\	
  $\bsym{\varepsilon}$  & small strain tensor\\
  $\widehat{\bsym{\varepsilon}}$ & element average of the small strain tensor \\
  $\bsym{\varepsilon}_I$  & nodal strain tensor\\ 
  $\partial E$  & boundary of the element  \\
  $\Gamma$  & boundary of the body \\
  $\Gamma_D$  & Dirichlet boundary \\           
        $\Gamma_N$ &  Neumann boundary \\       
        $\lambda$ &  Lam\'e's first parameter\\
        $\mu$ &  Lam\'e's second parameter (shear modulus of the linear elastic material)\\
        $\nu$ & Poisson's ratio of the linear elastic material \\
  $\Omega$ & domain of the  body\\
  $\bsym{\omega}$  & small rotation tensor\\
  $\widehat{\bsym{\omega}}$ & element average of the small rotation tensor \\
  $\bm{\sigma}$  & Cauchy stress tensor \\
\end{tabular}
\end{flushleft}

\bibliographystyle{elsarticle-num}

\end{document}